\documentclass[12pt]{amsart}
\usepackage{amsmath,amsthm,amssymb}
\usepackage{lscape}
\usepackage[all]{xy}
\usepackage{pifont}
\usepackage{bbding}
\usepackage{url}

\newtheorem{thm}{Theorem}[section]
\newtheorem{lem}[thm]{Lemma}
\newtheorem{dfn}[thm]{Definition}
\newtheorem{prop}[thm]{Proposition}
\newtheorem{cor}[thm]{Corollary}

\newtheorem{rem}[thm]{Remark}
\newenvironment{prf}{{\bf Proof.}}{\hfill $\diamond$}

\begin{document}

\def\id{ {\rm id } }
\def\dim{ {\rm dim } }
\def\Dim{ {\rm Dim } }
\def\Aut{ {\rm Aut } }
\def\Hom{ {\rm Hom} }
\def\Irr{ {\rm Irr} }
\def\Ind{ {\rm Ind} }
\def\End{ {\rm End} }
\def\Res{ {\rm Res} }
\def\Spin{ {\rm Spin} }
\def\khar{ {\rm char } }
\def\diag{ {\rm diag } }
\def\Stab{ {\rm Stab } }
\def\da{ \!\!\downarrow }
\def\NI{ \noindent }
\def\ol{ \overline }
\def\ua{ \!\!\uparrow }
\def\wh{ \widehat }
\def\wt{ \widetilde }
\def\mC{ {\mathbb C} }
\def\mF{ {\mathbb F} }
\def\bN{ {\bf N} }
\def\bff{ {\bf f} }
\def\mZ{ {\mathbb Z} }
\def\GL{ \text{\rm GL} }
\def\GU{ \text{\rm GU} }
\def\GO{ \text{\rm GO} }
\def\SO{ \text{\rm SO} }
\def\Sp{ \text{\rm Sp} }
\def\CSp{ \text{\rm CSp} }
\def\PCSp{ \text{\rm PCSp} }
\def\PSp{ \text{\rm PSp} }
\def\SL{ \text{\rm SL} }
\def\PSL{ \text{\rm PSL} }
\def\PGL{ \text{\rm PGL} }
\def\PGU{ \text{\rm PGU} }
\def\SU{ \text{\rm SU} }
\def\PSU{ \text{\rm PSU} }
\def\SO{ \text{\rm SO} }
\def\CO{ \text{\rm CO} }
\def\PSO{ \text{\rm PSO} }
\def\CSO{ \text{\rm CSO} }
\def\CGO{ \text{\rm CGO} }
\def\PCSO{ \text{\rm PCSO} }

%\frontmatter

\title{Imprimitive irreducible modules for finite quasisimple groups, II}

%    Remove any unused author tags.

%    author one information
\author{Gerhard Hiss}
\address{Lehrstuhl D f{\"u}r Mathematik, RWTH Aachen University,
52056 Aach\-en, Germany}
%\curraddr{}
\email{gerhard.hiss@math.rwth-aachen.de}

\author{Kay Magaard}
\address{School of Mathematics, University of Birmingham, Edgbaston,
Birmingham B15 2TT, UK}
%\curraddr{}
\email{k.magaard@bham.ac.uk}
%\thanks{}

%    \date is required; it is the date received by the editor.
%\date{02.08.2013}
%\date{06.01.2014}
%\date{17.10.2014}
%\date{18.11.2014}
%\date{12.02.2015}
%\date{23.09.2015}
%\date{25.09.2015}
%\date{30.09.2015}
%\date{12.10.2015}
%\date{02.11.2015}
%\date{20.11.2015}
%\date{09.12.2015}
%\date{16.02.2016}
%\date{18.02.2016}
%\date{23.03.2016}
%\date{30.03.2016}
%\date{16.05.2016}
%\date{28.06.2016}
%\date{30.09.2016}
%\date{17.10.2016}
%\date{21.10.2016}
%\date{28.10.2016}
%\date{31.10.2016}
\date{20.11.2016}



\begin{abstract}
This work completes the classification of the imprimitive irreducible
modules, over algebraically closed fields of characteristic~$0$, of 
the finite quasisimple groups.
\end{abstract}

\maketitle

\section{Introduction}

This is a continuation of our work in~\cite{HiHuMa}, where we began the
classification of the irreducible imprimitive modules of the finite
quasisimple groups over algebraically closed fields. We completed this
program for the quasisimple covering groups of the sporadic simple groups
and for the Tits simple group, as well as for the covering groups of the
finite groups of Lie type with an exceptional Schur multiplier or with
more than one defining characteristic.
For the covering groups of the alternating groups, the desired classification
is complete over fields of characteristic~$0$ by the work of Djokovi\'c and 
Malzan \cite{DM2}, respectively Nett and Noeske~\cite{NeNoe}. In this paper we
finish the classification for all quasisimple groups~$G$ and fields of
characteristic~$0$. 

Let~$G$ be a finite group. Write $\Irr(G)$ for the set of irreducible characters 
of~$G$ over some algebraically closed field of characteristic~$0$. An element 
$\chi \in \Irr(G)$ is \textit{imprimitive}, if it is induced from a proper 
subgroup of~$G$. Now let~$G$ be a finite group with a split $BN$-pair, for 
example a finite reductive group. We say that $\chi \in \Irr(G)$ is 
\textit{Harish-Chandra imprimitive}, if there is a proper Levi 
subgroup $L \leq G$ (in the sense of groups with a split $BN$-pair) and 
$\vartheta \in \Irr(L)$ such that $\chi = 
R_L^G( \vartheta )$. Here,~$R_L^G$ denotes Harish-Chandra induction from~$L$ 
to~$G$. Next, let~$G$ be a quasisimple finite reductive group. 
Suppose that~$G/Z(G)$ does not have 
an exceptional Schur multiplier and that~$G$ has a unique defining 
characteristic. Then~\cite[Theorem~$6.1$]{HiHuMa} states that $\chi \in \Irr(G)$ 
is imprimitive if and only if it is Harish-Chandra imprimitive. This leads to
the program to classify the Harish-Chandra imprimitive elements of $\Irr(G)$
for finite reductive groups~$G$.
In \cite[Theorem~$7.3$, Theorem~$8.4$]{HiHuMa}, such a classification has been 
achieved for finite reductive groups~$G$ arising from algebraic groups with a 
connected center (but~$G$ not necessarily quasisimple). Also, 
by \cite[Propositions $10.2$--$10.4$]{HiHuMa}, the classification is complete
for the Suzuki and Ree groups, so that we will not consider these groups 
in the following.

To formulate our main result, let us introduce some more notation. 
Let~$\mathbf{G}$ denote a connected reductive algebraic group defined over the 
algebraic closure of a finite field, and let~$F$ denote a Frobenius morphism 
of~$\mathbf{G}$. Put
$G := \mathbf{G}^F$, the finite group of $F$-fixed points on~$\mathbf{G}$.
We also let $(\mathbf{G}^*,F)$ denote a pair of a reductive group and a 
Frobenius morphism in duality to~$(\mathbf{G},F)$. Let~$s \in G^*$ be semisimple 
and write $\mathcal{E}(G,[s]) \subseteq \Irr(G)$ for the Lusztig series defined 
by the $G^*$-conjugacy class~$[s]$ of~$s$. By Lusztig's generalized Jordan 
decomposition of characters, there is an equivalence relation on 
$\mathcal{E}(G,[s])$ and a bijection between the equivalence classes on 
$\mathcal{E}(G,[s])$ and the $C_{\mathbf{G}^*}( s )^F$-orbits on 
$\mathcal{E}( C^\circ_{\mathbf{G}^*}( s )^F, [1] )$. (This bijection is, 
in general, not unique, but our results are true for any such bijection.)
We write $[\chi]$
for the equivalence class of $\chi \in \mathcal{E}(G,[s])$, and 
$[\lambda]$ for the $C_{\mathbf{G}^*}( s )^F$-orbit of $\lambda \in 
\mathcal{E}( C^\circ_{\mathbf{G}^*}( s )^F, [1] )$. For such a~$\lambda$ we let
$C_{\mathbf{G}^*}( s )^F_\lambda$ denote its stabilizer 
in~$C_{\mathbf{G}^*}( s )^F$.
We can now formulate one of the main results of our paper.

\begin{thm}
\label{MainResult}
Let~$\chi \in \Irr( G )$. Suppose that $\chi \in \mathcal{E}( G, [s] )$
for a semisimple element $s \in G^*$, and let $\lambda \in 
\mathcal{E}( C_{\mathbf{G}^*}^\circ( s )^F, [ 1 ] )$ such that $[\chi]$
corresponds to $[\lambda]$ under Lusztig's generalized Jordan decomposition of
characters. Then the following assertions hold.

{\rm (a)} If
\begin{equation}
\label{InclusionCondition}
C_{\mathbf{G}^*}( s )^F_\lambda\,\,C^\circ_{\mathbf{G}^*}( s ) \leq \mathbf{L}^*,
\end{equation}
for some split $F$-stable Levi subgroup~$\mathbf{L}^*$ of~$\mathbf{G}^*$, 
then~$\chi$ is Harish-Chandra induced from $L = \mathbf{L}^F$, where~$\mathbf{L}$
is an $F$-stable Levi subgroup of~$\mathbf{G}$ dual to~$\mathbf{L}^*$ (hence split).

{\rm (b)} Suppose that~$\mathbf{G}$ is simple and simply connected. If~$\chi$ is 
Harish-Chandra imprimitive, there is a proper split $F$-stable Levi subgroup 
$\mathbf{L}^*$ of~$\mathbf{G}^*$ such that 
Condition~{\rm (\ref{InclusionCondition})} is satisfied.
\end{thm}

\noindent
If, in the notation above, $C_{\mathbf{G}^*}( s )$ is connected, then 
$C_{\mathbf{G}^*}( s )^F_\lambda\,\,C^\circ_{\mathbf{G}^*}( s ) = 
C_{\mathbf{G}^*}( s )$. Thus Theorem~\ref{MainResult} generalizes parts of 
\cite[Theorem~$7.3$, Theorem~$8.4$]{HiHuMa}.

Let us now comment on the further results and the content of the individual
sections of our paper. Section~$2$ contains preliminaries and introduces some
notation. In Section~$3$ we prove Theorem~\ref{MainResult}(a) completely, as 
well as Theorem~\ref{MainResult}(b) under some restrictions. To remove these 
restricitions we have to 
resort to a case by case analysis, investigating all the possible 
groups~$\mathbf{G}$. This is achieved in Section~$5$, where we explicitly decide 
the Harish-Chandra imprimitivity of an element $\chi \in \mathcal{E}( G, [s] )$ 
in terms of the label of the unipotent character $\lambda \in 
\mathcal{E}( C_{\mathbf{G}^*}^\circ( s )^F, [ 1 ] )$, where~$[\chi]$ 
and~$[\lambda]$  correspond under Lusztig's generalized Jordan decomposition of
characters. For the classical groups~$\mathbf{G}$ this requires detailed
knowledge on the semisimple elements $s \in G^*$ with $C_{\mathbf{G}^*}(s)$
non-connected, as well as on the action of $C_{\mathbf{G}^*}(s)^F$ on the
unipotent characters of $C^\circ_{\mathbf{G}^*}(s)^F$. This information is
collected in Section~$4$ of our paper, and is given in terms of the natural
representation of a classical group~$\tilde{\mathbf{G}}^*$, 
respectively~$\check{\mathbf{G}}^*$, having~$\mathbf{G}^*$ as epimorphic image.
More precisely, we choose an element $\tilde{s} \in \tilde{G}^*$ mapping 
onto~$s$, and give the conditions for the Harish-Chandra imprimitivity of~$\chi$
in terms of the minimal polynomial of~$\tilde{s}$ in the natural
representation of~$\tilde{G}^*$ (and similarly for $\check{G}^*$). One could
as well work entirely within a maximally split torus of~$\mathbf{G}^*$ and the
action of the Weyl group on this torus, but this again is best achieved with the
help of a natural representation of a group surjecting onto~$\mathbf{G}^*$. The 
corresponding information for
the exceptional groups is taken from the tables~\cite{LL} computed by Frank 
L{\"u}beck. Although not visible any more in the final version of our paper,
we owe very much to explicit computations with substantial examples using 
GAP~\cite{GAP4}.

\section{Notation and preliminaries}

\subsection{Some general notation}
\label{SomeGeneralNotation}

The $n \times n$ identity matrix is denoted by~$I_n$, and by~$J_n$ we denote
the $n \times n$ anti-diagonal matrix. If $n = 2m$, we put
$$\tilde{J}_{n} := \left[ \begin{array}{cc} 0 & J_m \\ -J_m & 0 \end{array}\right].$$
In each case we omit the index~$n$ if
this is clear from the context. The transposed of a matrix~$a$ is denoted by~$a^T$.

By a character of a finite group~$H$ we mean a complex character and
$\Irr( H )$ denotes the set of irreducible characters of~$H$.

\subsection{Finite reductive groups}
\label{FiniteReductiveGroups}

Throughout this paper, we let~$\mathbb{F}$ denote an algebraic closure of 
the finite field with~$p$ elements, and we let~$q$ be a power of~$p$. 
Let~$\mathbf{G}$ be a connected reductive algebraic group defined over
$\mathbb{F}$, and let~$F$ be a Frobenius morphism of~$\mathbf{G}$ arising 
from an $\mathbb{F}_q$-structure on~$\mathbf{G}$. For any $F$-stable 
subgroup~$\mathbf{H}$ of~$\mathbf{G}$ we put $\mathbf{H}^F := \{ h \in
\mathbf{H} \mid F(h) = h \}$ for the finite group of $F$-fixed points 
in~$\mathbf{H}$. We call $\mathbf{G}^F$ a \textit{finite reductive group}.
We adopt the typographical convention to write $H := \mathbf{H}^F$ for
$F$-stable subgroups $\mathbf{H} \leq \mathbf{G}$ denoted by a single boldface 
letter. For reasons of clarity, we prefer the notation with exponent~$F$ for 
subgroups such as centralizers or normalizers. For example, we usually
write $C_{\mathbf{G}}(s)^F$ rather than $C_G(s)$, if~$s$ is an $F$-stable 
element of~$\mathbf{G}$.
The connected component of a closed subgroup $\mathbf{H}$ of~$\mathbf{G}$
is denoted by~$\mathbf{H}^\circ$. In order to avoid double exponents, we
put $Z^\circ( \mathbf{G} ) := Z( \mathbf{G} )^\circ$ and 
$C^\circ_{\mathbf{G}}(s) := C_{\mathbf{G}}(s)^\circ$, for $s \in \mathbf{G}$.

An $F$-stable Levi subgroup~$\mathbf{L}$ of~$\mathbf{G}$ is called \textit{split}, 
if~$\mathbf{L}$ is the Levi complement of an $F$-stable parabolic subgroup 
of~$\mathbf{G}$. This is the case, if and only if~$\mathbf{L}$ is the 
centralizer of a split torus of~$\mathbf{G}$. The 
finite reductive group~$G$ is, in particular, a finite group with a split 
$BN$-pair of characteristic~$p$. The Levi subgroups~$L$ of~$G$, in the sense 
of a group with a split $BN$-pair, are of the form $L = \mathbf{L}^F$ for 
$F$-stable split Levi subgroups~$\mathbf{L}$ of~$\mathbf{G}$.

Let~$L$ be a Levi subgroup of~$G$ and $\vartheta \in \Irr(L)$. Then 
$R_L^G( \vartheta )$ is the character of~$G$ obtained by Harish-Chandra 
inducing~$\vartheta$ from~$L$ to~$G$.

\subsection{The Weyl group}
\label{TheWeylGroup}

Fix an $F$-stable maximally split torus~$\mathbf{T}$ of~$\mathbf{G}$ and let
$W := N_{\mathbf{G}}( \mathbf{T} )/\mathbf{T}$ denote the Weyl group
of~$\mathbf{G}$ with respect to~$\mathbf{T}$. For every $w \in W$ we choose
an inverse image~$\dot{w}$ of~$w$ under the canonical epimorphism
$N_{\mathbf{G}}( \mathbf{T} ) \rightarrow W$. For $t \in \mathbf{T}$ and
$w \in W$ we write $t^w := \dot{w}^{-1} t \dot{w}$. This is independent
of the chosen representative~$\dot{w}$, and defines what we call the 
\textit{ conjugation action} of~$W$ on~$\mathbf{T}$.

\subsection{Semisimple elements and centralizers}
\label{SemisimpleElementsAndCentralizers}

An element $t \in \mathbf{T}$ is conjugate in~$\mathbf{G}$ to an $F$-stable
element, if and only if $t$ is conjugate to $F(t)$ by an element of~$W$. Indeed,
let $s \in \mathbf{G}^F$ be semisimple. Then~$s$ lies in an $F$-stable maximal
torus of~$\mathbf{G}$. In particular, there is $w \in W$ such that~$s$ is
conjugate to an element~$t$ of
$$\mathbf{T}^{Fw} := \{ t \in \mathbf{T} \mid F( t )^w = t \} = 
\{ t \in \mathbf{T} \mid \dot{w}^{-1} F( t ) \dot{w} = t \}.$$
In turn, $C_{\mathbf{G}}( s )^F$ is conjugate in~$\mathbf{G}$ to
$$C_{\mathbf{G}}( t )^{F\dot{w}} := \{ x \in C_{\mathbf{G}}( t ) \mid
\dot{w}^{-1} F( x ) \dot{w} = x \}.$$
Conversely, an element of $\mathbf{T}^{Fw}$ for some $w \in W$ is conjugate
in~$\mathbf{G}$ to an $F$-stable element.

\subsection{The component group}
\label{TheComponentGroup}

For $s \in \mathbf{G}$ semisimple, we put 
$$A_{\mathbf{G}}( s ) := C_{\mathbf{G}}( s )/C^\circ_{\mathbf{G}}( s ).$$ 
Clearly, $A_{\mathbf{G}}( s )$ is~$F$-stable, if~$s$ is $F$-stable. In this 
case we have $A_{\mathbf{G}}( s )^F \cong 
C_{\mathbf{G}}(s)^F/C^\circ_{\mathbf{G}}(s)^F$ as
$C^\circ_{\mathbf{G}}( s )$ is connected. By the Lang-Steinberg theorem,
there is a bijection between the $F$-conjugacy classes in
$A_{\mathbf{G}}( s )^F$ and the $G$-conjugacy classes in
the set of $F$-stable elements in the $\mathbf{G}$-conjugacy class of~$s$
(see \cite[(3.25)]{DiMi2}). If~$\mathbf{L}$ is a Levi subgroup of~$\mathbf{G}$
containing~$s$, then $C^\circ_{\mathbf{G}}(s) \cap \mathbf{L} = 
C^\circ_{\mathbf{L}}( s )$. In particular, 
the inclusion $\mathbf{L} \rightarrow \mathbf{G}$ induces an 
injective group homomorphism $A_\mathbf{L}(s) \rightarrow A_\mathbf{G}(s)$ 
(see \cite[8.B]{CeBo2}).

\subsection{Characters and Lusztig series}
Let $\mathbf{G}^*$ be a reductive group dual to~$\mathbf{G}$, and denote by~$F$
the Frobenius morphism of~$\mathbf{G}^*$ in duality with the Frobenius morphism
on~$\mathbf{G}$. For semisimple $s \in G^*$ we denote by~$[s]$ the 
$G^*$-conjugacy of~$s$ and by $\mathcal{E}( G, [s] ) \subseteq \Irr(G)$ the
rational Lusztig series of characters defined by~$[s]$ 
(see \cite[Definition~$8.23$]{CaEn}).

\subsection{Direct products}
For later use, we record an easy result for direct products of groups with 
a $BN$-pair.
\begin{lem}
\label{DirectProductHCInduction}
Let $G_0$ be a finite group with a split $BN$-pair.
Put $G_i := G_0$, $1 \leq i \leq n$, and 
$$G := G_1 \times G_2 \times \cdots \times G_n.$$
Then the symmetric group~$S_n$ on $n$ letters acts on~$G$ (by permuting the 
factors~$G_i$), hence~$S_n$ acts on $\Irr(G)$.

Let $\alpha_i \in \Irr(G_i)$, $1 \leq i \leq n$, and put
$$\alpha := \alpha_1 \otimes \alpha_2 \otimes \cdots \otimes \alpha_n.$$
For $1 \leq i \leq n$, let~$L_i$ denote a Levi subgoup of~$G_i$, and put 
$$L := L_1 \times L_2 \times \cdots \times L_n.$$
Let $\beta_i \in \Irr(L_i)$, $1 \leq i \leq n$,  and put 
$$\beta := \beta_1 \otimes \beta_2 \otimes \cdots \otimes \beta_n.$$
Suppose that $R_L^G( \beta )$ is multiplicity free and that~$\alpha$ 
is an irreducible constituent of $R_L^G( \beta )$. Finally, suppose that 
there is some $1 \leq j \leq n$ such that $R_{L_j}^{G_j}( \beta_j )$ is 
not irreducible.

Then there is an irreducible constituent of $R_L^G( \beta )$ which does not
lie in the $S_n$-orbit of~$\alpha$.
\end{lem}
\begin{prf}
We have
\begin{equation}
\label{HCInductionDirectProduct}
R_L^G( \beta ) = R_{L_1}^{G_1}( \beta_1 ) \otimes R_{L_2}^{G_2}( \beta_2 ) 
\otimes \cdots \otimes R_{L_n}^{G_n}( \beta_n ).
\end{equation}
As~$\alpha$ is an irreducible constituent of~$R_L^G( \beta )$, each~$\alpha_i$
is an irreducible constituent of $R_{L_i}^{G_i}( \beta_i )$ for 
$1 \leq i \leq n$. 
As $R_L^G( \beta )$ is multiplicity free, the same is true for each
$R_{L_i}^{G_i}( \beta_i )$, $1 \leq i \leq n$. By renumbering, we may assume 
that~$j = 1$. Let $\varphi \neq \psi$ be two irreducible constituents of
$R_{L_1}^{G_1}( \beta_1 )$. By Equation~(\ref{HCInductionDirectProduct}) and 
the subsequent remark, $\varphi \otimes \alpha_2 \otimes \cdots \otimes 
\alpha_n$ and $\psi \otimes \alpha_2 \otimes \cdots \otimes \alpha_n$ are 
distinct irreducible constituents of $R_L^G( \beta )$. At most one of these 
lies in the $S_n$-orbit of~$\alpha$.
\end{prf}

\section{Decent to commutator subgroups}\label{CommutatorSubgroups}
 
In this section we generalize the results of \cite[Chapters 8, 9]{HiHuMa}
to quasisimple groups. (Recall that in \cite[Chapters 8, 9]{HiHuMa} we
have assumed that our groups arise from algebraic groups with connected center.)
As in the proof of \cite[Theorem~8.4]{HiHuMa}, our approach uses regular 
embeddings, a standard technique introduced by Deligne and Lusztig 
in \cite[Corollary 5.18]{DeliLu}.

\subsection{Regular embeddings}
\label{RegularEmbeddings}
Let~$\mathbf{G}$ and~$F$ be as in Subsection~\ref{FiniteReductiveGroups}.
Then there is a connected reductive 
group $\tilde{\mathbf{G}}$, and a Frobenius morphism of $\tilde{\mathbf{G}}$, 
also denoted by~$F$, such that the following conditions are satisfied: The 
center of $\tilde{\mathbf{G}}$ is connected, $\mathbf{G}$ is an $F$-stable 
closed subgroup of $\tilde{\mathbf{G}}$ containing the derived subgroup 
$[\tilde{\mathbf{G}}, \tilde{\mathbf{G}}]$ of $\tilde{\mathbf{G}}$, and the
restriction of~$F$ from~$\tilde{\mathbf{G}}$ to~$\mathbf{G}$ is the original 
Frobenius morphism~$F$ on~$\mathbf{G}$. Let us 
denote by $i : \mathbf{G} \rightarrow \tilde{\mathbf{G}}$ the embedding
of~$\mathbf{G}$ into~$\tilde{\mathbf{G}}$. This induces a surjective morphism 
$i^* : \tilde{\mathbf{G}}^* \rightarrow \mathbf{G}^*$ of dual groups, 
compatible with~$F$, with kernel contained in the center 
of~$\tilde{\mathbf{G}}^*$.

Alternatively, one can start with a connected reductive algebraic 
group~$\tilde{\mathbf{G}}$ with connected center, defined over~$\mathbb{F}$, 
and equipped with a Frobenius morphism~$F$ with respect to some 
$\mathbb{F}_q$-structure on~$\tilde{\mathbf{G}}$. Then if we set $\mathbf{G} := 
[\tilde{\mathbf{G}}, \tilde{\mathbf{G}}]$, the pair $\mathbf{G}$, 
$\tilde{\mathbf{G}}$ satisfies the conditions above.

There is a bijection $\mathbf{L} \mapsto \tilde{\mathbf{L}}$ between the
$F$-stable Levi subgroups of~$\mathbf{G}$ and those of~$\tilde{\mathbf{G}}$,
such that $\mathbf{L} = i^{-1}( \tilde{\mathbf{L}} ) = \tilde{\mathbf{L}} 
\cap \mathbf{G}$ for an $F$-stable Levi subgroup~$\tilde{\mathbf{L}}$ 
of~$\tilde{\mathbf{G}}$. Moreover, the restriction of~$i$ to~$\mathbf{L}$ yields
a regular embedding $\mathbf{L} \rightarrow \tilde{\mathbf{L}}$, 
and~$\mathbf{L}$ is split if and only if~$\tilde{\mathbf{L}}$ is split.
Dually,~$(i^*)^{-1}$ induces a bijection $\mathbf{L}^* \mapsto 
\tilde{\mathbf{L}}^* := (i^*)^{-1}( \mathbf{L}^* )$ between the $F$-stable 
Levi subgroups of~$\mathbf{G}^*$ and those of~$\tilde{\mathbf{G}}^*$.
The $F$-stable Levi subgroups $\mathbf{L}$ and $\mathbf{L}^*$ 
of~$\mathbf{G}$ and $\mathbf{G}^*$ are in duality if and only if 
$\tilde{\mathbf{L}}$ and $\tilde{\mathbf{L}}^*$ are in duality.

Let $\tilde{\chi} \in \Irr(\tilde{G})$. By a result of Lusztig, the restriction
$\Res^{\tilde{G}}_G( \tilde{\chi} )$ is multiplicity free (see
\cite[Section 10]{Lu} and \cite[Proposition 15.11]{CaEn}).
The following easy observation is crucial.

%\addtocounter{thm}{1}
\begin{lem}
\label{ComparingCentralizers}
Let $s \in G^*$ be semisimple and let $\tilde{s} \in \tilde{G}^*$ with
$i^*( \tilde{s} ) = s$. Then ${(i^*)}^{-1}( C^\circ_{\mathbf{G}^*}( s ) ) =
C_{\tilde{\mathbf{G}}^*}( \tilde{s} )$. In particular $C^\circ_{\mathbf{G}^*}( s )$
is contained in a proper split $F$-stable Levi subgroup of~$\mathbf{G}^*$ if
and only if $C_{\tilde{\mathbf{G}}^*}( \tilde{s} )$ is contained in a proper
split $F$-stable Levi subgroup of~$\tilde{\mathbf{G}}^*$.
\end{lem}
\begin{prf}
We have $i^*( C_{\tilde{\mathbf{G}}^*}( \tilde{s} ) ) = C^\circ_{\mathbf{G}^*}( s )$
(see \cite[p.~36]{CeBo2}). As the kernel of~$i^*$ is
contained in $C_{\tilde{\mathbf{G}}^*}( \tilde{s} )$, the first result follows.

As~$i^*$ induces a bijection between the proper split $F$-stable Levi subgroups
of $\tilde{\mathbf{G}}^*$ and of $\mathbf{G}^*$, the second statement follows
from the first.
\end{prf}

\subsection{Restriction of characters}
An irreducible character $\chi \in \Irr(G)$ is called \textit{Harish-Chandra 
imprimitive}, if there is a proper split $F$-stable Levi subgroup~$\mathbf{L}$
of~$\mathbf{G}$, and $\vartheta \in \Irr(L)$ such that $\chi = R_L^G(\vartheta)$.
Otherwise,~$\chi$ is called \textit{Harish-Chandra primitive}.

\begin{lem}
\label{ClifforfAndImprimitivity}
Let~$\tilde{\chi}$ be an irreducible character of~$\tilde{G}$. If one irreducible 
constituent of $\Res^{\tilde{G}}_G( \tilde{\chi} )$ is imprimitive or
Harish-Chandra imprimitive, all of them are.
\end{lem}
\begin{prf}
The irreducible constituents of $\Res^{\tilde{G}}_G( \tilde{\chi} )$ 
are conjugate in~$\tilde{G}$. As conjugation of characters commutes with
induction and Harish-Chandra induction, respectively, our claim follows.
\end{prf}

\smallskip
\noindent
The following well-known result will also be very useful later on.
%\addtocounter{thm}{1}
\begin{lem}\label{7}
Let $\tilde{\mathbf{L}}$ be a split $F$-stable Levi subgroup 
of~$\tilde{\mathbf{G}}$ and $\mathbf{L} := \tilde{\mathbf{L}} \cap \mathbf{G}$. 
Then for every character $\tilde{\vartheta}$ of $\tilde{L}$ and 
every character $\vartheta$ of~$L$ we have
$${\rm Res}^{\tilde{G}}_{G}(R_{\tilde{L}}^{\tilde{G}}(\tilde{\vartheta}))
= R_{L}^{G}({\rm Res}_{L}^{\tilde{L}}(\tilde{\vartheta}))$$
and
$$\Ind_G^{\tilde{G}}( R_L^G( \vartheta ) ) = 
R_{\tilde{L}}^{\tilde{G}}( \Ind_L^{\tilde{L}}( \vartheta ) ).$$
\end{lem}
\begin{prf}
See \cite[Proposition 10.10]{CeBo2}.
\end{prf}

%\addtocounter{subsection}{3}
\subsection{Jordan decomposition of characters}
\label{JordanDecomposition}
Let $s \in G^*$ be semisimple. Choose a semisimple element
$\tilde{s} \in \tilde{G}^*$ with $i^*( \tilde{s} ) = s$. Let $\tilde{\chi}
\in \mathcal{E}( \tilde{G}, [\tilde{s}] )$. Then every irreducible constituent
of $\Res_G^{\tilde{G}}( \tilde{\chi} )$ is contained in $\mathcal{E}( G, [s] )$
(see \cite[Proposition 11.7(a)]{CeBo2}). Also, if $\chi \in 
\mathcal{E}( G, [s] )$, there is $\tilde{\chi} \in 
\mathcal{E}( \tilde{G}, [\tilde{s}] )$ such that $\chi$ is a constituent of
$\Res_G^{\tilde{G}}( \tilde{\chi} )$ (see \cite[Proposition 11.7(b)]{CeBo2}).
Moreover, 
$\Res_G^{\tilde{G}}( \tilde{\chi} )$ is multiplicity free (see 
\cite[Section 10]{Lu} and \cite[Proposition 15.11]{CaEn}).
The conjugation action of $\tilde{G}$
on~$G$ permutes these irreducible constituents, and~$G$ fixes all of these.
Thus~$\tilde{G}$ acts on $\mathcal{E}( G, [s] )$. We write $[\chi]$ for the
$\tilde{G}$-orbit of $\chi \in \mathcal{E}( G, [s] )$ and $c(\chi)$ for the
number of elements in~$[\chi]$. Thus $c(\chi)$ equals
the number of $\tilde{G}$-conjugates of~$\chi$. The orbit of~$\chi$ does not
depend on the chosen regular embedding of~$\mathbf{G}$ (see 
\cite[Corollary 15.14(i)]{CaEn}). On the other hand, $A_{\mathbf{G}^*}( s )^F$ 
acts on $\mathcal{E}( C^\circ_{\mathbf{G}^*}( s )^F, [1] )$ and we 
write~$[\lambda]$ for the $A_{\mathbf{G}^*}( s )^F$-orbit of $\lambda 
\in \mathcal{E}( C^\circ_{\mathbf{G}^*}( s )^F, [1] )$ and 
$A_{\mathbf{G}^*}( s )_\lambda^F$ for its stabilizer in 
$A_{\mathbf{G}^*}( s )^F$. Thus $A_{\mathbf{G}^*}( s )_\lambda^F =
C_{\mathbf{G}^*}(s)_\lambda^F/C^\circ_{\mathbf{G}^*}(s)^F$, where
$C_{\mathbf{G}^*}(s)_\lambda^F$ is the stabilizer of~$\lambda$ in
$C_{\mathbf{G}^*}(s)^F$.

The following notation will also be useful in the sequel. We put
$\tilde{C}_{\tilde{\mathbf{G}}^*}( \tilde{s} ) := (i^*)^{-1}( C_{\mathbf{G}^*}( s ) )$
and $\tilde{A}_{\tilde{\mathbf{G}}^*}( \tilde{s} ) := 
\tilde{C}_{\tilde{\mathbf{G}}^*}( \tilde{s} )/C_{\tilde{\mathbf{G}}^*}( \tilde{s} )$.
Then~$i^*$ induces an isomorphism $\tilde{A}_{\tilde{\mathbf{G}}^*}( \tilde{s} )^F
\rightarrow A_{\mathbf{G}^*}( s )^F$ which commutes with the actions of
$\tilde{A}_{\tilde{\mathbf{G}}^*}( \tilde{s} )^F$ on
$\mathcal{E}( C_{\tilde{\mathbf{G}}^*}( \tilde{s} )^F, [1] )$ and of
$A_{\mathbf{G}^*}( s )^F$ on $\mathcal{E}( C^\circ_{\mathbf{G}^*}( s )^F, [1] )$,
respectively. Let $\tilde{\lambda} \in \mathcal{E}( C_{\tilde{\mathbf{G}}^*}( \tilde{s} )^F, [1] )$.
Identify~$\tilde{\lambda}$ with an element~$\lambda$ of 
$\mathcal{E}( C^\circ_{\mathbf{G}^*}( s )^F, [1] )$ through the
surjective homomorphism $i^*: C_{\tilde{\mathbf{G}}^*}( \tilde{s} )^F \rightarrow
C^\circ_{\mathbf{G}^*}( s )^F$ (see \cite[Proposition 13.20]{DiMi2}).
Write $\tilde{C}_{\tilde{\mathbf{G}}^*}( \tilde{s} )^F_{\tilde{\lambda}}$ for the
stabilizer of~$\tilde{\lambda}$ in $\tilde{C}_{\tilde{\mathbf{G}}^*}( \tilde{s} )^F$ and
$\tilde{A}_{\tilde{\mathbf{G}}^*}( \tilde{s} )_{\tilde{\lambda}}^F := 
\tilde{C}_{\tilde{\mathbf{G}}^*}( \tilde{s} )^F_{\tilde{\lambda}}/C_{\tilde{\mathbf{G}}^*}( \tilde{s} )^F$.
Thus $i^*$ induces an isomorphism between
$\tilde{A}_{\tilde{\mathbf{G}}^*}( \tilde{s} )_{\tilde{\lambda}}^F$
and $A_{\mathbf{G}^*}( s )_\lambda^F$. 

Lusztig has shown in \cite{Lu}, that a Jordan decomposition of characters
between $\mathcal{E}( \tilde{G}, [\tilde{s}] )$ and
$\mathcal{E}( C_{\tilde{\mathbf{G}}^*}( \tilde{s} )^F, [1] )$ induces a bijection
between the $\tilde{G}$-orbits on $\mathcal{E}( G, [s] )$ and the
$A_{\mathbf{G}^*}( s )^F$-orbits on
$\mathcal{E}( C^\circ_{\mathbf{G}^*}( s )^F, [1] )$. Suppose that $[\chi] \subseteq 
\mathcal{E}( G, [s] )$ and $[\lambda] \subseteq 
\mathcal{E}( C^\circ_{\mathbf{G}^*}( s )^F, [1] )$ are corresponding orbits under this 
bijection. We then write $[\chi] \leftrightarrow [\lambda]$. 
The bijection is obtained as follows. First, choose $\tilde{\chi} \in 
\mathcal{E}( \tilde{G}, [\tilde{s}] )$ such that 
$\Res_G^{\tilde{G}}( \tilde{\chi} )$ contains~$\chi$ as a constituent. Next, 
let $\tilde{\lambda} \in \mathcal{E}( C_{\tilde{\mathbf{G}}^*}( \tilde{s} )^F, [1] )$ 
correspond to $\tilde{\chi}$ in the given Jordan decomposition of characters. 
As above, identify~$\tilde{\lambda}$ with an element~$\lambda$ of 
$\mathcal{E}( C^\circ_{\mathbf{G}^*}( s )^F, [1] )$. Then $[\chi]
\leftrightarrow [\lambda]$. An important property of any such bijection
is the fact that 
\begin{equation}
\label{OrbitLength}
c(\chi) = |A_{\mathbf{G}^*}( s )_\lambda^F|,
\end{equation}
if $[\chi] \leftrightarrow [\lambda]$ (see \cite[Proposition 5.1]{Lu}).
In the following, we will call any such bijection {\em Lusztig's generalized 
Jordan decomposition of characters}.

\subsection{The proof of Theorem~\ref{MainResult}(a)}
We begin by generalizing \cite[Theorem 7.3]{HiHuMa}.

%\addtocounter{thm}{1}
\begin{thm}
\label{SameComponentGroup}
Let $s \in G^*$ be semisimple such that $C^\circ_{\mathbf{G}^*}( s ) \leq 
\mathbf{L}^*$, where~$\mathbf{L}^*$ is a split $F$-stable Levi
subgroup of $\mathbf{G}^*$. Let~$\mathbf{L}$ be an $F$-stable Levi subgroup
of $\mathbf{G}$ dual to~$\mathbf{L}^*$.

Let $\chi \in \mathcal{E}( G, [s] )$. Then there is $\vartheta \in 
\mathcal{E}( L, [s] )$ such that~$\chi$ is an irreducible constituent 
of~$R_L^G( \vartheta )$. Any two such elements of $\mathcal{E}( L, [s] )$
are conjugate by an element of~$\tilde{L}$.

The number of $\tilde{L}$-conjugates of~$\vartheta$ is less than or equal to
the number of $\tilde{G}$-conjugates of~$\chi$. If equality holds, then
$\chi = R_L^G( \vartheta )$.
\end{thm}
\begin{prf}
Let $\tilde{s} \in \tilde{G}^*$ be semisimple with $i^*( \tilde{s} ) = s$.
Choose $\tilde{\chi} \in \mathcal{E}( \tilde{G}, [\tilde{s}] )$, such
that~$\chi$ is a constituent of $\Res^{\tilde{G}}_G( \tilde{\chi} )$.
By Lemma~\ref{ComparingCentralizers}, we have 
$C_{\tilde{\mathbf{G}}^*}( \tilde{s} ) = 
(i^*)^{-1}( C^\circ_{\mathbf{G}^*}( s ) ) \leq \tilde{\mathbf{L}}^*$, and
hence $C_{\tilde{\mathbf{G}}^*}( \tilde{s} ) = 
C_{\tilde{\mathbf{L}}^*}( \tilde{s} )$.
By a result of Lusztig (see \cite[(7.9.1)]{lusz2}), Harish-Chandra induction
yields a bijection $\mathcal{E}( \tilde{L}, [\tilde{s}] ) \rightarrow
\mathcal{E}( \tilde{G}, [\tilde{s}] )$. Thus there is $\tilde{\vartheta} \in
\mathcal{E}( \tilde{L}, [\tilde{s}] )$ with 
$R_{\tilde{L}}^{\tilde{G}}( \tilde{\vartheta} ) = \tilde{\chi}$.

By Lemma~\ref{7} we have
\begin{equation}
\label{EquationInSameComponentGroup}
\Res^{\tilde{G}}_G( \tilde{\chi} ) = 
\Res^{\tilde{G}}_G( R_{\tilde{L}}^{\tilde{G}}( \tilde{\vartheta} ) ) =
R_L^G( \Res^{\tilde{L}}_L( \tilde{\vartheta} ) ).
\end{equation}
As~$\chi$ is an irreducible constituent of $\Res^{\tilde{G}}_G( \tilde{\chi} )$,
this gives our first claim. 

Let $\rho \in \mathcal{E}( L, [s] )$ such that~$\chi$ is an irreducible 
constituent of~$R_L^G( \rho )$. By Lemma~\ref{7} we have
$\Ind_G^{\tilde{G}}( R_L^G( \rho ) ) = 
R_{\tilde{L}}^{\tilde{G}}( \Ind_L^{\tilde{L}}( \rho ) )$. As~$\tilde{\chi}$
is a constituent of $\Ind_G^{\tilde{G}}( \chi )$, there is an irreducible
constituent $\tilde{\rho}$ of $\Ind_L^{\tilde{L}}( \rho )$, such 
that~$\tilde{\chi}$ occurs in $R_{\tilde{L}}^{\tilde{G}}( \tilde{\rho} )$.
But then $\tilde{\rho} \in \mathcal{E}( \tilde{L}, [\tilde{s}] )$ (see 
\cite[Proposition 15.7]{CaEn}), and thus 
$R_{\tilde{L}}^{\tilde{G}}( \tilde{\rho} ) = \tilde{\chi}$. It follows that
$\tilde{\rho} = \tilde{\vartheta}$ which gives our second claim.

Since the characters 
$\Res^{\tilde{G}}_G( \tilde{\chi} )$ and 
$\Res^{\tilde{L}}_L( \tilde{\vartheta} )$ are multiplicity free
(see above), the last two assertions follow 
from Equation~(\ref{EquationInSameComponentGroup}).
\end{prf}

\smallskip

\noindent
We also have a kind of converse to Theorem~\ref{SameComponentGroup}.

\begin{thm}
\label{SameComponentGroupConverse}
Let $s \in G^*$ be semisimple, let~$\mathbf{L}^*$ be an $F$-stable split Levi
subgroup of~$\mathbf{G}^*$ containing~$s$, and let~$\mathbf{L}$ be an $F$-stable
split Levi subgroup of~$\mathbf{G}$ dual to~$\mathbf{L}^*$.

Suppose that there is $\chi \in \mathcal{E}( G, [ s ] )$ and $\vartheta 
\in \mathcal{E}( L, [ s ] )$ such that $R_L^G( \vartheta ) = \chi$. 
Let~$\tilde{\vartheta}$ be an irreducible constituent 
of~$\Ind_{L}^{\tilde{L}}( \vartheta )$. Then all irreducible constituents of 
$R_{\tilde{L}}^{\tilde{G}}( \tilde{\vartheta} )$ have the same degree, and 
the number of these  constituents equals~$c(\vartheta)/c(\chi)$.
In particular, $c(\vartheta) \geq c(\chi)$.
If $c(\chi) = c(\vartheta)$, then $C^\circ_{\mathbf{G}^*}( s ) \leq 
\mathbf{L}^*$.
\end{thm}
\begin{prf}
By Lemma~\ref{7} we have
\begin{equation}
\label{EquationFromLemma7}
\Ind_G^{\tilde{G}}( \chi ) = \Ind_G^{\tilde{G}}( R_L^G( \vartheta ) ) = 
R_{\tilde{L}}^{\tilde{G}}( \Ind_L^{\tilde{L}}( \vartheta ) ).
\end{equation}
The restriction of every irreducible character of~$\tilde{G}$ to~$G$ (and 
of~$\tilde{L}$ to~$L$) is multiplicity free (see \cite[Section 10]{Lu} and
\cite[Proposition 15.11]{CaEn}). Clifford theory then implies that the 
irreducible constituents of $\Ind_G^{\tilde{G}}( \chi )$ and of 
$\Ind_L^{\tilde{L}}( \vartheta )$ have degrees $c(\chi)\chi(1)$ and 
$c(\vartheta)\vartheta(1)$, respectively. 
In particular, all constituents of
$R_{\tilde{L}}^{\tilde{G}}( \tilde{\vartheta} )$ have the same degree.
Using the fact that $\tilde{G}/G \cong \tilde{L}/L$ (see
\cite[Corollaire 2.6, D\'emonstration]{CeBo2}) and comparing the degree
of $\vartheta(1)$ with the degree of 
$R_{\tilde{L}}^{\tilde{G}}( \tilde{\vartheta} )$, we obtain
the asserted number of constituents of the latter character.

If $c(\chi) = c(\vartheta)$, 
then~$R_{\tilde{L}}^{\tilde{G}}( \tilde{\vartheta} )$ is irreducible. It 
follows from the first part of the proof of \cite[Theorem 8.4]{HiHuMa} that
$C_{\tilde{\mathbf{G}}^*}( \tilde{s} ) \leq \tilde{\mathbf{L}}^*$. This gives
the third claim by Lemma~\ref{ComparingCentralizers}.
\end{prf}

\begin{cor}
\label{SufficientImprimitivityCondition}
Let $s \in G^*$ be semisimple such that $C^\circ_{\mathbf{G}^*}( s ) \leq 
\mathbf{L}^*$, where~$\mathbf{L}^*$ is a split $F$-stable Levi
subgroup of $\mathbf{G}^*$. Let~$\mathbf{L}$ be an $F$-stable Levi subgroup
of $\mathbf{G}$ dual to~$\mathbf{L}^*$.

Let $\chi \in \mathcal{E}( G, [s] )$ and let~$\lambda$ be a unipotent character 
of $C_{\mathbf{G}^*}^\circ( s )^F$ such that $[\chi]$ corresponds to $[\lambda]$ in 
Lusztig's generalized Jordan decomposition of characters.
Then the following two assertions are equivalent. 
 
{\rm (a)}  We have $C_{\mathbf{G}^*}( s )_\lambda^F \leq \mathbf{L}^*$.

{\rm (b)} There is $\vartheta \in \mathcal{E}( L, [s] )$ with $R_L^G( \vartheta ) 
= \chi$.
\end{cor}

\begin{prf}
Let~$\tilde{s} \in \tilde{G}^*$ be semisimple with $s = i^*( \tilde{s} )$. Let 
$\tilde{\chi} \in \mathcal{E}( \tilde{G}, [\tilde{s}] )$ such that~$\chi$ is an 
irreducible constituent of $\Res^{\tilde{G}}_G( \tilde{\chi} )$. As in the proof 
of Theorem~\ref{SameComponentGroup}, we let $\tilde{\vartheta} \in 
\mathcal{E}( \tilde{L}, [\tilde{s}] )$ with 
$R_{\tilde{L}}^{\tilde{G}}( \tilde{\vartheta} ) = \tilde{\chi}$. Moreover, we
choose an irreducible constituent~$\vartheta$ of 
$\Res^{\tilde{L}}_L( \tilde{\vartheta} )$ such that~$\chi$ occurs 
in~$R_L^G( \vartheta )$.

Composing the Jordan decomposition of characters between
$\mathcal{E}( \tilde{L}, [\tilde{s}] )$ and
$\mathcal{E}( C_{\tilde{\mathbf{L}}^*}( \tilde{s} )^F, [1] )$ with Harish-Chandra 
induction, we obtain the Jordan decomposition between
$\mathcal{E}( \tilde{G}, [\tilde{s}] )$ and
$\mathcal{E}( C_{\tilde{\mathbf{G}}^*}( \tilde{s} )^F, [1] )$. Thus $\tilde{\vartheta}
\in \mathcal{E}( \tilde{L}, [\tilde{s}] )$ and
$R_{\tilde{L}}^{\tilde{G}}( \vartheta ) = \tilde{\chi} \in 
\mathcal{E}( \tilde{G}, [\tilde{s}] )$ correspond to the same unipotent
character of $C_{\tilde{\mathbf{L}}^*}( \tilde{s} )^F = 
C_{\tilde{\mathbf{G}}^*}( \tilde{s} )^F$. Viewed as a character 
of~$C_{\mathbf{G}^*}^\circ(s)^F$,
this lies in the $A_{\mathbf{G}^*}(s)^F$-orbit $[\lambda]$. 
Moreover, $[\vartheta] \leftrightarrow [\lambda]$.

Suppose that~(a) holds. Then $|A_{\mathbf{G}^*}( s )_\lambda^F| = 
|A_{\mathbf{L}^*}( s )_\lambda^F|$, and (b) follows from~(\ref{OrbitLength})
and the last statement of Theorem~\ref{SameComponentGroup}.

Now suppose that~(b) is satisfied. It follows from 
Theorem~\ref{SameComponentGroupConverse} and ~(\ref{OrbitLength}) that 
$|A_{\mathbf{G}^*}( s )_\lambda^F| \geq |A_{\mathbf{L}^*}( s )_\lambda^F|$. 
We thus have equality, implying~(a).
\end{prf}

Notice that assertion~(a) of Theorem~\ref{MainResult} follows from 
Corollary~\ref{SufficientImprimitivityCondition}.
Assertion~(b) of that theorem will be proved case by case.
We start by investigating a special situation.

\begin{cor}
\label{SufficientPrimitivityConditionNew}
Let $s \in G^*$ be semisimple, let $\chi \in \mathcal{E}( G, [ s ] )$ and 
$\lambda \in \mathcal{E}( C_{\mathbf{G}^*}^\circ( s )^F, [ 1 ] )$ such that 
$[\chi] \leftrightarrow [\lambda]$. Assume that $A_{\mathbf{G}^*}( s )^F$ 
fixes~$\lambda$, i.e.\ $A_{\mathbf{G}^*}( s )_\lambda^F = 
A_{\mathbf{G}^*}( s )^F$.

{\rm (a)} Suppose that there is a proper split $F$-stable Levi 
subgroup~$\mathbf{L}^*$ of~$\mathbf{G}^*$ with $s \in \mathbf{L}^*$ and 
$\vartheta \in \mathcal{E}( L, [s] )$ such that $\chi = R_L^G( \vartheta )$
(where~$\mathbf{L}$ is an $F$-stable split Levi subgroup of~$\mathbf{G}$
dual to~$\mathbf{L}^*$). Then $C_{\mathbf{G}^*}(s)^F C^\circ_{\mathbf{G}^*}(s)
\leq \mathbf{L}^*$.

{\rm (b)} If $C_{\mathbf{G}^*}( s )^F$ is not contained in any proper split 
$F$-stable Levi subgroup of~$\mathbf{G}^*$, then~$\chi$ is Harish-Chandra 
primitive.
\end{cor}
\begin{prf}
It suffices to prove~(a).
By~(\ref{OrbitLength}) and our assumption, we have $c(\vartheta) \leq 
|A_{\mathbf{L}^*}( s )^F|$ and $c(\chi) = |A_{\mathbf{G}^*}( s )_\lambda^F| = 
|A_{\mathbf{G}^*}( s )^F|$. 
Theorem~\ref{SameComponentGroupConverse} implies $c(\chi) \leq c(\vartheta)$,
and thus $|A_{\mathbf{G}^*}( s )^F| \leq |A_{\mathbf{L}^*}( s )^F|$. As the
reverse inequality also holds (see~\ref{TheComponentGroup}), we have $c(\chi) 
= c(\vartheta)$, and then $C_{\mathbf{G}^*}^\circ( s ) \leq \mathbf{L}^*$, again 
by Theorem~\ref{SameComponentGroupConverse}. Finally, $C_{\mathbf{G}^*}(s)^F =
C_{\mathbf{G}^*}(s)_\lambda^F$ by assumption, and $C_{\mathbf{G}^*}(s)_\lambda^F
\leq \mathbf{L}^*$ by Corollary~\ref{SufficientImprimitivityCondition}.
\end{prf}

%\addtocounter{subsection}{5}
\subsection{Some results on Weyl groups.}
\label{PrimitivityUnderRestriction}
Before we proceed, we need to strengthen \cite[Lemma~$8.2$]{HiHuMa}. In the 
lemma below we adopt the ususal exponential notation for partitions.

\begin{lem}
\label{WeylGroupInduction}
Let~$W$ be a finite irreducible Coxeter group and let~$W_0$ be a parabolic 
subgroup of~$W$. Furthermore, let $\rho \in \Irr(W_0)$.

{\rm (a)} Suppose that $\Ind_{W_0}^W( \rho ) = \psi_1 + \psi_2$ with 
$\psi_1, \psi_2 \in \Irr(W)$. Then one of the following
holds.

{\rm (i)} We have $W = S_n$ and $W_0 = S_{n-1}$ for some $n \geq 1$, and
$\rho$ is labelled by one of the partitions $(a^b)$ of $n - 1$ (where 
$(a,b) = (0,1)$ if $n = 1$). The characters~$\psi_1$ 
and~$\psi_2$ are then labelled by $(a+1,a^{b-1})$ and $(a^b,1)$, respectively.

{\rm (ii)} We have $W = S_n$ and $W_0 = S_{n-k} \times S_k$ for some $n \geq 4$
and some $1 < k < n - 1$, and $\rho$ is labelled by $(1^{n-k}) \times (k)$.
In this case,~$\psi_1$ and~$\psi_2$ are labelled by $(k+1,1^{n-k-1})$ and 
$(k, 1^{n-k})$, respectively.

{\rm (iii)} We have~$W$ of type $D_m$ and $W_0$ of type~$D_{m-1}$ for some odd 
$m \geq 5$, and $\rho$ is one of the two characters labelled by an 
unordered pair of partitions of the form $\{ (a^b), (a^b) \}$, where 
$ab = (m-1)/2$. Moreover, the unordered pairs of partitions 
$\{ (a^b), (a+1,a^{b-1}) \}$ and $\{ (a^b), (a^b,1) \}$ label~$\psi_1$ 
and~$\psi_2$, respectively.

{\rm (b)} Assume that~$W$ is of type $A_1$,~$A_2$,~$B_2$,~$G_2$,~$B_3$ 
or~$D_4$.
Suppose further that $\Ind_{W_0}^W( \rho ) = \psi_1 + \psi_2 + \psi_3$ with 
(not necessarily distinct) $\psi_i \in \Irr(W)$, $i = 1, 2, 3$. Then~$W$ is of 
type~$B_2$,~$B_3$ or~$D_4$ and the labels of $\psi_1$, $\psi_2$, $\psi_3$ 
are as given in the following lists, where each line corresponds to one
pair~$(W_0,\rho)$.
$$ B_2: \begin{array}{ccc} \\ \hline\hline 
                           \psi_1 & \psi_2 & \psi_3 \\ \hline\hline
                           (1^2,-) & (1,1) & (-,1^2) \rule[  0pt]{0pt}{ 13pt} \\
                           (1,1) & (2,-) & (-,2) \\
                           (1,1) & (-,1^2) & (-,2) \\
                           (1^2,-) & (1,1) & (2,-) \rule[- 5pt]{0pt}{13pt} \\ \hline\hline
\end{array}
\quad B_3: \begin{array}{ccc} \\ \hline\hline 
                           \psi_1 & \psi_2 & \psi_3 \\ \hline\hline
                           (1^3,-) & (1^2,1) & (21,-) \rule[  0pt]{0pt}{ 13pt} \\
                           (1,1^2) & (-,1^3) & (-,21) \\
                           (21,-) & (2,1) & (3,-) \\
                           (1,2) & (-,21) & (-,3) \rule[- 5pt]{0pt}{13pt} \\ \hline\hline
\end{array}$$
$$ D_4: \begin{array}{ccc} \\ \hline\hline 
                           \psi_1 & \psi_2 & \psi_3 \\ \hline\hline
                           (1^2,1^2)^- & (1,1^3) & (-,1^4) \rule[  0pt]{0pt}{ 13pt} \\
                           (2,2)^- & (1,3) & (-,4) \\
                           (1^2,1^2)^+ & (1,1^3) & (-,1^4) \\
                           (2,2)^+ & (1,3) & (-,4) \\
                           (1,1^3) & (-,1^4) & (-,21^2) \\
                           (1,3) & (-,31) & (-,4) \rule[- 5pt]{0pt}{13pt} \\ \hline\hline
\end{array}$$

\medskip

{\rm (c)} Suppose that $(W, \psi_1, \psi_2)$ is as in {\rm (a)}.
For $1 \leq i \leq 2$ let~$a_i$ denote the value of Lusztig's 
$\mathbf{a}$-function associated to the irreducible representation~$\psi_i$ 
of~$W$ (see \cite[($4.1.1$)]{luszbuch}).
Then $a_1 \neq a_2$.

{\rm (d)} Suppose that $(W, \psi_1,\psi_2,\psi_3)$ is as in {\rm (b)}.
For $1 \leq i \leq 3$ let~$d_i$ denote the generic degree associated to the 
irreducible representation~$\psi_i$ (see \cite[p.~$61$]{luszbuch} or 
\cite[$8.1.8$]{GP}). Then~$d_i$ is a polynomial in one variable if~$W$ is 
of type~$D_4$, and in two variables, otherwise. Let~$k$ and~$\ell$ be positive 
integers. Then, in the first case, $|\{ d_i( q^k ) \mid i = 1, 2, 3 \}| > 1$, 
and in the second case, $|\{ d_i( q^k, q^\ell ) \mid i = 1, 2, 3 \}| > 1$.
\end{lem}
\begin{prf}
(a) Let $(W,W_0,\rho)$ satisfy the hypothesis of (a). If~$W$ is a dihedral 
group of order~$2d$, then $d \leq 3$, as~$W_0$ has index at least~$d$, and the 
character degrees of~$W$ are~$1$ and~$2$. Thus~$W$ is of type~$A_2$ and cannot
be of type~$G_2$. Using CHEVIE (see~\cite{chevie}), it is easy to check that~$W$ 
is not one of the Coxeter groups~$H_3$ or~$H_4$. Now \cite[Lemma~$8.2$]{HiHuMa} 
implies that~$W_0$ is a maximal parabolic subgroup of~$W$. 

A further application of CHEVIE shows that~$W$ is not one of the 
exceptional Weyl groups~$F_4$, or~$E_i$, $i = 6, 7, 8$. Next,
\cite[Lemma~$8.1$]{HiHuMa} rules out the case that~$W$ is of type~$B_m$ for
some $m \geq 2$. It thus remains to consider the cases that~$W$ is of 
type~$A_{n-1}$, i.e.\ $W = S_n$ for $n > 1$, or of type~$D_m$ for $m \geq 4$. 
We begin by investigating the case $W = S_n$, $n \geq 1$, and $W_0 = S_{n-1}$.
In this case, the branching rule (see e.g.\ \cite[$2.4.3$]{jake}) easily gives 
our claim. Now suppose that $W = S_n$, $n \geq 4$, and $W_0 = S_{n-k} \times S_k$
for some $1 < k < n - 1$. Here, we use the Littlewood-Richardson rule (see the 
version of \cite[Corollary~$2.8.14$]{jake}) to show that exactly the cases 
listed in~(ii) satisfy the hypothesis. 
%%%%%%%%%%%%%%%%%%%%%%%%%%%%%%%%%%%%%%%%%%%%%%%%%%%%%%%%%%%%%%%%%%%%%%%%%%%%%%%%
%%%%%%%%%%%%%%%%%%%%%%%%%%%%%%%%%%%%%%%%%%%%%%%%%%%%%%%%%%%%%%%%%%%%%%%%%%%%%%%%
%%
%% Für mehr Details zu (a) siehe Notizen vom 07.04.2016, Seiten 1 - 5.
%%
%%%%%%%%%%%%%%%%%%%%%%%%%%%%%%%%%%%%%%%%%%%%%%%%%%%%%%%%%%%%%%%%%%%%%%%%%%%%%%%%
%%%%%%%%%%%%%%%%%%%%%%%%%%%%%%%%%%%%%%%%%%%%%%%%%%%%%%%%%%%%%%%%%%%%%%%%%%%%%%%%
Finally, suppose that~$W$ is of 
type~$D_m$ for some $m \geq 4$. We proceed as in the last two paragraphs of 
the proof of \cite[Lemma~$8.1$]{HiHuMa}. Namely, we embed~$W$ into a Weyl 
group~$\hat{W}$ of type~$B_m$. Then~$W_0$ is of the form 
$W_0 = W' \times S_k$, for some $1 \leq k \leq m-1$, where~$W'$ denotes a Weyl 
group of type~$D_{m-k}$, naturally embedded into~$W$. Further, $\hat{W}_0 :=
\hat{W}' \times S_k$ is a maximal parabolic subgroup of~$\hat{W}$, where~$W'$
is embedded into~$\hat{W}'$ in the same way as~$W$ is embedded into~$\hat{W}$. 
Now
\begin{equation}
\label{InductionEquation}
\Ind_{\hat{W}_0}^{\hat{W}}( \Ind_{W_0}^{\hat{W_0}}( \rho ) ) =
\Ind_{W_0}^{\hat{W}}( \rho ) = \Ind_W^{\hat{W}}( \Ind_{W_0}^W( \rho ) )
= \Ind_W^{\hat{W}}( \psi_1 + \psi_2 ).
\end{equation}
Hence $\Ind_{\hat{W}_0}^{\hat{W}}( \Ind_{W_0}^{\hat{W_0}}( \rho ))$ has at 
most four irreducible constituents, as~$W$ has index~$2$ in~$\hat{W}$. It 
follows from \cite[Lemma~$8.1$]{HiHuMa}
that $k = 1$ and also that $\hat{\rho} := 
\Ind_{W_0}^{\hat{W_0}}( \rho )$ is irreducible. Thus~$\hat{\rho}$ is 
labelled by an unordered pair of partitions of the form $\{ \pi, \pi \}$. In 
particular, $m - 1$ is even. If~$\pi$ is not of the form $\pi = (a^b)$ with 
$ab = (m - 1)/2$, then $\Ind_{\hat{W}_0}^{\hat{W}}( \hat{\rho} )$ has more 
than four irreducible constituents by the branching rule for irreducible 
characters of~$\hat{W}$. This contradicts~(\ref{InductionEquation}), proving 
our claim.

\begin{table}[t]
\caption{\label{GenericDegreesB3} Some generic degrees in $B_3$}
$$\begin{array}{cc} \hline\hline
\text{Bipartition} & \text{Generic degree} \rule[- 7pt]{0pt}{ 22pt} \\ \hline\hline
(3, - ) & 1 \rule[ 0pt]{0pt}{ 13pt} \\
(21,-)  & \frac{X^2(X^2Y+1)(X+1)}{X+Y} \\
(1^3,-) & \frac{X^6(XY+1)(X^2Y+1)}{(X+Y)(X^2+Y)} \\
(2,1)   & \frac{XY(XY+1)(X^2+X+1)}{X+Y} \\
(1^2,1) & \frac{X^3Y(X^2Y+1)(X^2+X+1)}{X^2+Y} \\
(1,2)   & \frac{XY^2(X^2Y+1)(X^2+X+1)}{X^2+Y} \\
(1,1^2) & \frac{X^3Y(XY + 1)(X^2+X+1)}{X+Y} \\
(-,3)   & \frac{Y^3(XY+1)(X^2Y+1)}{(X+Y)(X^2+Y)} \\
(-,21)  & \frac{X^2Y^3(X^2Y + 1)(X+1)}{X+Y} \\
(-,1^3) & X^6Y^3 \rule[- 8pt]{0pt}{18pt} \\ \hline\hline
\end{array}
$$
\end{table}

(b) This is easily checked with CHEVIE. 

(c) The~$a_i$ can be computed from the labels of the $\psi_i$ (see 
\cite[($4.4$)]{luszbuch} for~$W$ of type~$A$, and \cite[($4.6$)]{luszbuch} 
for~$W$ of type~$D$).
%%%%%%%%%%%%%%%%%%%%%%%%%%%%%%%%%%%%%%%%%%%%%%%%%%%%%%%%%%%%%%%%%%%%%%%%%%%%%%%%
%%%%%%%%%%%%%%%%%%%%%%%%%%%%%%%%%%%%%%%%%%%%%%%%%%%%%%%%%%%%%%%%%%%%%%%%%%%%%%%%
%%
%% Für mehr Details zu (c) siehe Notizen vom 13.04.2016, Seiten 1 - 4.
%%
%%%%%%%%%%%%%%%%%%%%%%%%%%%%%%%%%%%%%%%%%%%%%%%%%%%%%%%%%%%%%%%%%%%%%%%%%%%%%%%%
%%%%%%%%%%%%%%%%%%%%%%%%%%%%%%%%%%%%%%%%%%%%%%%%%%%%%%%%%%%%%%%%%%%%%%%%%%%%%%%%

(d) The first part for the possibilities for $\{ \psi_1, \psi_2, \psi_3 \}$ is
again checked with CHEVIE. The generic degrees for~$W$ can be found, for 
example, in \cite[$13.5$]{cart}. If~$W$ is of type~$B_2$, the generic degrees 
for the characters labelled by $(2,-)$, $(1,1)$ and $(-,1^2)$ are 
$1$, $XY(X+1)(Y+1)/(X+Y)$ and $X^2Y^2$, respectively. Substituting $(X,Y)$ by
$(q^k,q^\ell)$ in the second of these polynomials will neither evaluate to~$1$ 
nor to $q^{2(k+\ell)}$. This yields our claim for~$W$ of type~$B_2$.
Now assume that $W$ is of type~$B_3$. The generic degrees of~$W$ can be computed
with CHEVIE and are given in Table~\ref{GenericDegreesB3}. The proof is now 
completed as in the precious case. We omit the details.
%%%%%%%%%%%%%%%%%%%%%%%%%%%%%%%%%%%%%%%%%%%%%%%%%%%%%%%%%%%%%%%%%%%%%%%%%%%%%%%%
%%%%%%%%%%%%%%%%%%%%%%%%%%%%%%%%%%%%%%%%%%%%%%%%%%%%%%%%%%%%%%%%%%%%%%%%%%%%%%%%
%%
%% Für mehr Details zu (d) im Fall B3 siehe Notizen vom 21.04.2016, Seiten 2 - 5.
%%
%%%%%%%%%%%%%%%%%%%%%%%%%%%%%%%%%%%%%%%%%%%%%%%%%%%%%%%%%%%%%%%%%%%%%%%%%%%%%%%%
%%%%%%%%%%%%%%%%%%%%%%%%%%%%%%%%%%%%%%%%%%%%%%%%%%%%%%%%%%%%%%%%%%%%%%%%%%%%%%%%
Finally, let~$W$ be of type~$D_4$. If $d_i \in \mathbb{Z}[X]$, write 
$$d_i(X) = c_iX^{a_i} + \text{\ sum of higher powers of\ } X$$ 
with $a_i$ a non-negative integer and $0 \neq c_i \in \mathbb{Z}$.
Computing the generic degrees for the characters of~$W$, we observe: For each 
possible set $\{ \psi_1, \psi_2, \psi_3 \}$ there are $i \neq j \in 
\{ 1, 2, 3 \}$ such that $d_i, d_j \in \mathbb{Z}[X]$ with $a_i 
\neq a_j$ and $c_i = c_j = 1$. Hence
$d_i(q^k)$ and $d_j(q^k)$ are divisible by different powers of~$q$, and thus 
cannot be equal.
%%%%%%%%%%%%%%%%%%%%%%%%%%%%%%%%%%%%%%%%%%%%%%%%%%%%%%%%%%%%%%%%%%%%%%%%%%%%%%%%
%%%%%%%%%%%%%%%%%%%%%%%%%%%%%%%%%%%%%%%%%%%%%%%%%%%%%%%%%%%%%%%%%%%%%%%%%%%%%%%%
%%
%% Für mehr Details zu (d) im Fall D4 siehe Notizen vom 03.04.2016, Seite 1, 
%% sowie vom 12.04.2016, Seiten 2 - 3.
%%
%%%%%%%%%%%%%%%%%%%%%%%%%%%%%%%%%%%%%%%%%%%%%%%%%%%%%%%%%%%%%%%%%%%%%%%%%%%%%%%%
%%%%%%%%%%%%%%%%%%%%%%%%%%%%%%%%%%%%%%%%%%%%%%%%%%%%%%%%%%%%%%%%%%%%%%%%%%%%%%%%
\end{prf}

\begin{cor}
\label{CorollaryToWeylGroupInduction}
Resume the notation of~{\rm \ref{RegularEmbeddings}}.
Let~$\tilde{\mathbf{L}}$ denote a split $F$-stable Levi subgroup 
of~$\tilde{\mathbf{G}}$, and let $\tilde{\vartheta} \in \Irr( \tilde{L} )$.

{\rm (a)} If $q$ is odd, then $R_{\tilde{L}}^{\tilde{G}}( \tilde{\vartheta} )$ 
is not a sum of two irreducible characters of the same degree.

{\rm (b)} Suppose that $\tilde{\mathbf{G}}$ is of type~$E_6$. Let 
$\tilde{s} \in \tilde{L}^*$ be semisimple such that $\tilde{\vartheta} \in 
\mathcal{E}( \tilde{L}, [\tilde{s}] )$, and put $s := i^*(\tilde{s}) \in L^*$. 
If $A_{\mathbf{G}^*}( s )^F \neq 1$, then
$R_{\tilde{L}}^{\tilde{G}}( \tilde{\vartheta} )$ is not a sum of three 
irreducible characters of equal degrees.
\end{cor}
\begin{prf}
Let $(\tilde{L}_0, \tilde{\vartheta}_0 )$ be a cuspidal pair
below~$\tilde{\vartheta}$. Then
$W_{\tilde{G}}( \tilde{L}_0, \tilde{\vartheta}_0 )$ is a Coxeter group
containing $W_{\tilde{L}}( \tilde{L}_0, \tilde{\vartheta}_0 )$ as a parabolic
subgroup (see \cite[Theorem 8.6]{luszbuch}).
Let~$\rho$ be an irreducible character
of~$W_{\tilde{L}}( \tilde{L}_0, \tilde{\vartheta}_0 )$ corresponding
to~$\tilde{\vartheta}$ via Harish-Chandra theory. By the Howlett-Lehrer 
comparison theorem~\cite[Theorem~$5.9$]{HowLeh2}, there is a multiplicity
preserving bijection between the irreducible constituents of
 $\Ind_{W_{\tilde{L}}( \tilde{L}_0, \tilde{\vartheta}_0 )}
     ^{W_{\tilde{G}}( \tilde{L}_0, \tilde{\vartheta}_0 )}( \rho )$
and those of~$R_{\tilde{L}}^{\tilde{G}}( \tilde{\vartheta} )$. The formula of
\cite[Corollary~$8.7$]{luszbuch} gives the degrees of the latter characters
in terms of the generic degrees associated to the irreducible constituents of
$\Ind_{W_{\tilde{L}}( \tilde{L}_0, \tilde{\vartheta}_0 )}
     ^{W_{\tilde{G}}( \tilde{L}_0, \tilde{\vartheta}_0 )}( \rho )$.
More precisely, the generic degrees are evaluated at the parameters of the 
Iwahori-Hecke algebra 
$\End_{\mathbb{C}\tilde{G}}( R_{\tilde{L}_0}^{\tilde{G}}( \tilde{\vartheta}_0 ) )$.
By \cite[Theorem $8.6$]{luszbuch}, these parameters are positive powers of~$q$.
If two generic degrees of $W_{\tilde{G}}( \tilde{L}_0, \tilde{\vartheta}_0 )$ 
evaluate to different integers, then the corresponding characters of
$R_{\tilde{L}_0}^{\tilde{G}}( \tilde{\vartheta}_0 )$ have different degrees.

(a) Suppose that $R_{\tilde{L}}^{\tilde{G}}( \tilde{\vartheta} )$ 
has exactly two constituents. Under this hypothesis, one of the irreducible 
Weyl groups occuring in
the decomposition of $W_{\tilde{G}}( \tilde{L}_0, \tilde{\vartheta}_0 )$
satisfies the hypothesis of Lemma~\ref{WeylGroupInduction}(a). As the generic 
degrees behave well with respect to the factorisation of a Weyl group into
irreducibles, we may assume, for the sake of the argument, that 
$W_{\tilde{G}}( \tilde{L}_0, \tilde{\vartheta}_0 )$ is irreducible. Then
$(W_{\tilde{G}}( \tilde{L}_0, \tilde{\vartheta}_0 ), 
  W_{\tilde{L}}( \tilde{L}_0, \tilde{\vartheta}_0, \rho ) )$ is one of the 
pairs listed in Lemma~\ref{WeylGroupInduction}(a)(i)--(iii). In this
case, the generic degree associated to a character~$\psi$ of 
$W_{\tilde{G}}( \tilde{L}_0, \tilde{\vartheta}_0 )$ is of the form
$$X^{a_\psi}/f_\psi + \text{\ sum of higher powers of\ } X,$$
where $f_\psi = 1$ if $W_{\tilde{G}}( \tilde{L}_0, \tilde{\vartheta}_0 )$ is
of type~$A$, and~$f_\psi$ is a power of~$2$, otherwise 
(see \cite[Corollary~$9.3.6$]{GP} and \cite[(4.1.1),(4.14.2)]{luszbuch}).
Here,~$a_\psi$ is the value of Lusztig's $\mathbf{a}$-function associated 
to~$\psi$.  As~$q$ is odd, two generic degrees with different 
$\mathbf{a}$-values evaluate to different integers under substituting the 
indeterminate $X$ by a positive power of~$q$. An application of 
Lemma~\ref{WeylGroupInduction}(c) proves our claim.

(b) Suppose that $R_{\tilde{L}}^{\tilde{G}}( \tilde{\vartheta} )$   
has exactly three constituents.
By \cite[(8.5.4),(8.5.7),(8.5.8)]{luszbuch}, the group 
$W_{\tilde{G}}( \tilde{L}_0, \tilde{\vartheta}_0 )$ is isomorphic to the $F$-fixed
points of a relative Weyl group in $C_{\tilde{\mathbf{G}}^*}( \tilde{s} )$.
Let us write $W(\tilde{s})$ for the Weyl group of 
$C_{\tilde{\mathbf{G}}^*}( \tilde{s} )$. Using the assumption that 
$A_{\mathbf{G}^*}(s)^F \neq 1$, the tables of L{\"u}beck~\cite{LL} show that 
the irreducible components 
of~$W(\tilde{s})^F$ have rank at most~$2$ or else have type~$B_3$ or~$D_4$. 
It is easy to check with CHEVIE that the same statement then holds for the
relative Weyl groups occuring, i.e.\ for 
$W_{\tilde{G}}( \tilde{L}_0, \tilde{\vartheta}_0 )$. As in the proof of~(a) we may
assume that $W_{\tilde{G}}( \tilde{L}_0, \tilde{\vartheta}_0 )$ is irreducible and
that $W_{\tilde{L}}( \tilde{L}_0, \tilde{\vartheta}_0 )$ is a maximal parabolic
subgroup of $W_{\tilde{G}}( \tilde{L}_0, \tilde{\vartheta}_0 )$. Our claim now
follows from Lemma~\ref{WeylGroupInduction}(b),(d).
\end{prf}

%\addtocounter{thm}{1}
\subsection{Special cases of Theorem~\ref{MainResult}(b)}
\label{TowardsTheorem11b}
Let~$\mathbf{G}$ and~$F$ be such that~$G$ is one of the following groups:
\begin{itemize}
\item[(a)] a symplectic group $\Sp_{n}(q)$ with~$q$ odd and $n \geq 4$ even;
\item[(b)] a spin group $\Spin_{n+1}(q)$ with~$q$ odd and $n \geq 6$ even;
\item[(c)] the universal covering group $E_6(q)_{\rm sc}$ of the simple 
Chevalley group $E_6(q)$;
\item[(d)] the covering group ${^2\!E}_6(q)_{\rm sc}$ (which is the universal 
covering group if $q > 2$) of the simple twisted Steinberg group ${^2\!E}_6(q)$;
\item[(e)] the universal covering group $E_7(q)_{\rm sc}$ of the simple
Chevalley group $E_7(q)$.
\end{itemize}
For the remainder of this subsection we will assume that $(\mathbf{G},F)$ is
one of the pairs introduced above.
Choose $(\tilde{\mathbf{G}},F)$ and the notation as in~\ref{RegularEmbeddings}.

\begin{thm}\label{RestrictionToCommutator}
Let~$\tilde{\mathbf{G}}$ be one of the groups introduced above.
If~$\tilde{\chi}$ is an irreducible Harish-Chandra primitive character 
of~$\tilde{G}$, then every irreducible constituent of $\Res^{\tilde{G}}_
G( \tilde{\chi} )$ is Harish-Chandra primitive.
\end{thm}
\begin{prf}
Suppose that an irreducible 
constituent~$\chi$ of $\Res^{\tilde{G}}_G( \tilde{\chi} )$ is Harish-Chandra
imprimitive. There is a semisimple element $\tilde{s} \in \tilde{G}^*$ such
that $\tilde{\chi} \in \mathcal{E}( \tilde{G}, [\tilde{s}] )$ and $\chi \in
\mathcal{E}( G, [s] )$, where $s := i^*( \tilde{s} ) \in G^*$.
Let~$\mathbf{L}$ be a proper split $F$-stable Levi subgroup of~$\mathbf{G}$ 
such that~$L$ contains an irreducible character~$\vartheta$ with 
$R_L^G( \vartheta ) = \chi$. We may assume that $s \in L^*$ for some $F$-stable 
Levi subgroup $\mathbf{L}^*$ of~$\mathbf{G}^*$ dual to~$\mathbf{L}$, and that 
$\vartheta \in \mathcal{E}( L, [s] )$ (see \cite[Proposition 15.7]{CaEn}). 
Using the notation of Theorem~\ref{SameComponentGroupConverse}, we have $c(\chi)
< c(\vartheta)$. (Otherwise, $C^\circ_{\mathbf{G}^*}( s ) \leq \mathbf{L}^*$, 
implying that $C_{\tilde{\mathbf{G}}^*}( \tilde{s} ) = 
(i^*)^{-1}( C^\circ_{\mathbf{G}^*}( s ) ) \leq \tilde{\mathbf{L}}^*$, and then 
$\mathcal{E}( \tilde{G}, [\tilde{s}] )$ would contain only Harish-Chandra 
imprimitive characters by \cite[Theorem~$7.3$]{HiHuMa}.) Now~$c(\vartheta)$
divides~$|A_{\mathbf{L}^*}( s )|$ by~(\ref{OrbitLength}), 
and~$A_{\mathbf{L}^*}( s )$ is isomorphic to a subgroup of 
$Z(\mathbf{G})/Z^\circ(\mathbf{G})$ (see~\ref{TheComponentGroup} and
\cite[Lemme~$8.3$]{CeBo2}). This is a group of order~$2$ if~$\mathbf{G}$ is as 
in~(a) or~(b), a group of order~$1$ or~$2$ if~$\mathbf{G}$ is as in~(e) and~$q$ 
is even or odd, respectively, and a group of order~$1$ or~$3$
if~$\mathbf{G}$ is of type~$E_6$. It follows that $1 = c(\chi) < c(\vartheta) 
= |Z(\mathbf{G})/Z^\circ(\mathbf{G})|$. In 
particular,~$q$ is odd if~$\mathbf{G}$ is of type~$E_7$. The last 
assertion of Theorem~\ref{SameComponentGroupConverse} 
states that~$c(\vartheta)$ equals the number of irreducible constituents of 
$R_{\tilde{L}}^{\tilde{G}}( \tilde{\vartheta} )$, and that all of these 
constituents have the same degree. This contradicts
Corollary~\ref{CorollaryToWeylGroupInduction}.
\end{prf}

\begin{cor}
\label{ThirdLusztigSeriesResult}
Let $s \in G^*$ be semisimple. If $\mathcal{E}(G,[s])$ contains a Harish-Chandra
imprimitive element, then~$C^\circ_{\mathbf{G}^*}(s)$ is contained in a proper 
split $F$-stable Levi subgroup of~$\mathbf{G}^*$.
\end{cor}
\begin{prf}
Let $\chi \in \mathcal{E}(G,[s])$ be Harish-Chandra imprimitive.
Choose $\tilde{s} \in \tilde{G}^*$ with $i^*( \tilde{s} ) = s$. 
By \cite[Proposition 11.7(b)]{CeBo2}, there is an irreducible character 
$\tilde{\chi} \in \mathcal{E}(\tilde{G},[\tilde{s}])$ such that $\chi$ is a 
constituent of the restriction of $\tilde{\chi}$ to~$G$. It follows from 
Theorem~\ref{RestrictionToCommutator} that~$\tilde{\chi}$ is Harish-Chandra
imprimitive. 

By \cite[Theorem~8.4]{HiHuMa}, we have $C_{\tilde{\mathbf{G}}^*}( \tilde{s} ) 
\leq \tilde{\mathbf{L}}^*$ for a proper split $F$-stable Levi 
subgroup~$\tilde{\mathbf{L}}^*$ of $\tilde{\mathbf{G}}^*$. Now
$i^*(C_{\tilde{\mathbf{G}}^*}( \tilde{s} )) = C^\circ_{\mathbf{G}^*}(s)$ by
Lemma~\ref{ComparingCentralizers}, and thus $C^\circ_{\mathbf{G}^*}(s) \leq 
i^*( \tilde{\mathbf{L}}^* )$. The latter is a proper split $F$-stable Levi 
subgroup of $\mathbf{G}^*$.
\end{prf}

\smallskip

\noindent
We can now prove Theorem~\ref{MainResult}(b) for the groups considered in this 
subsection.

\begin{cor}
\label{ProofOfTheorem11bInSpecialCaseNew}
Theorem~{\rm \ref{MainResult}(b)} holds for~$G$.
\end{cor}
\begin{prf}
Let $s \in G^*$ and let $\chi \in \mathcal{E}( G, [s] )$ be Harish-Chandra
imprimitive. Let $\lambda \in \mathcal{E}( C^\circ_{\mathbf{G}^*}(s)^F, [ 1 ] )$
with $[\chi] \leftrightarrow [\lambda]$. By
Corollary~\ref{ThirdLusztigSeriesResult}, there is a proper
split $F$-stable Levi subgroup $\mathbf{L}^* \leq \mathbf{G}^*$ with
$C^\circ_{\mathbf{G}^*}(s) \leq \mathbf{L}^*$. Suppose first that
$A_{\mathbf{G}^*}(s)_\lambda^F = \{1\}$. Then $C_{\mathbf{G}^*}(s)^F_\lambda
 = C^\circ_{\mathbf{G}^*}(s)^F \leq \mathbf{L}^*$
and our claim follows.

Now assume that $A_{\mathbf{G}^*}(s)_\lambda^F \neq \{1\}$. 
As~$A_{\mathbf{G}^*}(s)$ is isomorphic to a subgroup of 
$Z(\mathbf{G})/Z^\circ(\mathbf{G})$ (see \cite[Lemme~$8.3$]{CeBo2}) and this is 
a group of prime order, we conclude that $A_{\mathbf{G}^*}(s)_\lambda^F = 
A_{\mathbf{G}^*}(s)^F$. By Corollary~\ref{SufficientPrimitivityConditionNew}(a),
Theorem~\ref{MainResult}(b) holds in this case as well.
\end{prf}

\smallskip

\noindent
Although we know that Theorem~\ref{MainResult} holds for the groups
considered here, it remains to determine the pairs $(s, \lambda)$ for
which Condition~(\ref{InclusionCondition}) is satisfied. It also remains to
prove Theorem~\ref{MainResult} for the other quasisimple groups not treated
here, namely $\SL_n(q)$, $\SU_n(q)$ and $\Spin^\pm_{2m}(q)$. Both tasks will 
be achieved in Section~\ref{QuasisimpleGroupsOfLieType}.

\section{Preliminary results on some classical groups}

Recall that~$p$ denotes a prime number,~$q$ a power of~$p$, and~$\mathbb{F}$ 
an algebraic closure of the field with~$p$ elements. The finite fields of
characteristic~$p$ are viewed as subfields of~$\mathbb{F}$.

\subsection{Polynomials and linear transformations}

We collect some notations and results on polynomials, linear transformations 
and centralizers of classical groups needed later on.

\subsubsection{Polynomials}
\label{SubsectionPreliminaries}
Let $\mathbb{F}[X]^0$ denote the set of monic,
separable polynomials over $\mathbb{F}$ in the indeterminate~$X$, with
non-zero constant coefficient. An element of $\mathbb{F}[X]^0$ is uniquely
determined by its set of roots.
For $\mu \in \mathbb{F}[X]^0$ of degree~$d$ and $\alpha \in \mathbb{F}^*$
we put
$$\mu^{\alpha} := \alpha^{d} \mu( X/\alpha ).$$
This defines an action of~$\mathbb{F}^*$ on $\mathbb{F}[X]^0$.
Notice that $\zeta \in \mathbb{F}^*$ is a root of~$\mu$ if and only if
$\alpha\zeta$ is a root of $\mu^{\alpha}$. In case $\alpha = -1$, we
put $\mu' := \mu^\alpha$. Thus the roots of~$\mu'$ are the negatives of 
the roots of~$\mu$. We write $\mu^*$ for the element of 
$\mathbb{F}[X]^0$ whose roots are the
inverses of the roots of~$\mu$. Finally, we put
$$\mu^{*\alpha} := (\mu^*)^\alpha.$$
Thus $\zeta \in \mathbb{F}^*$ is a root of~$\mu$ if and only if
$\alpha\zeta^{-1}$ is a root of~$\mu^{*\alpha}$.
Then $(\mu^{*\alpha})^{*\alpha} = \mu$.
In addition, for $\mu \in \mathbb{F}[X]^0$, we write $\mu^{\dagger}$ for the
element of~$\mathbb{F}[X]^0$ whose roots are the $(-q)$ths powers of the roots
of~$\mu$.

%\addtocounter{thm}{1}
\begin{lem}
\label{MuLemma}
Let~$\alpha \in \mathbb{F}_q^*$, and let $\mu \in \mathbb{F}_q[X] \cap
\mathbb{F}[X]^0$ be of degree~$d$.

{\rm (a)} Suppose that $\mu = \mu^{*\alpha}$. 
If $\zeta^2 \neq \alpha$ for every root $\zeta \in \mathbb{F}$ 
of~$\mu$, then $d = 2e$ is even, and the product of the roots of~$\mu$
equals $\alpha^e$.

If $\mu = \mu^{*\alpha}$ and~$\mu$ is irreducible in $\mathbb{F}_q[X]$, one of 
the following cases occurs.

\begin{itemize}
\item[{\rm (I)}] There is $\zeta \in \mathbb{F}_q$ with $\alpha = \zeta^2$,
$\mu = X - \zeta$ and $\mu' = X + \zeta$.
\item[{\rm (II)}] There is $\zeta \in \mathbb{F}_{q^2} \setminus \mathbb{F}_q$
with $\alpha = \zeta^2$ and $\mu = X^2 - \alpha = \mu'$.
\item[{\rm (III)}] There is $\zeta \in \mathbb{F}_{q^2} \setminus \mathbb{F}_q$
with $\alpha = - \zeta^2$ and $\mu = X^2 + \alpha = \mu'$.
\item[{\rm (IV)}] For all roots $\zeta \in \mathbb{F}$ of~$\mu$ we have
$\alpha \neq \pm \zeta^2$. In this case, $\zeta^{q^e} = \alpha\zeta^{-1}$. 
Moreover, $\mu \neq \mu'$ if $q$ is odd.
\end{itemize}
Cases~{\rm (II)} and~{\rm (III)} only occur if~$q$ is odd.

{\rm (b)} If~$\mu$ is irreducible,~$d$ and~$q$ are odd, and $\mu' = \mu^{*\alpha}$, 
then~$-\alpha$ is a square in~$\mathbb{F}_q$.
\end{lem}
\begin{prf}
(a) Let~$\zeta$ be a root of~$\mu$ in $\mathbb{F}$. If $\zeta^2 \neq \alpha$, 
then $\alpha \zeta^{-1} \neq \zeta$. This implies the first statement.

Assume now that~$\mu$ is irreducible.
Then the roots of~$\mu$ are $\zeta^{q^i}$, $0 \leq i \leq 
d - 1$. If $d = 1$, then $\zeta \in \mathbb{F}_q^*$ and $\zeta = 
\alpha \zeta^{-1}$, i.e.\ $\alpha = \zeta^2$ is a square in $\mathbb{F}_q$.
Thus~$\mu$ is as in Case~(I).
Now suppose that $d > 1$. If $\alpha = \zeta^2$, then $d = 2$, $q$ is odd,
$\mu = \mu'$, and $\mu$ is as in~(II).

Suppose then that $\zeta^2 \neq \alpha$. Then, by the first statement,
$d = 2e$ is even and $\zeta^{q^e} = \alpha \zeta^{-1}$ by an easy argument.
Suppose that~$q$ is odd and $\mu = \mu'$. A similar argument as above shows 
that $\zeta^{q^e} = -\zeta$. We obtain $\zeta^2 = -\alpha \in \mathbb{F}_q$. 
Thus $d = 2$ and $\alpha \zeta^{-1} = \zeta^q = -\zeta$, i.e.~$\mu$ is as 
in~(III). If~$q$ is even or if~$q$ is odd and $\mu \neq \mu'$, then~$\mu$ is 
as in~(IV).

(b) This is an elementary exercise.
%%%%%%%%%%%%%%%%%%%%%%%%%%%%%%%%%%%%%%%%%%%%%%%%%%%%%%%%%%%%%%%%%%%%%%%%%%%%%%%%
%%%%%%%%%%%%%%%%%%%%%%%%%%%%%%%%%%%%%%%%%%%%%%%%%%%%%%%%%%%%%%%%%%%%%%%%%%%%%%%%
%%
%% Zum Beweis von (a) siehe die Notizen vom 28.09.2016, Seiten 1 - 2.
%% Für einen Beweis von (b) siehe Notizen vom 11.03.2016, Seiten 1 - 2.
%%
%%%%%%%%%%%%%%%%%%%%%%%%%%%%%%%%%%%%%%%%%%%%%%%%%%%%%%%%%%%%%%%%%%%%%%%%%%%%%%%%
%%%%%%%%%%%%%%%%%%%%%%%%%%%%%%%%%%%%%%%%%%%%%%%%%%%%%%%%%%%%%%%%%%%%%%%%%%%%%%%%
\end{prf}

%\addtocounter{subsection}{1}
\subsubsection{Linear transformations}
\label{RemarksOnTransformations}
Let $\mathbf{V}$ be a finite dimensional vector space over~$\mathbb{F}$, and
let~$s$ be a semisimple element of $\GL( \mathbf{V} )$.
If $\mu \in \mathbb{F}[X]^0$ and $U \subseteq \mathbf{V}$, we put
$U_\mu(s) := \ker( \mu( s ) ) \cap U$. Then $\mathbf{V}_\mu(s)$ is invariant
under the centralizer of~$s$ in $\GL( \mathbf{V} )$.
More generally, suppose that $h \in \GL( \mathbf{V} )$ such that 
$hsh^{-1} = \alpha s$ for some $\alpha \in \mathbb{F}$. Then 
$\mathbf{V}_\mu(s) h = \mathbf{V}_{\mu^\alpha}(s)$.

%\addtocounter{subsection}{2}
\subsection{The general unitary groups}
\label{GeneralUnitaryGroups}
In this subsection we introduce the general unitary groups and investigate some
of their semisimple elements.

\subsubsection{The groups}
\label{SpecialUnitaryGroupsI}
Let $\mathbf{G} := \GL_n( \mathbb{F} )$, acting on the natural vector space 
$\mathbf{V} := \mathbb{F}^n$ on the right. 
For a matrix $b = (\beta_{ij}) \in \mathbb{F}^{d \times e}$ we write $\bar{b} 
:= ( \beta_{ij}^q )$. 
We define the Frobenius morphism~$F$ on $\mathbf{G}$ by
\begin{equation}
\label{UnitaryF}
F(a) := J \left( \bar{a}^T \right)^{-1} J^{-1},\quad\quad 
a \in \mathbf{G}.
\end{equation}
Then $G = \mathbf{G}^F = \GU_n(q) \leq \GL_n(q^2)$.
Indeed, the sesqui-linear form 
\begin{equation}
\label{SesquiLinear}
\mathbf{V} \times \mathbf{V} \rightarrow
\mathbf{V}, (v,w) \mapsto vJ{\bar{w}}^T \text{\  for\ } v,w \in \mathbf{V}
\end{equation}
is non-degenerate, its restriction to $V := (\mathbb{F}_{q^2})^n$ is hermitian
and~$G$ is the unitary group of this hermitian form, acting naturally on~$V$.

For a semisimple element $s \in G$ we define $\mathcal{F}_{s} 
\subseteq \mathbb{F}_{q^2}[X] \cap \mathbb{F}[X]^0$ as the set of monic 
irreducible factors of the minimal polynomial of~$s$ (viewed as a linear 
transformation on $V = (\mathbb{F}_{q^2})^n$). Let $\mu, \nu \in 
\mathcal{F}_{s}$. Then $\mathbf{V}_\mu(s) \leq \mathbf{V}_\nu(s)^\perp$ if 
and only if $\nu \neq \mu^\dagger$. In particular, $\mathbf{V}_\mu(s)$ is 
totally isotropic if and only if $\mu \neq \mu^\dagger$. In the latter case, 
$\mathbf{V}_\mu(s) \oplus \mathbf{V}_{\mu^\dagger}(s)$ is non-degenerate. If 
$\mu = \mu^\dagger$, then $\mathbf{V}_\mu(s)$ is non-degenerate. Analogous 
statements hold for the finite spaces $V_\mu( s ) = \mathbf{V}_\mu(s) \cap V$.
Finally, as~$s$ is conjugate 
to $({\bar{s}^T})^{-1}$, we have $\dim \mathbf{V}_\mu(s) = 
\dim \mathbf{V}_{\mu^\dagger}(s)$ for all $\mu \in \mathcal{F}_{s}$.
Let~$\mathbf{L}$ denote an $F$-stable Levi subgroup of~$\mathbf{G}$ 
containing~$s$. We then put
$$\tilde{C}_{\mathbf{L}}( s ) := \{ g \in \mathbf{L} \mid gsg^{-1} = \gamma s
\text{\ for some\ } \gamma \in \mathbb{F}^* \}$$
and
$$\tilde{A}_{\mathbf{L}}( s ) := 
\tilde{C}_{\mathbf{L}}( s )/C_{\mathbf{L}}( s ).$$
Then $\tilde{C}_{\mathbf{L}}( s )$ is $F$-stable and thus there is a natural
action of~$F$ as an endomorphism on $\tilde{A}_{\mathbf{L}}( s )$. Moreover,
$\tilde{A}_{\mathbf{L}}( s )^F = 
\tilde{C}_{\mathbf{L}}( s )^F/C_{\mathbf{L}}( s )^F$, as $C_{\mathbf{L}}( s )$
is connected.

For $\zeta_1, \ldots , \zeta_n \in \mathbb{F}^*$, let $h(\zeta_1, \zeta_2, 
\ldots, \zeta_n) \in \mathbf{G}$ denote the diagonal matrix with entry~$\zeta_i$
at the $(i,i)$-position, and let
$$\mathbf{T} := \{ h(\zeta_1, \ldots , \zeta_n) \mid \zeta_i \in \mathbb{F}^*
\text{\ for all\ } 1 \leq i \leq n \}.$$
Then~$\mathbf{T}$ is a maximally split torus of~$\mathbf{G}$ and the Weyl group
$W := N_{\mathbf{G}}( \mathbf{T} )/\mathbf{T}$ may and will be identified with
the set of permutation matrices in~$\mathbf{G}$.

\subsubsection{Semisimple elements}
We begin by constructing certain semisimple elements and determine
their centralizers.
%\addtocounter{thm}{1}
\begin{lem}
\label{CentralizersInUnitaryGroups}
Let $\mu \in \mathbb{F}_{q^2}[X] \cap \mathbb{F}[X]^0$ be irreducible 
over~$\mathbb{F}_{q^2}$ and of degree~$d$. If $\mu = \mu^\dagger$, put 
$\tilde{\mu} := \mu$, otherwise, put $\tilde{\mu} := \mu\mu^\dagger$. 
Assume that $n = \mbox{\rm deg}({\tilde{\mu}})k$ for some integer~$k$.

Then there exists a semisimple element $s \in G$ with characteristic
polynomial~$\tilde{\mu}^k$.

Let~$s \in G$ be semisimple with characteristic 
polynomial~$\tilde{\mu}^k$. Put $C := C_{\mathbf{G}}( s )^F$. Then one of the 
following cases occurs.

{\rm (u)} If $\mu = \mu^\dagger$, then~$d$ is odd and  we have 
$C \cong \GU_k( q^{d} )$.

{\rm (l)} If $\mu \neq \mu^\dagger$, we have $C \cong \GL_k( q^{2d} )$.
\end{lem}
\begin{prf}
Let $t := h( \zeta_1 , \ldots , \zeta_n) \in \mathbf{T}$, where $\{ \zeta_1, 
\ldots , \zeta_n \}$ is the multiset of zeros of~$\tilde{\mu}^k$. Then~$F(t) 
= h(\zeta_n^{-q}, \ldots , \zeta_1^{-q})$ and~$t$ are conjugate by an element 
of~$W$ as $\tilde{\mu} = \tilde{\mu}^\dagger$. Thus there is an $F$-stable 
conjugate~$s$ of~$t$ (see~\ref{SemisimpleElementsAndCentralizers}).

Now let~$s$ be a semisimple, $F$-stable element with characteristic 
polynomial~$\tilde{\mu}^k$. Then~$s$ is conjugate to~$t$ in~$\mathbf{G}$. 
Let us sketch the arguments to
derive the isomorphism type of~$C$ and the remaining claims in case $\mu = 
\mu^\dagger$. Let $\zeta$ be a root of~$\mu$. Then $\zeta^{-q} 
= \zeta^{q^{2i}}$ for some $0 \leq i < d$. As~$d$ is the smallest positive
integer with $\zeta^{q^{2d}} = \zeta$, it follows that~$d$ is odd. 
Next, let $F'$ denote the Frobenius endomorphism of~$\mathbf{G}$ defined by
$F'(a) = {(\bar{a}^T)}^{-1}$ for $a \in \mathbf{G}$, so that $F(a) = 
w_0F'(a)w_0^{-1}$, where $w_0 \in W$ denotes the longest element of~$W$.
Then, if $F'(t)^w = t$ for some $w \in W$, we have $C_{\mathbf{G}}(s)^F \cong 
C_{\mathbf{G}}(t)^{Fw_0w} = C_{\mathbf{G}}(t)^{F'w}$ 
(see~\ref{SemisimpleElementsAndCentralizers}). Next, we may assume that
$$t = h( \zeta, \ldots, \zeta, \zeta^{q^2}, \ldots , \zeta^{q^2}, \ldots,
\zeta^{q^{2d-2}}, \ldots, \zeta^{q^{2d-2}}),$$
where each eigenvalue of~$t$ occurs exactly~$k$ times. For $1 \leq i \leq d$,
let $\mathbf{V}_i \leq \mathbf{V}$ denote the eigenspace of~$t$ for the 
eigenvalue $\zeta^{q^{2i-2}}$. Then 
$$C_{\mathbf{G}}( t ) = \GL( \mathbf{V}_1 ) \times \cdots \times 
\GL( \mathbf{V}_d ),$$
where $\GL( \mathbf{V}_i )$ is viewed as a subgroup of~$\mathbf{G}$ in the
natural way, $i = 1, \ldots , d$.
Let~$c$ denote the $d$-cycle on $\{ 1, \ldots , d \}$ such that
$(\zeta^{q^{2i - 2}})^{-q} = \zeta^{q^{2c(i) - 2}}$ for $1 \leq i \leq d$.
Let $w \in W$ be the corresponding block permutation so that $F'(t)^w = t$. 
Then
$$C_{\mathbf{G}}( t )^{F'w} = \{ (x_1, \ldots, x_d ) \mid
x_i \in \GL( \mathbf{V}_i ), F'(x_i) = x_{c(i)}, 1 \leq i \leq d \}.$$
As~$c$ is a $d$-cycle, and as~$d$ is odd, if follows that 
$$C_{\mathbf{G}}( t )^{F'w} \cong 
\{ x \in \GL( \mathbf{V}_1 ) \mid (F')^d(x) = x \} \cong \GU_k( q^d ).$$
\end{prf}

%\addtocounter{thm}{1}

We next consider the case of a general semisimple element in~$G$.
\begin{lem}
\label{CriticalCaseInGUn}
Let $s \in G$ be semisimple. Then there is $\alpha \in \mathbb{F}^*$
with $\alpha^n = 1 = \alpha^{q+1}$, and an isomorphism
\begin{equation}
\label{AlphaIsomorphism}
\langle \alpha \rangle \rightarrow \tilde{A}_{\mathbf{G}}( s )^F,
\end{equation}
such that the following conditions hold.

{\rm (a)} For all $\mu \in \mathcal{F}_{s}$ we have $\mu^\alpha \in 
\mathcal{F}_{s}$, $\mu^\dagger \in \mathcal{F}_{s}$ and
${(\mu^\dagger)}^\alpha = {(\mu^\alpha)}^\dagger$. In particular,
$\tilde{A}_{\mathbf{G}}( s )^F$ acts on $\mathcal{F}_{s}$ through the
isomorphism~{\rm (\ref{AlphaIsomorphism})}, and this action is equivalent
to the action of
$\tilde{A}_{\mathbf{G}}( s )^F$ on $\{ \mathbf{V}_\mu( s ) \mid 
\mu \in \mathcal{F}_{s} \}$.

{\rm (b)} Let $\mu \in \mathcal{F}_{s}$ be of degree~$d$ and let $\mathcal{O}$ 
denote the $\langle \alpha \rangle$-orbit of~$\mu$, and $\tilde{\mathcal{O}}$ the
union of~$\mathcal{O}$ with the $\langle \alpha \rangle$-orbit of~$\mu^\dagger$. 
Let $\mathbf{V}_1 := \sum_{\nu \in \tilde{\mathcal{O}}} \mathbf{V}_\nu(s)$. 
Then~$\mathbf{V}_1$ is non-degenerate with respect to the form defined 
in~{\rm (\ref{SesquiLinear})}. Write~$s_1$ for the element of~$\mathbf{G}$ that 
acts as~$s$ on~$\mathbf{V}_1$, and as the identity on its orthogonal complement. 
There is a natural embedding of $\GL( \mathbf{V}_1 )$ into~$\mathbf{G}$ such that 
$\GL( \mathbf{V}_1 )$ is $F$-stable, $s_1 \in \GL( \mathbf{V}_1 )$, and 
$F(s_1) = s_1$. Put $C := C_{\GL( \mathbf{V}_1 )}( s_1 )^F$.

If $\mu = \mu^\dagger$, we define $e := |\mathcal{O}|$. If $\mu \neq \mu^\dagger$,
then $|\tilde{\mathcal{O}}|$ is even, and we put $e := |\tilde{\mathcal{O}}|/2$. 

{\rm (u)} If $\mu = \mu^\dagger$, then~$d$ is odd and  we have
\begin{equation}
\label{Caseu}
C \cong \GU_k( q^{d} ) \times \cdots \times \GU_k( q^{d} )\quad\quad( e \text{ factors}),
\end{equation}
and $\tilde{A}_{\mathbf{G}}( s )^F$ acts on~$C$ by transitively permuting the
factors $\GU_k( q^{d} )$.

{\rm (l)} If $\mu \neq \mu^\dagger$, we have
\begin{equation}
\label{Casel}
C \cong \GL_k( q^{2d} ) \times \cdots \times \GL_k( q^{2d} )\quad\quad( e \text{ factors}),
\end{equation}
and $\tilde{A}_{\mathbf{G}}( s )^F$ acts on~$C$ by transitively
permuting the factors $\GL_k( q^{2d} )$. Moreover, the following two subcases
occur.
\begin{itemize}
\item[{\rm (ls)}] The $\langle \alpha \rangle$-orbit of~$\mu$ does not
contain~$\mu^\dagger$.
\item[{\rm (lt)}] The $\langle \alpha \rangle$-orbit of~$\mu$
contains~$\mu^\dagger$. In this case,~$q$ is odd.
\end{itemize}
\end{lem}
\begin{prf}
The map $\tilde{C}_{\mathbf{G}}( s ) \rightarrow \{ \gamma I_n \mid 
\gamma \in \mathbb{F}^* \} \leq \mathbf{G}$, $h \mapsto [h,s]$ is a group 
homomorphism with kernel $C_{\mathbf{G}}( s )$. If $\gamma I_n$ lies in the 
image of this map, then~$s$ and $\gamma s$ have the same determinant, hence 
$\gamma^n = 1$. It follows that $\tilde{A}_{\mathbf{G}}( s )$ is cyclic of
order dividing~$n$. Choose $g \in \tilde{C}_{\mathbf{G}}( s )^F$ such
that $\bar{g} := gC_{\mathbf{G}}(s)^F$ is a generator of 
$\tilde{A}_{\mathbf{G}}( s )^F$, and suppose that $[g,s] = \alpha I_n$.
Then, as $F(g) = g$, we have $\alpha^{-q} = \alpha$, i.e.\ $\alpha^{q+1} = 1$.
Moreover, $\langle \alpha \rangle \rightarrow \tilde{A}_{\mathbf{G}}( s )^F$,
$\alpha^i \mapsto \bar{g}^i$, $i \in \mathbb{Z}$, is an isomorphism.

(a) As $F( s ) = s$, we have $\mu^\dagger \in 
\mathcal{F}_{s}$ for all $\mu \in \mathcal{F}_{s}$.
As $\alpha s = g s g^{-1}$, the minimal
polynomial of $\alpha s$ is the same as that of~$s$, and thus
$\mu^{\alpha} \in \mathcal{F}_{s}$ for all $\mu \in 
\mathcal{F}_{s}$. From $\alpha^{-q} = \alpha$, it follows that
${(\mu^\dagger)}^\alpha = {(\mu^\alpha)}^\dagger$ for all
$\mu \in \mathcal{F}_{s}$. Finally,
$\mathbf{V}_\mu( s )g = \mathbf{V}_{\mu^\alpha}( s )$
for all $\mu \in \mathcal{F}_{s}$, and all parts of~(a) are proved.

(b) Clearly,~$\mathbf{V}_1$ is non-degenerate, and its orthogonal complement 
equals $\mathbf{V}_2 := 
\sum_{\nu \in (\mathcal{F}_s \setminus \tilde{\mathcal{O}})} \mathbf{V}_\nu(s)$
(see the second paragraph in~\ref{SpecialUnitaryGroupsI}).
Let $\mathbf{H} \leq \mathbf{G}$ denote the stabilizer of~$\mathbf{V}_1$ 
and~$\mathbf{V}_2$. Then~$\mathbf{H}$ factors as~$\mathbf{H}
= \mathbf{H}_1 \times \mathbf{H}_2$, where~$\mathbf{H}_i$ is naturally 
isomorphic to $\GL( \mathbf{V} )_i$, $i = 1, 2$. Now~$\mathbf{H}_i$ is
$F$-stable, as the setwise stabilizer and the pointwise stabilizer 
in~$\mathbf{G}$ of~$\mathbf{V}_i$ are $F$-stable, $i = 1, 2$. Clearly, 
$s_1 \in \mathbf{H}_1$, and writing $s = s_1 s_2$ with $s_2 \in \mathbf{H}_2$, 
we obtain $F(s_i) = s_i$ from the fact that~$s$ and $\mathbf{H}_i$ are 
$F$-stable, $i = 1, 2$.
By choosing an appropriate basis for~$\mathbf{V}_1$, the action of~$F$ on 
$\GL( \mathbf{V}_1 )$ is as in~(\ref{UnitaryF}). We may thus assume that 
$\mathbf{V} = \mathbf{V}_1$.

If $\mu = \mu^\dagger$, the claims on~$d$ and~$C$ follow by applying 
Lemma~\ref{CentralizersInUnitaryGroups} to each element of~$\mathcal{O}$. 
Suppose that $\mu \neq \mu^\dagger$. Then $\nu \neq \nu^\dagger \in 
\tilde{\mathcal{O}}$ for each $\nu \in 
\tilde{\mathcal{O}}$ and hence $|\tilde{\mathcal{O}}|$ is even. If 
$\mu^\dagger \in \mathcal{O}$, then $\tilde{\mathcal{O}} = \mathcal{O}$,
and as $|\mathcal{O}|$ divides the order of~$\alpha$, the latter is even as
well. Thus~$q$ is odd, since $\alpha^{q+1} = 1$. The other statements follow
from Lemma~\ref{CentralizersInUnitaryGroups} applied to $\nu\nu^\dagger$
for $\nu \in \mathcal{O}$.
\end{prf}

It is not hard to see that all three cases in part~(b) of the above lemma occur
for suitable~$q$.

\subsection{Conformal groups} 
\label{ConformalGroupsSection}
Here, we introduce the conformal groups of symplectic and quadratic forms and
investigate certain of their semisimple elements and their centralizers.

\subsubsection{Forms and groups} \label{ConformalGroups} 
Assume that $\mbox{\rm dim}( \mathbf{V} ) = 2m$ with $m \geq 1$. 
Let~$\mathbf{V}$ be equipped with a non-degenerate symplectic form~$\beta$ 
or a non-degenerate quadratic form~$Q$. There is a basis
\begin{equation}
\label{StandardBasis}
v_1, \ldots, v_m, v_m', \ldots, v_1'
\end{equation}
of~$\mathbf{V}$ such that if $(a_1, \ldots, a_m,a_m', \ldots , a_1'), 
(b_1, \ldots, b_m,b_m', \ldots , b_1') \in \mathbb{F}^n$ are coordinates of
$v, w \in \mathbf{V}$ with respect to~(\ref{StandardBasis}), we have
\begin{equation}\label{Beta}
\beta(v,w) = \sum_{i=1}^m a_ib_i'-a_i'b_i
\end{equation}
and 
\begin{equation} \label{Q}
Q(v) = \sum_{i=1}^ma_ia_i'.
\end{equation}
If we identify~$\mathbf{V}$ with $\mathbb{F}^n$ using the 
basis~(\ref{StandardBasis}), we view~$\beta$ and~$Q$ as forms on~$\mathbb{F}^n$
defined by the formulae~(\ref{Beta}) respectively~(\ref{Q}).
Let~$\hat{\mathbf{G}}$ denote the conformal group of~$\mathbf{V}$ with respect
to the form~$\beta$ or~$Q$, and put $\mathbf{G} := \hat{\mathbf{G}}^\circ$. If 
the form is symplectic, we have $\mathbf{G} = \hat{\mathbf{G}}$, 
whereas~$\mathbf{G}$ has index~$2$ in $\hat{\mathbf{G}}$, if the form is 
quadratic. In this case we call~$\mathbf{G}$ the special conformal group of~$Q$.
Usually we identify~$\hat{\mathbf{G}}$ and~$\mathbf{G}$ with their groups of 
matrices with respect to the basis~(\ref{StandardBasis}), so that 
$\hat{\mathbf{G}} = \mathbf{G} = \CSp_{2m}( \mathbb{F} )$ or $\hat{\mathbf{G}} 
= \CO_{2m}( \mathbb{F} )$ and $\mathbf{G} = \CSO_{2m}( \mathbb{F} )$, where
$\CSp_n( \mathbb{F} )$ denotes the group of all $g \in \GL_n( \mathbb{F} )$
such that $\beta(vg,wg) = \alpha_g \beta(v,w)$ for some $\alpha_g \in 
\mathbb{F}^*$ and all $v,w \in \mathbb{F}^n$. Similarly, $\CO_{n}( \mathbb{F} )$ 
is the group of all $g \in \GL_n( \mathbb{F} )$ such that 
$Q(vg) = \alpha_g Q(v)$ for some $\alpha_g \in \mathbb{F}^*$ and all
$v \in \mathbb{F}^n$. If~$q$ is odd, $\CSO_{n}( \mathbb{F} ) = 
\{ g \in \CO_n( \mathbb{F} ) \mid \det(g) = \alpha_g^m \}$ (whereas 
$\CSO_{n}( \mathbb{F} )$ simply denotes the connected component of
$\CO_{2m}( \mathbb{F} )$ if~$q$ is even). 
The element~$\alpha_g$ in the
above description is called the multiplier of~$g$. The map $\mathbf{G}
\rightarrow \mathbb{F}^*$, $g \mapsto \alpha_g$ is a surjective
homomorphism of groups with kernel $[\mathbf{G},\mathbf{G}]$. We have
$[\mathbf{G},\mathbf{G}] = \Sp_{2m}( \mathbb{F} )$, respectively
$[\mathbf{G},\mathbf{G}] = \SO_{2m}( \mathbb{F} )$ (where 
$\SO_{2m}( \mathbb{F} )$ is defined to be the connected component of the full 
orthogonal group $\GO_{2m}( \mathbb{F} )$ if~$q$ is even). If~$q$ is even, we
have $\mathbf{G} = Z(\mathbf{G}) \times [\mathbf{G},\mathbf{G}]$ as
algebraic groups, and in view of our intended applications we could as 
well just work with~$[\mathbf{G},\mathbf{G}]$. We have chosen our approach 
for the sake of a uniform treatment.

Let~$F'$ denote the Frobenius map on~$\GL( \mathbf{V} )$ which raises 
every matrix entry (with respect to the basis~(\ref{StandardBasis})) to its
$q$th power. Then $\mathbf{G}^{F'} = \CSp_{2m}( q )$ or $\CSO^+_{2m}( q )$, 
respectively. For $1 \leq i \leq m$ we let $\sigma_i$ denote the 
automorphism of~$\mathbf{V}$ which swaps the basis vectors $v_i$ and $v_i'$
(and fixes the other basis vectors). The Frobenius morphism $F''$ of
$\mathbf{G} = \CSO_{2m}( \mathbb{F} )$ is defined by $F''( g ) := 
\sigma_m^{-1}F'(g)\sigma_m$ for $g \in \mathbf{G}$. Then $\mathbf{G}^{F''}
= \CSO^-_{2m}( q )$. 
The pairs $(\CSO_{2m}( \mathbb{F} ), F' )$ and 
$(\CSO_{2m}( \mathbb{F} ), F'' )$ or the corresponding finite groups
are also called orthogonal groups of plus-type and minus-type, respectively.
In the following we let~$F$ be one of~$F'$ or~$F''$, where $F = F''$ 
implicitly assumes $\mathbf{G} = \CSO_{2m}( \mathbb{F} )$.
By~$V$ we denote the natural $\mathbb{F}_q$-vector space for~$\mathbf{G}^F$
inside~$\mathbf{V}$. This is $V = \langle v_1, \ldots , v_m, v_m', \ldots , 
v_1' \rangle_{\mathbb{F}_q}$ if $F = F'$, and $V = \{ \sum_{i=1}^m a_iv_i + 
a_i' v_i' \mid a_i, a_i' \in \mathbb{F}_q, 1 \leq i \leq m - 1, a_m \in 
\mathbb{F}_{q^2}, a_m' = a_m^q \}$ if $F = F''$.

\subsubsection{A maximally split torus and the Weyl group}
\label{MaximalllySplitTorus}
Let~$\mathbf{T}$ denote the torus of diagonal matrices (with respect to the
basis~(\ref{StandardBasis})) in~$\mathbf{G}$. We write
$h( \zeta_1, \ldots, \zeta_m; \alpha )$ for the element of~$\GL( \mathbf{V} )$
which acts by multiplication with $\zeta_i$ on~$v_i$, and by multiplication
with $\alpha\zeta_i^{-1}$ on~$v_i'$, $1 \leq i \leq m$. Then
$h( \zeta_1, \ldots, \zeta_m; \alpha ) \in \mathbf{T}$ has
multiplier~$\alpha$. Moreover,
$$\mathbf{T} = \{ h( \zeta_1, \ldots, \zeta_m; \alpha ) \mid
\zeta_i, \alpha \in \mathbb{F}^*, 1 \leq i \leq m \}.$$
Let $W := N_{\mathbf{G}}( \mathbf{T} )/\mathbf{T}$ denote the Weyl group
of~$\mathbf{G}$ with respect to~$\mathbf{T}$. We will identify~$W$ with a
group of permutations of the basis $\{ v_i, v_i' \mid 1 \leq i \leq m \}$
as follows. First, the $\sigma_i$ defined in~\ref{ConformalGroups} are viewed 
as permutations of this basis. Next, for $1 \leq i < m$ let~$\tau_i$ denote 
the double transposition $(v_i, v_{i+1})(v_i', v_{i+1}')$. Define
$$\hat{W} := \langle \tau_1, \ldots , \tau_{m-1}, 
\sigma_m \rangle \leq \mbox{\rm Sym}( \{ v_i, v_i' \mid 1 \leq i \leq m \} ).$$
Notice that $\sigma_i \in \hat{W}$ for all $1 \leq i \leq m$. If~$\mathbf{G}$ is 
symplectic, $\hat{W} = W$. If~$\mathbf{G}$ is orthogonal,~$W$ is of index~$2$ 
in~$\hat{W}$, namely
$$W = \langle \tau_1, \ldots , \tau_{m-1}, 
\sigma_m^{-1} \tau_{m-1} \sigma_m \rangle.$$
Notice that a product of $\sigma_i$ is contained in~$W$ if and only if it 
consists of an even number of factors.
In either case,~$\hat{W}$ is a Coxeter group of type~$B_m$ which naturally acts
on~$\mathbf{T}$; in the orthogonal case,~$W$ is a Coxeter group of type~$D_m$
(which is trivial if~$m = 1$). In the latter case, we also have $\hat{W} = 
N_{\hat{\mathbf{G}}}( \mathbf{T} )/{\mathbf{T}}$ and $\hat{\mathbf{G}} = 
\langle \mathbf{G}, \dot{w} \rangle$ for any $w \in \hat{W} \setminus W$.

Suppose that $t \in \mathbf{T}$ such that $C_{\hat{W}}( t ) \leq W$ and let
$w \in \hat{W} \setminus W$ with $F(t)^w = t$. Then there is no element in~$W$
conjugating~$t$ to~$F(t)$. Indeed, if $w' \in \hat{W}$ is such that $t^{w'} = 
F(t)$, then $w'w \in C_{\hat{W}}( t ) \leq W$, and thus $w' \in \hat{W} 
\setminus W$. This argument will be used frequently in the sequel.

\subsubsection{Semisimple elements}\label{SemisimpleInConformal}
Let $s \in \mathbf{G}$ be semisimple with multiplier~$\alpha$. Then~$s$ is
conjugate in~$\mathbf{G}$ to an element $h(\zeta_1, \ldots, \zeta_m;\alpha) 
\in \mathbf{T}$. In particular, the multiplicities of $\zeta, \alpha \zeta^{-1}
\in \mathbb{F}^*$ as eigenvalues of~$s$ are the same. If $\zeta^2 = \alpha$,
then $\alpha\zeta^{-1} = \zeta$, and hence~$\zeta$ occurs with even multiplicity
as eigenvalue of~$s$. More generally, let~$\mu$ be a monic factor of the minimal 
polynomial of~$s$. Then each root of~$\mu$ occurs with the same multiplicity in 
the characteristic polynomial of~$s$ as the corresponding root 
of~$\mu^{*\alpha}$. In particular, $\mu^{*\alpha}$ also divides the minimal 
polynomial of~$s$.  Moreover, if~$\mu$ and~$\mu^{*\alpha}$ are relatively prime, 
then $\mathbf{V}_\mu( s )$ is totally isotropic, and $\mathbf{V}_\mu( s )
\oplus \mathbf{V}_{\mu^{*\alpha}}( s )$ is non-degenerate. On the other hand, if 
$\mu = \mu^{*\alpha}$, then $\mathbf{V}_\mu( s )$ is non-degenerate. Again, 
corresponding statements also hold for the finite vector spaces $V_\mu(s) =
\mathbf{V}_\mu(s) \cap V$. (As in \cite{HiHuMa} we call a subset~$U$ 
of~$\mathbf{V}$ \textit{totally isotropic}, if the form vanishes on~$U \times U$
respectively~$U$. In case of a quadratic form, a set $U \subseteq \mathbf{V}$ 
with this property is sometimes called \textit{totally singular}.)

Suppose that $s \in \mathbf{G}$. We then put
$$\tilde{C}_{\mathbf{G}}( s ) := \{ g \in \mathbf{G} \mid g s g^{-1} = \pm s \}$$
and
$$\tilde{A}_{\mathbf{G}}( s ) := \tilde{C}_{\mathbf{G}}( s )/C^\circ_{\mathbf{G}}( s ).$$

\subsubsection{Connectedness of centralizers}\label{ConnectednessOfCentralizers}
Notice that centralizers of semi\-simple elements in~$\mathbf{G}$ are connected 
if $\mathbf{G} = \CSp_{2m}( \mathbb{F} )$ or if $\mathbf{G} = 
\CSO_{2m}( \mathbb{F} )$ and $q$ is even: In the first case, 
$[\mathbf{G}, \mathbf{G}] = \Sp_{2m}( \mathbb{F} )$ is simply connected, and in
the second case $Z( \mathbf{G}^* )$ is connected.
We next recall the classification of the semisimple elements with a 
non-connected centralizer in case $\mathbf{G} = \CSO_{2m}( \mathbb{F} )$
and~$q$ odd.

%\addtocounter{thm}{2}
\begin{lem}
\label{NonConnectedCentralizersInCSO}
Let~$q$ be odd.
Suppose that~$\mathbf{G} = \CSO_{2m}( \mathbb{F} )$ and let $s \in \mathbf{G}$ 
be semisimple with multiplier~$\alpha$. Let~$\zeta \in \mathbb{F}^*$ with 
$\zeta^2 = \alpha$. If~$\zeta$ is not an eigenvalue of~$s$, then 
$C_{\mathbf{G}}( s )$ and $C_{[\mathbf{G},\mathbf{G}]}( s )$ are connected.

Suppose 
that~$\zeta$ and~$-\zeta$ are eigenvalues of~$s$. Then 
$|A_{\mathbf{G}}( s )| = 2$ and there is $h \in \mathbf{G}$ 
stabilizing $\mathbf{V}_{X -\zeta}( s ) \oplus \mathbf{V}_{X + \zeta}( s )$ 
and acting trivially on its orthogonal complement, such that
$C_{\mathbf{G}}( s ) = \langle C^\circ_{\mathbf{G}}( s ), h \rangle$.
\end{lem}
\begin{prf}
We have $\mathbf{G} = Z(\mathbf{G})[\mathbf{G},\mathbf{G}]$ and thus 
$C_{\mathbf{G}}( s ) = Z(\mathbf{G})C_{[\mathbf{G},\mathbf{G}]}( s )$. It 
follows that $C_{\mathbf{G}}( s )$ is connected if 
$C_{[\mathbf{G},\mathbf{G}]}( s )$ is, as $Z(\mathbf{G})$ is a torus.
We may assume that $s  = h(\zeta_1, \ldots , \zeta_m; \alpha) \in \mathbf{T}$. 
We have an orthogonal decomposition
$$\mathbf{V} = \mathbf{V}_{X - \zeta}( s )  \oplus 
\mathbf{V}_{X + \zeta}( s ) \oplus \left( \bigoplus_{\xi} 
\tilde{\mathbf{V}}_\xi \right),$$ 
with $\tilde{\mathbf{V}}_\xi := \mathbf{V}_{X - \xi}( s ) \oplus 
\mathbf{V}_{X - \alpha\xi^{-1}}( s )$, 
where~$\xi$ runs through a suitable subset of $\{ \zeta_1, \ldots, \zeta_m \} 
\setminus \{ \zeta, -\zeta \}$. Now $C_{\mathbf{G}}( s )$ fixes all these
spaces, and $C_{[\mathbf{G},\mathbf{G}]}( s )$ induces the full linear group 
$\GL( \mathbf{V}_{X - \xi}( s ) )$ on each of $\tilde{\mathbf{V}}_\xi$.

We may thus assume that $\mathbf{V} = \mathbf{V}_{X - \zeta}( s )  \oplus
\mathbf{V}_{X + \zeta}( s )$. Then $C_{[\mathbf{G},\mathbf{G}]}( s )$ is 
connected if $\mathbf{V}_{X - \zeta}( s ) = 0$.
If none of the above two spaces vanishes, then
$$C_{\mathbf{G}}( s ) = \langle \CSO( \mathbf{V}_{X - \zeta}( s ) ) \times 
\CSO( \mathbf{V}_{X + \zeta}( s ) ), h \rangle,$$
for an element $h \in \mathbf{G}$ of order~$2$ which simultaneously swaps the
elements of a hyperbolic pair of~$\mathbf{V}_{X - \zeta}( s )$ and of one of 
$\mathbf{V}_{X + \zeta} ( s )$.
Thus $C_{\mathbf{G}}( s )$ is not connected, and $A_{\mathbf{G}}( s )$
is of order~$2$, generated by the image of~$h$.
\end{prf}

%\addtocounter{thm}{1}

\subsubsection{Special semisimple elements}
Here we investigate certain critical semisimple elements and their centralizers.
In the following, we will fix an irreducible, monic polynomial~$\mu \in 
\mathbb{F}_q[X]$ of degree~$d$ with non-zero constant term, a root 
$\zeta \in \mathbb{F}^*$ of~$\mu$ and an element $\alpha \in \mathbb{F}_q^*$. We 
will distinguish various types of~$\mu$. If $\mu = \mu^{*\alpha}$, we label the 
types of~$\mu$ from (I) to (IV) as in Lemma~\ref{MuLemma}(a). We say that~$\mu$ 
has Type~(V), if $\mu \neq \mu^{*\alpha}$. Recall that Types~(II) and~(III) only
occur if~$q$ is odd.

\begin{lem}
\label{SpecialMinimalPolynomial}
Let~$\mu \in \mathbb{F}_q[X]$ and $\alpha$ be as above. If $\mu = 
\mu^{*\alpha}$, put $\bar{\mu} := \mu$. Otherwise, $\bar{\mu} := 
\mu\mu^{*\alpha}$. Notice that Lemma~{\rm \ref{MuLemma}(a)} implies that 
$\deg(\bar{\mu})$ is even, unless~$\mu$ has Type~{\rm (I)}. Suppose that 
$2m = n = \deg(\bar{\mu})k$ for some integer~$k$.

Then there exists a semisimple element $s \in G$ with 
multiplier~$\alpha$ and characteristic polynomial $\bar{\mu}^k$ if and only 
if $(\mathbf{G},F)$ and~$k$ satisfy the conditions displayed in the following 
table. The centralizer $C^\circ_{[\mathbf{G},\mathbf{G}]}(s)^F$ of such an
element~$s$ is as given in the table.
$$
\begin{array}{ccccc} \\ \hline\hline
\mbox{\rm Type} & m &(\mathbf{G},F) & k & C^\circ_{[\mathbf{G},\mathbf{G}]}(s)^F 
\rule[- 7pt]{0pt}{ 20pt} \\ \hline\hline
\mbox{\rm I} & k/2 &\begin{array}{c} (\CSp_{2m}( \mathbb{F} ), F' ) \\
                                      (\CSO_{2m}( \mathbb{F} ), F') \\ 
                                      ( \CSO_{2m}( \mathbb{F} ), F'') \end{array} 
             & \mbox{\rm even} 
             & \begin{array}{c} \Sp_{2m}(q) \\ 
                                \SO_{2m}^+( q ) \\ 
                                \SO_{2m}^-(q) \end{array} \rule[ 0pt]{0pt}{ 26pt} \\ \hline
\mbox{\rm II} & k & \begin{array}{c} (\CSp_{2m}( \mathbb{F} ), F' ) \\
                                      (\CSO_{2m}( \mathbb{F} ), F') \\ 
                                      ( \CSO_{2m}( \mathbb{F} ), F'') \end{array} 
              & \mbox{\rm even} 
              & \begin{array}{c} \Sp_{m}(q^2) \\
                                 \SO_{m}^+( q^2 ) \\
                                 \SO_{m}^-(q^2) \end{array} \rule[ 0pt]{0pt}{ 26pt} \\ \hline
\mbox{\rm III, IV} & ke & \begin{array}{c} (\CSp_{2m}( \mathbb{F} ), F' ) \\
                                      (\CSO_{2m}( \mathbb{F} ), F') \\ 
                                      ( \CSO_{2m}( \mathbb{F} ), F'') \end{array} &
                     \begin{array}{c} \text{\rm\ any} \\ 
                                      \mbox{\rm\ even} \\ 
                                      \mbox{\rm\ odd} \end{array}
                                         & \GU_k(q^e) \rule[ 0pt]{0pt}{ 26pt} \\ \hline
\mbox{\rm V} & kd & (\mathbf{G},F') & \mbox{\rm any} & \GL_k(q^d) \rule[- 0pt]{0pt}{13pt} \\ \hline\hline
\end{array}
$$

\bigskip

\noindent
In particular, there is no such element if $\mu \neq \mu^{*\alpha}$ and
$(\mathbf{G},F) = (\CSO_{2m}( \mathbb{F} ), F'')$.
\end{lem}
\begin{prf}
The proof is given separately for each type of~$\mu$. 
If~$\mu$ has Type~(I), the assertion is trivially satisfied.

Suppose that~$q$ is odd and that~$\mu$ has Type~(II). Assume that~$G$ contains 
an element with characteristic polynomial $\mu^k$. As already observed in 
\ref{SemisimpleInConformal}, the square roots $\zeta$ and $-\zeta$ 
of~$\alpha$ occur with even multiplicity in~$\mu^k$ and thus~$k$ is even.
Conversely, if~$k$ is even, the element
$$
t := h( \zeta, -\zeta, \zeta, -\zeta, \ldots, \zeta, -\zeta; \alpha ) \in 
\mathbf{T}
$$
is conjugate to an element of~$G$ with multiplier~$\alpha$ and characteristic 
polynomial~$\mu^k$. Indeed, as $\zeta^q = -\zeta$, we have
$F(t)^w = t$ for $w = \tau_1 \tau_3 \cdots \tau_{k-1} \in W$. 

Suppose that~$\mu$ has Type~(IV) or that~$q$ is odd and~$\mu$ has Type~(III). By 
Lemma~\ref{MuLemma}(a), the degree of~$\mu$ is even, say $d = 2e$, and the roots 
of~$\mu$ are 
$\zeta, \zeta^q, \ldots, \zeta^{q^{e-1}}, \alpha\zeta^{-q^{e-1}}, 
\ldots, \alpha\zeta^{-q}, \alpha\zeta^{-1}$
with $\zeta^{q^e} = \alpha\zeta^{-1}$.
Hence
$$
\label{ElementInTofTypeIV}
t := h(\zeta, \ldots, \zeta, \zeta^q, \ldots, \zeta^q, \ldots, \zeta^{q^{e-1}},
\ldots, \zeta^{q^{e-1}}; \alpha),
$$
where each eigenvalue of~$t$ occurs exactly~$k$ times, has characteristic 
polynomial~$\mu^k$. Every element of~$\mathbf{T}$ with characteristic 
polynomial~$\mu^k$ is conjugate in~$W$ to~$t$ or to~$t^{\sigma_m}$. In
particular, there is $s \in G$ with characteristic 
polynomial~$\mu^k$, if and only if~$s$ is conjugate in~$\mathbf{G}$ to~$t$ or 
to~$t^{\sigma_m}$. As $\xi \neq \alpha\xi^{-1}$ for all roots~$\xi$ of~$\mu$, 
we have $C_{\hat{W}}( t ) \leq W$. 
Clearly, $F(t)^w = t$ for an element $w \in \hat{W}$ of the form 
$w = \sigma_{(e-1)k + 1} \cdots \sigma_{m'} \cdot \tau$ for some $\tau \in 
\langle \tau_{1} , \ldots , \tau_{m-1} \rangle$, where $m' = m$ if $F = F'$, 
and $m' = m-1$ if $F = F''$. Thus $w \in W$ if and only if $(\mathbf{G},F,k)$ 
are as in the assertion. The same argument works for~$t^{\sigma_m}$.

Finally suppose that~$\mu$ has Type~(V). If such an element~$s$ exists, we have 
$\mathbf{V} = \mathbf{V}_\mu(s) \oplus \mathbf{V}_{\mu^{*\alpha}}(s)$, a direct 
sum of two totally isotropic subspaces. Thus $(\mathbf{G},F) \neq 
(\CSO_{2m}( \mathbb{F} ), F'')$. Hence assume that $F = F'$ in the following.
The roots of~$\mu$ are $\zeta, \zeta^q, \ldots, \zeta^{q^{d-1}}$, and the roots 
of $\mu^{*\alpha}$ are
$\alpha \zeta^{-q^{d - 1}}, \alpha \zeta^{-q^{d-2}}, \ldots, \alpha \zeta^{-1}$.
Thus
\begin{equation}
\label{ElementInTofTypeV}
t := h( \zeta, \ldots, \zeta, \zeta^q, \ldots, \zeta^q, \ldots ,
\zeta^{q^{d-1}}, \ldots , \zeta^{q^{d-1}};\alpha), 
\end{equation}
where each eigenvalue of~$t$
occurs exactly~$k$ times, has characteristic polynomial~$(\mu\mu^{*\alpha})^k$.
Since $F = F'$, the element $F(t)$ is conjugate to~$t$ by an element $w \in 
\langle \tau_{1} , \ldots , \tau_{m-1} \rangle \leq W$, proving the existence 
of~$s$ as claimed. 

To prove the claims on the structure of the centralizers, we employ the method 
introduced in~\ref{SemisimpleElementsAndCentralizers}. For each $(t,w)$ as above
with $F(t)^w = t$, choose $u \in \mathbf{G}$ with $F(u)u^{-1} = \dot{w}$; then
$s := t^u$ is $F$-stable. To determine the structure of 
$C^\circ_{[\mathbf{G},\mathbf{G}]}( s )^F$, we have to compute 
$C^\circ_{[\mathbf{G},\mathbf{G}]}( t )^{F\dot{w}}$, which amounts to a routine 
calculation. We omit the details.
\end{prf}

\subsubsection{Critical semisimple elements} \label{CriticalSemisimpleElements}
In the investigations below we will use the following notation. Let $s \in G$ be 
semisimple and let~$\nu$ be a monic factor of the minimal polynomial of~$s$ 
with $\nu = \nu^{*\alpha}$, so that $\mathbf{V}_\nu(s)$ is a non-degenerate 
subspace of~$\mathbf{V}$. We then write $\hat{\mathbf{G}}_{\nu}(s)$ for the set 
of restrictions to $\mathbf{V}_{\nu}(s)$ of the elements of the setwise 
stabilizer of $\mathbf{V}_{{\nu}}(s)$ in~$\hat{\mathbf{G}}$.
Then~$\hat{\mathbf{G}}_{\nu}(s)$ is the conformal group on $\mathbf{V}_{\nu}(s)$
with respect to the restricted symplectic, respectively quadratic form on
$\mathbf{V}_{\nu}(s)$. If~$\kappa$ is another monic factor of the minimal 
polynomial of~$s$ with $\kappa = \kappa^{*\alpha}$ and $\nu$ and~$\kappa$ are
relatively prime, then $\hat{\mathbf{G}}_{\nu\kappa}(s) \cong 
\hat{\mathbf{G}}_\nu(s) \times \hat{\mathbf{G}}_\kappa(s)$, and if $\nu\kappa$
equals the minimal polynomial of~$s$, then $\hat{\mathbf{G}}_\nu(s) \times 
\hat{\mathbf{G}}_\kappa(s)$ can naturally be identified with the stabilizer
of~$\mathbf{V}_\nu(s)$ and $\mathbf{V}_\kappa(s)$ in~$\hat{\mathbf{G}}$.
The subgroup of $\hat{\mathbf{G}}_{\nu}(s)$ constituting
the special conformal group on $\mathbf{V}_{\nu}(s)$ is denoted by
$\mathbf{G}_{\mu}(s)$ (by convention, $\hat{\mathbf{G}}_{\nu}(s) = 
\mathbf{G}_{\nu}(s)$ if $\mathbf{G} = \CSp_{2m}( \mathbb{F} )$).
Finally,~$s_{\nu}$ denotes the restriction of~$s$ to $\mathbf{V}_{\nu}(s)$, so
that $s_{\nu} \in \hat{\mathbf{G}}_{\nu}(s)$. As the setwise stabilizer and the
pointwise stabilizer of $\mathbf{V}_{\nu}(s)$ in~$\hat{\mathbf{G}}$ are
$F$-stable, $F$~induces a well-defined Frobenius morphism on
$\hat{\mathbf{G}}_{\nu}(s)$, which we also denote by~$F$. Then
$s_{\nu} \in \hat{\mathbf{G}}_{\nu}(s)^F$.

We now consider centralizers of critical elements $s \in G$ and the action 
of $\tilde{A}_{\mathbf{G}}(s)$ on these in case~$q$ is odd (recall the 
definition of $\tilde{A}_{\mathbf{G}}(s)$ in {\rm \ref{SemisimpleInConformal}}). 
The main results here are given in 
Tables~{\rm \ref{ConjugateToNegativeSymplecticSpecial}}
and~{\rm \ref{ConjugateToNegativeOrthogonalSpecial}}, where we adopt the 
following notation. The cyclic group of order~$2$
is denoted by~$C_2$. We write $\text{\rm ``f''}$ and $\text{\rm ``g''}$ for a
field and a graph automorphism, respectively. Also, $\text{\rm ``sg''}$ denotes
an automorphism which simultaneously acts as a graph automorphism on each of
the two factors of $C^\circ_{[\mathbf{G},\mathbf{G}]}( s )^F$. (If $\mathbf{G}
= \CSO_{2m}( \mathbb{F} )$ and $\mathbf{H} \leq \mathbf{G}$ is $F$-stable and
invariant under conjugation by~$\sigma_m$, the restriction of this conjugation
to~$\mathbf{H}^F$ is called a graph automorphism of~$\mathbf{H}^F$.) We write
$\text{\rm ``flip''}$ for the automorphism which swaps the two factors of
$C^\circ_{[\mathbf{G},\mathbf{G}]}( s )^F$ (more precisely, an element $xy \in
C^\circ_{[\mathbf{G},\mathbf{G}]}( s )^F$, where~$x$ and~$y$ lie in distinct
direct factors of the factorization of $C^\circ_{[\mathbf{G},\mathbf{G}]}( s )^F$
given in Tables~\ref{ConjugateToNegativeSymplecticSpecial} 
and~\ref{ConjugateToNegativeOrthogonalSpecial}, is mapped to $yx$). 
The indices in $\text{\rm I}.1$, $\text{\rm I}.2$,
$\text{\rm II}.1$ and $\text{\rm II}.2$
distinguish the two $\CSO_{2m}(\mathbb{F})^F$-conjugacy classes with
characteristic polynomial $(X^2 - \alpha)^m$.

\begin{lem}
\label{ElementsInConformalGroups}
Assume that~$q$ is odd.
Suppose that
$\mu = \mu^{*\alpha}$ or $\mu' = \mu^{*\alpha}$.  If $\mu = \mu'$, put 
$\tilde{\mu} := \mu$, otherwise put $\tilde{\mu} := \mu\mu'$. Assume that 
$2m = n = \mbox{\rm deg}(\tilde{\mu})k$ for 
some integer~$k$. 

{\rm (a)} There exists a semisimple $s \in G$ with 
multiplier~$\alpha$ and characteristic polynomial $\tilde{\mu}^k$, if and 
only if $F$ and~$m$ satisfy the conditions displayed in 
Tables~{\rm \ref{ConjugateToNegativeSymplecticSpecial}}
and~{\rm \ref{ConjugateToNegativeOrthogonalSpecial}}.

{\rm (b)} Suppose that $s \in G$ is semisimple with 
multiplier~$\alpha$ and characteristic polynomial $\tilde{\mu}^k$.
Then $C_{[\mathbf{G},\mathbf{G}]}( s )$ is connected, unless
$\mu$ has Type~{\rm (I)} or~{\rm (II)}.
The structure of the groups $C^\circ_{[\mathbf{G},\mathbf{G}]}( s )^F$
and $\tilde{A}_{\mathbf{G}}(s)^F$, as well as the action of the latter 
on the former, are as displayed in
Tables~{\rm \ref{ConjugateToNegativeSymplecticSpecial}}
and~{\rm \ref{ConjugateToNegativeOrthogonalSpecial}}. 

{\rm (c)} Let~$s$ be as in~{\rm (b)} and suppose that $\mathbf{G} = 
\CSO_{2m}( \mathbb{F} )$. Let $S \in \hat{\mathbf{G}}$ denote the set of 
elements conjugating~$s$ to~$-s$. Then~$S$ is $F$-stable and non-empty.
Moreover, $S^F$ is non-empty, unless $\mu$ has Type~{\rm (I)} or~{\rm (II)} 
and $F = F''$. If~$S^F$ is non-empty, then~$S^F$ contains an element with 
multiplier~$1$. Table~{\rm \ref{ConjugatingSToMinusSInCSO}} displays the 
location of~$S$ in $\hat{\mathbf{G}}$.
\end{lem}

\begin{table}[t]
\caption{\label{ConjugateToNegativeSymplecticSpecial} The critical elements in
$\CSp_{2m}(q)$, $q$ odd (explanations in the paragraph preceding 
Lemma~\ref{ElementsInConformalGroups})}
$$
\begin{array}{ccccc}\\ \hline\hline
\text{Type} & m & C_{[\mathbf{G},\mathbf{G}]}( s )^F & 
\tilde{A}_{\mathbf{G}}( s )^F & \text{Action}
\rule[- 7pt]{0pt}{ 22pt} \\ \hline\hline
\text{I} & \text{\ even} & \Sp_{m}( q ) \times \Sp_{m}( q ) & 
C_2 & \text{flip} \rule[ 0pt]{0pt}{ 13pt} \\ \hline
\text{II} & \text{\ even} & \Sp_{m}( q^2 ) & C_2 & \text{\ f}
\rule[ 0pt]{0pt}{ 13pt} \\ \hline
\text{III} & \text{\ any} & \GU_{m}( q ) & C_2 & \text{\ g}
\rule[ 0pt]{0pt}{ 13pt} \\ \hline
\text{IV} & 2ke & \GU_{k}( q^{e} ) \times \GU_{k}( q^{e} ) & C_2 & 
\text{\ flip}
\rule[- 0 pt]{0pt}{ 13pt} \\ \hline
\text{V} & kd & \GL_{k}( q^{d} ) & 
C_2 & \text{\ g}
\rule[- 5pt]{0pt}{18pt} \\ \hline\hline
\end{array}
$$
\end{table}

\begin{table}[t]
\caption{\label{ConjugateToNegativeOrthogonalSpecial} The critical elements in
$\CSO^\pm_{2m}(q)$, $q$ odd (explanations in the paragraph preceding 
Lemma~\ref{ElementsInConformalGroups})}
$$
\begin{array}{cccccc}\\ \hline\hline
\text{Type} & F & m & C^\circ_{[\mathbf{G},\mathbf{G}]}( s )^F & 
\tilde{A}_{\mathbf{G}}( s )^F & \text{Action}
\rule[- 7pt]{0pt}{ 22pt} \\ \hline\hline
\text{I}.1 & \begin{array}{c} F' \\ F'' \end{array} & \text{\ even} & 
           \begin{array}{c} \SO^+_{m}( q ) \times \SO^+_{m}( q ) \\ 
           \SO^+_{m}( q ) \times \SO^-_{m}( q ) \end{array} & 
           \begin{array}{c} C_2 \times C_2 \\ C_2 \end{array} & 
           \begin{array}{c} \text{flip}\times\text{sg} \\
           \text{\ sg} \end{array} \rule[ 0pt]{0pt}{ 18pt} \\ \hline
\text{I}.2 & \begin{array}{c} F' \\ F'' \end{array} & \text{\ even} & 
           \begin{array}{c} \SO^-_{m}( q ) \times \SO^-_{m}( q ) \\ 
           \SO^-_{m}( q ) \times \SO^+_{m}( q ) \end{array} & 
           \begin{array}{c} C_2 \times C_2 \\ C_2 \end{array} & 
           \begin{array}{c} \text{flip}\times\text{sg} \\
           \text{\ sg} \end{array} \rule[ 0pt]{0pt}{ 18pt} \\ \hline
\text{II}.1,2 & \begin{array}{c} F' \\ F'' \end{array} & \text{\ even} & 
            \begin{array}{c} \SO^+_{m}( q^2 ) \\ \SO^-_m( q^2 ) \end{array} & 
            \begin{array}{c} C_2 \times C_2 \\ C_2 \end{array} & 
            \begin{array}{c} \text{f} \times\text{g} \\
           \text{\ g} \end{array} \rule[ 0pt]{0pt}{ 20pt} \\ \hline
\text{III} & \begin{array}{c} F' \\ F'' \end{array} & 
             \begin{array}{c} \text{\ even} \\ \text{\ odd} \end{array} & 
             \GU_{m}( q ) & C_2 & \text{\ g} \rule[ 0pt]{0pt}{ 20pt} \\ \hline
\text{IV} & F' & 2ke & \GU_{k}( q^{e} ) \times \GU_{k}( q^{e} ) & C_2 & 
\text{\ flip}
\rule[- 0 pt]{0pt}{ 13pt} \\ \hline
\text{V} & F' & kd & \GL_{k}( q^{d} ) & 
C_2 & \text{\ g}
\rule[- 5pt]{0pt}{18pt} \\ \hline\hline
\end{array}
$$
\end{table}

\begin{table}[t]
\caption{\label{ConjugatingSToMinusSInCSO} The location of the set~$S$ of
elements of~$\hat{\mathbf{G}} = \CO_{2m}( \mathbb{F} )$ conjugating~$s$ to~$-s$}
$$
\begin{array}{cccc}\\ \hline\hline
\text{Type} & F & m & S 
\rule[- 7pt]{0pt}{ 22pt} \\ \hline\hline
\text{I} & 
\begin{array}{c} F' \\ F'' \end{array} & \text{\ even} &
            \begin{array}{c} S \cap \mathbf{G}^F \neq \emptyset \neq 
                             S \cap ( \hat{\mathbf{G}}^F \setminus \mathbf{G}^F) \\ 
                             S \subset \mathbf{G} \setminus \hat{\mathbf{G}}^F 
            \end{array}
\rule[ 0pt]{0pt}{ 20pt} \\ \hline
\text{II} & 
\begin{array}{c} F' \\ F'' \end{array} & \text{\ even} &
            \begin{array}{c} S \subset \mathbf{G}, S^F \neq \emptyset \\ 
                             S \subset \mathbf{G} \setminus \hat{\mathbf{G}}^F 
            \end{array}
\rule[ 0pt]{0pt}{ 20pt} \\ \hline
\text{III} & \begin{array}{c} F' \\ F'' \end{array} & 
            \begin{array}{c} \text{\rm even} \\ \text{\rm odd} \end{array} &
            \begin{array}{c} S \subset \mathbf{G}, S^F \neq \emptyset \\ 
                             S \subset \hat{\mathbf{G}} \setminus \mathbf{G}, S^F \neq \emptyset
            \end{array}
\rule[ 0pt]{0pt}{ 20pt} \\ \hline
\text{IV} & F' & \text{\ any} & S \subset \mathbf{G}, S^F \neq \emptyset \rule[ 0pt]{0pt}{ 13pt} \\ \hline
\text{V} & F' & 
            \begin{array}{c} \text{\rm even} \\ \text{\rm odd} \end{array} &
            \begin{array}{c} S \subset \mathbf{G}, S^F \neq \emptyset \\ 
                             S \subset \hat{\mathbf{G}} \setminus \mathbf{G}, S^F \neq \emptyset
            \end{array}
\rule[- 4pt]{0pt}{25pt} \\ \hline\hline
\end{array}
$$
\end{table}

\begin{prf}
Let us begin with some preliminary remarks. Suppose that $s \in \mathbf{G}^F$ 
is semisimple with multiplier~$\alpha$ and characteristic polynomial 
$\tilde{\mu}^k$ and let $t \in \mathbf{T}$ be conjugate to~$s$ in~$\mathbf{G}$.
Then~$s$ is conjugate to~$-s$ in $\hat{\mathbf{G}}$, respectively~$\mathbf{G}$, 
if and only if~$t$ and~$-t$ lie in the same $\hat{W}$-orbit, respectively 
$W$-orbit, on~$\mathbf{T}$. Assume now that $\tilde{\mu} \neq X^2 - \alpha$. 
Then $C_{[\mathbf{G},\mathbf{G}]}( s )$ is connected by 
Lemma~\ref{NonConnectedCentralizersInCSO}. Moreover, $\Stab_{\hat{W}}( t ) 
\leq W$, as $\xi \neq \alpha\xi^{-1}$ for all eigenvalues~$\xi$ of~$t$.
This implies $C_{\hat{\mathbf{G}}}( s ) \leq C_{\mathbf{G}}( s )$.
Then if~$s$ is conjugate to~$-s$ by an element of~$\hat{\mathbf{G}}$ with 
multiplier~$1$, there is such a conjugating element which is $F$-stable. 
Moreover, either~$S$ is contained in~$\mathbf{G}$ or in $\hat{\mathbf{G}} 
\setminus \mathbf{G}$.
%%%%%%%%%%%%%%%%%%%%%%%%%%%%%%%%%%%%%%%%%%%%%%%%%%%%%%%%%%%%%%%%%%%%%%%%%%%%%%%%
%%%%%%%%%%%%%%%%%%%%%%%%%%%%%%%%%%%%%%%%%%%%%%%%%%%%%%%%%%%%%%%%%%%%%%%%%%%%%%%%
%%
%% Für einen Beweis dieser "preliminary remarks" siehe Notizen
%% vom 16.02.2016, Seiten 1 - 3.
%%
%%%%%%%%%%%%%%%%%%%%%%%%%%%%%%%%%%%%%%%%%%%%%%%%%%%%%%%%%%%%%%%%%%%%%%%%%%%%%%%%
%%%%%%%%%%%%%%%%%%%%%%%%%%%%%%%%%%%%%%%%%%%%%%%%%%%%%%%%%%%%%%%%%%%%%%%%%%%%%%%%

We prove (a), (b) and (c) simultaneously, distinguishing cases according to the 
type of~$\tilde{\mu}$. For each of these types, we give representatives of the
$W$-orbits of the elements of~$\mathbf{T}$ with characteristic polynomial
$\tilde{\mu}^k$. If~$t$ is such a representative, we give elements $w \in W$ 
with $F(t)^w = t$, one for each $G$-conjugacy class of $F$-stable 
elements conjugate to~$t$ in~$\mathbf{G}$. For $\tilde{\mu} = X^2 - \alpha$ 
and $\mathbf{G} = \CSO_{2m}( \mathbb{F} )$ there are two such classes if~$m$ 
is even; in all other cases, there is at most one such class. 
The structure of $C^\circ_{[\mathbf{G},\mathbf{G}]}( s )^F$ has already been
determined in Lemma~\ref{SpecialMinimalPolynomial}.
Finally, the action of 
$\tilde{A}_{\mathbf{G}}( s )^F$ on $C^\circ_{[\mathbf{G},\mathbf{G}]}( s )^F$
can be derived from the action of $\tilde{A}_{\mathbf{G}}( t )^{F\dot{w}} =
\tilde{C}_{\mathbf{G}}( t )^{F\dot{w}}/C^\circ_{\mathbf{G}}( t )^{F\dot{w}}$
on $C^\circ_{[\mathbf{G},\mathbf{G}]}( t )^{F\dot{w}}$.

We put $a := \tau_1 \tau_3 \cdots \tau_{m-1} \in W$
and $b := \sigma_{m-1} \sigma_m \in W$ (if $m \geq 2$). We also choose
commuting lifts $\dot{a}, \dot{b} \in N_{\mathbf{G}}( \mathbf{T} )$ of~$a$ 
and~$b$ with multipliers~$1$ and $F'( \dot{a} ) = \dot{a}$, $F( \dot{b} ) 
= \dot{b}$. Finally, we put $c := \sigma_1 \sigma_2 \cdots \sigma_m \in \hat{W}$, 
and choose a lift $\dot{c} \in N_{\hat{\mathbf{G}}}( \mathbf{T} )$ with 
multiplier~$1$ and $F( \dot{c} ) = \dot{c}$. If~$\mathbf{G} = 
\CSO_{2m}( \mathbb{F} )$, we choose~$\dot{a}$, $\dot{b}$ and $\dot{c}$ 
to be of order~$2$, and such that $F''( \dot{a} )\dot{a}^{-1} = \dot{b}$.

Suppose first that~$\mu$ has Type~(I) or~(II). Thus $\zeta \in \mathbb{F}_{q^2}$ 
and $\zeta^2 = \alpha$. As~$\zeta$ and~$-\zeta$ occur with even multiplicity as 
eigenvalues of any element of~$\mathbf{G}$ (see \ref{SemisimpleInConformal}), 
there is $s \in G^F$ with characteristic polynomial $\tilde{\mu}^m$ if and only 
if~$m$ is even. Suppose that~$m$ is even and let $s \in G$ have characteristic 
polynomial $\tilde{\mu}^m$. Then~$s$ is conjugate in~$\mathbf{G}$ 
to 
$$t = h( \zeta, -\zeta, \ldots, \zeta, -\zeta; \alpha) \in \mathbf{T}.$$
Notice that $t^a = -t$, $t^b = t$ and that $\langle a, b \rangle$ is a Klein 
four group. We have $F(t) = t$ if~$\zeta \in \mathbb{F}_q$, and $F(t)^a = t$,
if not.
Let $x, y, z \in \mathbf{G}$ with $F(x)x^{-1} = \dot{a}$, $F(y)y^{-1} = \dot{b}$
and $F(z)z^{-1} = \dot{a}\dot{b}$.
If~$\mathbf{G}$ is symplectic, $C_{\mathbf{G}}(t)$ is connected, and thus there
is a unique $G$-conjugacy classes in the $\mathbf{G}$-class of~$t$, 
with representative $s = t$ if~$\zeta \in \mathbb{F}_q$ and $s = t^x$, otherwise.
If~$\mathbf{G}$ is is orthogonal, $C_{\mathbf{G}}(t) = 
\langle C^\circ_{\mathbf{G}}(t), \dot{b} \rangle$ by 
Lemma~\ref{NonConnectedCentralizersInCSO}. Thus there are two 
$G$-conjugacy classes in the $\mathbf{G}$-class of~$t$, with
representatives $s \in \{ t, t^y \}$ if~$\zeta \in \mathbb{F}_q$, and 
$s \in \{ t^x, t^z \}$, otherwise. We thus may take ${w} \in 
\{ 1, {b} \}$ if $\zeta \in \mathbb{F}_q$, and ${w} \in 
\{ {a}, {a}{b} \}$, otherwise, and where ${w} = {b}$
and ${w} = {a}{b}$ only occur if $\mathbf{G} = \CSO_{2m}( \mathbb{F} )$.
We are done if~$\mathbf{G}$ is symplectic.
Suppose then that $\mathbf{G} = \CSO_{2m}( \mathbb{F} )$. 
Lemma~\ref{NonConnectedCentralizersInCSO} implies that
$\tilde{C}_{\mathbf{G}}( t ) = \langle C^\circ_{\mathbf{G}}( t ),
\dot{a}, \dot{b} \rangle$. If $F = F'$, then the elements of 
$\langle \dot{a}, \dot{b} \rangle$ are fixed by~$F'\dot{w}$ and thus
$\tilde{A}_{\mathbf{G}}( t )^{F'\dot{w}} \cong 
\langle \dot{a}, \dot{b} \rangle$. In particular there is an element of
$\mathbf{G}^{F'}$ with multiplier~$1$ conjugating~$s$ to~$-s$. 
If~$\mu$ has Type~(I), then $\dot{\sigma}_m \in 
C_{\hat{\mathbf{G}}}(t)^{F'\dot{w}}$, which means that there is an $F'$-stable 
conjugate of $\dot{\sigma}_m$ centralizing~$s$, and thus 
there is also an element in $\hat{\mathbf{G}}^{F'} \setminus \mathbf{G}^{F'}$ 
with multiplier~$1$ conjugating~$s$ to~$-s$. 
If~$\mu$
has Type~(II), then $C_{\hat{\mathbf{G}}}(t)^{F'\dot{w}} \leq \mathbf{G}^{F'}$.
%%%%%%%%%%%%%%%%%%%%%%%%%%%%%%%%%%%%%%%%%%%%%%%%%%%%%%%%%%%%%%%%%%%%%%%%%%%%%%%%
%%%%%%%%%%%%%%%%%%%%%%%%%%%%%%%%%%%%%%%%%%%%%%%%%%%%%%%%%%%%%%%%%%%%%%%%%%%%%%%%
%%
%% Einen Beweis für die Aussage $C_{\hat{\mathbf{G}}}(t)^{F'\dot{w}} \leq \mathbf{G}^F$
%% findet sich auf Blatt 1 vom 15.03.2016
%%
%%%%%%%%%%%%%%%%%%%%%%%%%%%%%%%%%%%%%%%%%%%%%%%%%%%%%%%%%%%%%%%%%%%%%%%%%%%%%%%%
%%%%%%%%%%%%%%%%%%%%%%%%%%%%%%%%%%%%%%%%%%%%%%%%%%%%%%%%%%%%%%%%%%%%%%%%%%%%%%%%
Thus there is no element of $\hat{\mathbf{G}}^{F'} \setminus \mathbf{G}^{F'}$
conjugating~$s$ to~$-s$. If $F = F''$, then 
$\tilde{A}_{\mathbf{G}}( t )^{F''\dot{w}} = 
C_{\mathbf{G}}(t)^{F''\dot{w}}/C^\circ_{\mathbf{G}}(t)^{F''\dot{w}}$,
as $F''( \dot{a} )\dot{a}^{-1} = \dot{b} \not\in C^\circ_{\mathbf{G}}(t)$. 
In particular, there is no element in $\mathbf{G}^{F''}$ conjugating~$s$ to~$-s$.

%%%%%%%%%%%%%%%%%%%%%%%%%%%%%%%%%%%%%%%%%%%%%%%%%%%%%%%%%%%%%%%%%%%%%%%%%%%%%%%%
%%%%%%%%%%%%%%%%%%%%%%%%%%%%%%%%%%%%%%%%%%%%%%%%%%%%%%%%%%%%%%%%%%%%%%%%%%%%%%%%
%%
%% Mehr Details zu diesen Überlegungen finden sich in meinen Notizen
%% vom 25.02.2016, Seiten 1 - 4.
%%
%%%%%%%%%%%%%%%%%%%%%%%%%%%%%%%%%%%%%%%%%%%%%%%%%%%%%%%%%%%%%%%%%%%%%%%%%%%%%%%%
%%%%%%%%%%%%%%%%%%%%%%%%%%%%%%%%%%%%%%%%%%%%%%%%%%%%%%%%%%%%%%%%%%%%%%%%%%%%%%%%
%%%%%%%%%%%%%%%%%%%%%%%%%%%%%%%%%%%%%%%%%%%%%%%%%%%%%%%%%%%%%%%%%%%%%%%%%%%%%%%%
%%%%%%%%%%%%%%%%%%%%%%%%%%%%%%%%%%%%%%%%%%%%%%%%%%%%%%%%%%%%%%%%%%%%%%%%%%%%%%%%
%%
%% Für Details zum Beweis für die Typen (I) und (II) siehe Notizen
%% vom 17.02.2016, Seiten 1 - 5.
%%
%%%%%%%%%%%%%%%%%%%%%%%%%%%%%%%%%%%%%%%%%%%%%%%%%%%%%%%%%%%%%%%%%%%%%%%%%%%%%%%%
%%%%%%%%%%%%%%%%%%%%%%%%%%%%%%%%%%%%%%%%%%%%%%%%%%%%%%%%%%%%%%%%%%%%%%%%%%%%%%%%

Suppose now that~$\mu$ has Type~(III). As $\tilde{\mu} = \mu$ in this case, the 
assertions of~(a) are contained in Lemma~\ref{SpecialMinimalPolynomial}.
To prove part~(b), assume that~$m$ is even in case $(\mathbf{G},F) = 
( \CSO_{2m}( \mathbb{F} ), F' )$, and that $m$ is odd in case $(\mathbf{G},F) =
( \CSO_{2m}( \mathbb{F} ), F'' )$. We have $\alpha\zeta^{-1} = -\zeta = \zeta^q$.
Put
$$t_1 := h( \zeta, \ldots, \zeta; \alpha) \in \mathbf{T}.$$
Any $s \in G$ with characteristic polynomial $\mu^m$ is conjugate 
in~$\mathbf{G}$ to~$t_1$ 
or to $t_2 := t_1^{\sigma_m}$. Let $t \in \{ t_1, t_2 \}$. Then $F(t) = -t$
if $F = F'$, and $F(t) = -t^{\sigma_m}$ if $F = F''$. As $t^c = -t$, we have
$F(t)^w = t$ for $w = c$ or $w = \sigma_m \cdot c = \sigma_1 \sigma_2 \cdots 
\sigma_{m-1}$.
As $\dot{c} \in \hat{\mathbf{G}}$ 
has multiplier~$1$, we conclude that~$s$ and~$-s$ are conjugate by an element of 
$\hat{G}$ by the preliminary remarks. This element can be chosen to 
lie in $G$, unless $(\mathbf{G},F) = ( \CSO_{2m}( \mathbb{F} ), F'' )$. 
This completes the proof of all claims in this case.
%%%%%%%%%%%%%%%%%%%%%%%%%%%%%%%%%%%%%%%%%%%%%%%%%%%%%%%%%%%%%%%%%%%%%%%%%%%%%%%%
%%%%%%%%%%%%%%%%%%%%%%%%%%%%%%%%%%%%%%%%%%%%%%%%%%%%%%%%%%%%%%%%%%%%%%%%%%%%%%%%
%%
%% Für Details zum Beweis für Typ (III) siehe Notizen
%% vom 18.02.2016, Seiten 1 - 2.
%%
%%%%%%%%%%%%%%%%%%%%%%%%%%%%%%%%%%%%%%%%%%%%%%%%%%%%%%%%%%%%%%%%%%%%%%%%%%%%%%%%
%%%%%%%%%%%%%%%%%%%%%%%%%%%%%%%%%%%%%%%%%%%%%%%%%%%%%%%%%%%%%%%%%%%%%%%%%%%%%%%%

Suppose next that~$\mu$ has Type~(IV). Then $\mu = \mu^{*\alpha} \neq \mu'$.
Assume that there is a semisimple element of~$G$ with characteristic 
polynomial~$\tilde{\mu}^k$. Then $\mathbf{V} = \mathbf{V}_{\mu}(s) 
\oplus \mathbf{V}_{\mu'}(s)$ with non-degenerate, orthogonal $F$-stable
subspaces $\mathbf{V}_{\mu}(s)$ and $\mathbf{V}_{\mu'}(s)$ of~$\mathbf{V}$.
Let~$s_1$ denote the restriction of~$s$ to~$\mathbf{V}_{\mu}(s)$ and~$s_1'$
the restriction of~$s$ to~$\mathbf{V}_{\mu'}(s)$. Then the characteristic
polynomials of~$s_1$ and~$s_1'$ equal $\mu^k$ and $(\mu')^k$, respectively.
Now apply Lemma~\ref{SpecialMinimalPolynomial} to $\mathbf{V}_{\mu}( s )$ 
and $\mathbf{V}_{\mu'}(s)$. If $\mathbf{G} = \CSO_{2m}( \mathbb{F} )$, then
$(\mathbf{G}_{\mu}(s),F)$ and $(\mathbf{G}_{\mu'}(s),F)$ are 
both of plus-type or both of minus-type and thus $(\mathbf{G},F) =
(\CSO_{2m}( \mathbb{F} ),F')$ (recall the definition of $\mathbf{G}_{\mu}(s)$
at the beginning of~\ref{CriticalSemisimpleElements}). 
Suppose from now on that $(\mathbf{G},F) \in 
\{ (\CSp_{2m}( \mathbb{F} ),F'), (\CSO_{2m}( \mathbb{F} ),F') \}$. Again
by Lemma~\ref{SpecialMinimalPolynomial}, there is an element~$s$ as required.
Put
$$t_1 := h(\zeta, -\zeta, \ldots, \zeta, -\zeta, \ldots, \zeta^{q^{e-1}}, 
 -\zeta^{q^{e-1}}, \ldots , \zeta^{q^{e-1}}, -\zeta^{q^{e-1}}; 
\alpha) \in \mathbf{T},$$
where each eigenvalue of~$t$ occurs exactly~$k$ times. Then~$s$ is conjugate 
in~$\mathbf{G}$ to~$t_1$ or to $t_2 := t_1^{\sigma_m}$. Let 
$t \in \{t_1, t_2 \}$. Then $F(t)^w = t$ for some $w \in W$, similar to the 
one taken in the proof of Lemma~\ref{SpecialMinimalPolynomial}. Also 
$t^a = -t$. 
All assertions in this case
follow from our preliminary considerations.
%%%%%%%%%%%%%%%%%%%%%%%%%%%%%%%%%%%%%%%%%%%%%%%%%%%%%%%%%%%%%%%%%%%%%%%%%%%%%%%%
%%%%%%%%%%%%%%%%%%%%%%%%%%%%%%%%%%%%%%%%%%%%%%%%%%%%%%%%%%%%%%%%%%%%%%%%%%%%%%%%
%%
%% Für Details zum Beweis für Typ (III) siehe Notizen
%% vom 18.02.2016, Seiten 3 - 4.
%%
%%%%%%%%%%%%%%%%%%%%%%%%%%%%%%%%%%%%%%%%%%%%%%%%%%%%%%%%%%%%%%%%%%%%%%%%%%%%%%%%
%%%%%%%%%%%%%%%%%%%%%%%%%%%%%%%%%%%%%%%%%%%%%%%%%%%%%%%%%%%%%%%%%%%%%%%%%%%%%%%%

Suppose finally that~$\mu$ has Type~(V). Here, $\mu \neq \mu^{*\alpha} = \mu'$, 
and~(a) follows from Lemma~\ref{SpecialMinimalPolynomial}. 
Let $t \in \mathbf{T}$ be the element defined in~(\ref{ElementInTofTypeV}), 
and let $w \in W$ with $F(t)^w = t$.
As the sequence $(\alpha\zeta^{-q^{d-1}}, \ldots , \alpha\zeta^{-1})$ is a
permutation of the sequence $(-\zeta, \ldots , -\zeta^{q^{d-1}})$, there is
an element $\tau \in \langle \tau_1, \ldots , \tau_{m-1} \rangle$ such that
$\delta := \tau \cdot a$ conjugates~$t$ to~$-t$. Clearly,~$\delta$ has an 
$F'$-stable lift $\dot{\delta}$ to $N_{\mathbf{G}}( \mathbf{T} )$ with 
multiplier~$1$. Moreover,~$\delta$ lies in~$W$, if and only if $kd = m$ is 
even. 
Again we are done with the preliminary remarks.
%%%%%%%%%%%%%%%%%%%%%%%%%%%%%%%%%%%%%%%%%%%%%%%%%%%%%%%%%%%%%%%%%%%%%%%%%%%%%%%%
%%%%%%%%%%%%%%%%%%%%%%%%%%%%%%%%%%%%%%%%%%%%%%%%%%%%%%%%%%%%%%%%%%%%%%%%%%%%%%%%
%%
%% Für Details zum Beweis für Typ (V) siehe Notizen
%% vom 22.02.2016, Seite 1.
%%
%%%%%%%%%%%%%%%%%%%%%%%%%%%%%%%%%%%%%%%%%%%%%%%%%%%%%%%%%%%%%%%%%%%%%%%%%%%%%%%%
%%%%%%%%%%%%%%%%%%%%%%%%%%%%%%%%%%%%%%%%%%%%%%%%%%%%%%%%%%%%%%%%%%%%%%%%%%%%%%%%
\end{prf}

\subsubsection{General semisimple elements} \label{GeneralSemisimpleElements}
Here we consider the most general semisimple elements of~$G$ relevant for our 
investigation. Let $s \in G$ be semisimple with multiplier~$\alpha$. Then the 
minimal polynomial of~$s$ lies in $\mathbb{F}_q[X]$, and we denote 
by~$\mathcal{F}_s$ the set of irreducible monic factors in $\mathbb{F}_q[X]$ 
of this minimal polynomial. Assume now that~$q$ is odd and that every $\mu \in 
\mathcal{F}_s$ satisfies $\mu = \mu^{*\alpha}$ or $\mu' = \mu^{*\alpha}$. If 
$\mu = \mu'$, put $\tilde{\mu} := \mu$, otherwise put $\tilde{\mu} := \mu\mu'$. 
For $\mu \in \mathcal{F}_s$ write $d_\mu$ and $k_\mu(s)$ for the degree of~$\mu$ 
and the multiplicity of~$\mu$ as a factor of the characteristic polynomial 
of~$s$, respectively. 

%%%%%%%%%%%%%%%%%%%%%%%%%%%%%%%%%%%%%%%%%%%%%%%%%%%%%%%%%%%%%%%%%%%%%%%%%%%%%%%%
%%%%%%%%%%%%%%%%%%%%%%%%%%%%%%%%%%%%%%%%%%%%%%%%%%%%%%%%%%%%%%%%%%%%%%%%%%%%%%%%
%%
%% Für Details zum Beweis für die Einschraenkung des Frobenius Morphismus siehe 
%% Notizen vom 29.03.2016, Seiten 1 - 3.
%%
%%%%%%%%%%%%%%%%%%%%%%%%%%%%%%%%%%%%%%%%%%%%%%%%%%%%%%%%%%%%%%%%%%%%%%%%%%%%%%%%
%%%%%%%%%%%%%%%%%%%%%%%%%%%%%%%%%%%%%%%%%%%%%%%%%%%%%%%%%%%%%%%%%%%%%%%%%%%%%%%%
As $\mathbf{V}_{\tilde{\mu}}(s)$ is a non-degenerate subspace of~$\mathbf{V}$, 
the notation introduced at the beginning of~\ref{CriticalSemisimpleElements}
will be applied with $\nu = \tilde{\mu}$. As~$s$
is fixed in the following we put $\hat{\mathbf{G}}_{\tilde{\mu}} := 
\hat{\mathbf{G}}_{\tilde{\mu}}(s)$ and $\mathbf{G}_{\tilde{\mu}} := 
\mathbf{G}_{\tilde{\mu}}(s)$.
Notice that if $\mathbf{G} = \CSO_{2m}( \mathbb{F} )$, then 
$(\mathbf{G}_{\tilde{\mu}},F)$ may be of plus-type or of minus-type 
(see~\ref{ConformalGroups}).

Let 
$\mathcal{F}_s' \subseteq \mathcal{F}_s$ denote a set of representatives
of the partition $\{ \{ \mu, \mu' \} \mid \mu \in \mathcal{F}_s \}$.
Let $I_{\text{(IV)}}$ and $I_{\text{(V)}}$ be two disjoint index sets labelling 
the elements of $\mathcal{F}_s'$ of Type~(IV) and of Type~(V), respectively. Then 
for $j \in I_{\text{(IV)}}$ and $\mu = \mu_j$ we write $e_j := d_{\mu}/2$ and 
$\ell_j := k_{\mu}( s )$.  Similarly, for $j \in I_{\text{(V)}}$ and $\mu = 
\mu_j$ we write $d_j := d_{\mu}$ and $k_j := k_{\mu}( s )$.  

\begin{prop}
\label{CorollaryElementsInConformalGroups}
Let the notation and the assumptions be as above. In particular,~$q$ is odd.
Then the following assertions hold.

{\rm (a)} Suppose that $\mathbf{G} = \CSp_{2m}( \mathbb{F} )$. Then~$s$ is
conjugate to~$-s$ in~$G$, if and only if $k_\mu(s) = k_{\mu'}(s)$
for all $\mu \in \mathcal{F}_s$, and these multiplicities satisfy the conditions 
displayed in {\rm Table~\ref{ConjugateToNegativeSymplectic}}. (The part of this
table after the vertical line is only relevant for 
Remark~{\rm \ref{RemarkForOddDimensionalSpinGroups}}.)

If these conditions are satisfied, $\tilde{A}_{\mathbf{G}}(s)^F$ has order~$2$.

{\rm (b)} Suppose that $\mathbf{G} = \CSO_{2m}( \mathbb{F} )$.
If~$s$ is conjugate to~$-s$ in~$G$, then $k_\mu(s) = k_{\mu'}(s)$ 
for all $\mu \in \mathcal{F}_s$, and one of the following cases occurs.

{\rm (i)} If $F = F'$ and~$m$ is even, the conditions displayed in 
{\rm Table~\ref{ConjugateToNegativeOrthogonalA}} are satisfied.

{\rm (ii)} If $F = F'$ and~$m$ is odd, then $q \equiv 1\,(\text{mod\ } 4)$,
all elements of $\mathcal{F}_s$ are of Types~{\rm(I),~(IV)} or~{\rm(V)}, and
the conditions displayed in
{\rm Table~\ref{ConjugateToNegativeOrthogonalB}} are satisfied.

{\rm (iii)} If $F = F''$, then~$m$ is odd, $q \equiv 3\,(\text{mod\ } 4)$,
there are no elements in $\mathcal{F}_s$ of Type~{\rm (II)}, $X^2 + \alpha
\in \mathcal{F}_s$ and 
the conditions displayed in
{\rm Table~\ref{ConjugateToNegativeOrthogonalC}} are satisfied.
(The respective parts of Tables~{\rm \ref{ConjugateToNegativeOrthogonalA},
\ref{ConjugateToNegativeOrthogonalB}} 
and~{\rm \ref{ConjugateToNegativeOrthogonalC}}
after the vertical lines are only relevant for
Remark~{\rm \ref{RemarkForEvenDimensionalSpinGroups}}.)

Suppose that~$s$ is as in one of the cases~{\rm (i)},~{\rm (ii)} or~{\rm (iii)}.
Then~$s$ is conjugate to~$-s$ in~$G$. If~$m$ is odd, 
then~$\tilde{A}_{\mathbf{G}}(s)^F$ is cyclic of order~$4$, and if $m$ is even, 
then~$\tilde{A}_{\mathbf{G}}(s)^F$ is a Klein four group if $X^2 - \alpha$ 
divides the minimal polynomial of~$s$, and is of order~$2$, otherwise. 
\end{prop}
\begin{prf}
We only prove~(b). The proof of~(a) is similar but much simpler. Let,~$k$ 
and~$\ell$ denote the multiplicity of $X^2 - \alpha$, respectively 
$X^2 + \alpha$, in the characteristic polynomial of~$s$. 
Suppose that~$s$ is conjugate to~$-s$ in~$G$. Then $k_\mu(s) = k_{\mu'}(s)$
for all $\mu \in \mathcal{F}_s$. 
Now $\mathbf{V}$ is the orthogonal direct sum of the non-degenerate subspaces 
$\mathbf{V}_{\tilde{\mu}}(s)$, where~$\mu$ runs through $\mathcal{F}_s'$.
If $\tilde{\mu} \neq X^2 - \alpha$, then $s_{\tilde{\mu}} \in 
\mathbf{G}_{\tilde{\mu}} = \CSO( \mathbf{V}_{\tilde{\mu}}(s) )$, as the 
eigenvalues of~$s_{\tilde{\mu}}$ come in pairs $\xi, \alpha\xi^{-1}$ with 
$\xi \neq \alpha\xi^{-1}$. As $s \in \CSO( \mathbf{V} )$, we also have
$s_{\tilde{\mu}} \in \CSO( \mathbf{V}_{\tilde{\mu}}(s) )$ for $\tilde{\mu} = 
X^2 - \alpha$ (if $k > 0$). Each of the pairs 
$(\mathbf{V}_{\tilde{\mu}}(s), s_{\tilde{\mu}})$ thus satisfies the hypotheses 
of Lemma~\ref{ElementsInConformalGroups}. In particular,~$k$ is even.
Moreover, an element of~$G$ conjugating~$s$ to~$-s$ must fix the 
spaces $\mathbf{V}_{\tilde{\mu}}(s)$, and thus induces an element of
$\hat{\mathbf{G}}_{\tilde{\mu}}^F$ conjugating~$s_{\tilde{\mu}}$ 
to~$-s_{\tilde{\mu}}$. If $\mu \in \mathcal{F}_s'$ is not of Type~(III),
then $(\mathbf{G}_{\tilde{\mu}},F)$ is of plus-type by 
Lemma~\ref{ElementsInConformalGroups}(c). If $\mu = X^2 + \alpha \in 
\mathcal{F}_s'$, then $(\mathbf{G}_{\tilde{\mu}},F)$ is of plus-type or of 
minus-type if~$\ell$ is even or odd, respectively. In particular, 
$(\mathbf{G},F)$ is of plus-type if and only if~$\ell$ is even. We now 
consider the individual cases.

(i) Suppose that $F = F'$ and that~$m$ is even. Then 
$\sum_{j \in I_{\text{\rm (V)}}} {d_jk_j}$ is even, 
and thus all conditions of Table~\ref{ConjugateToNegativeOrthogonalA} are 
satisfied.

(ii) Now suppose that $F = F'$ and that~$m$ is odd. Then~$\ell$ is even and
$\sum_{j \in I_{\text{\rm (V)}}} {d_jk_j}$ is odd, and thus an odd number of 
the $d_jk_j, j \in I_{\text{\rm (V)}}$, are odd. Lemma~\ref{MuLemma}(b) now 
implies that $-\alpha$ is a square in $\mathbb{F}_q$, and thus $X^2 + \alpha 
\not\in \mathcal{F}_s$, as it is not irreducible. 
Lemma~\ref{ElementsInConformalGroups}(c) implies that if~$\mu$ is of Type~(V) 
and $d_\mu k_{\mu}(s)$ is odd, then $s_{\tilde{\mu}}$ is conjugate to 
$-s_{\tilde{\mu}}$ by an element of $\hat{\mathbf{G}}_{\tilde{\mu}}^F$, but not
by an element of $\mathbf{G}_{\tilde{\mu}}^F$. Thus there must exist $\mu \in 
\mathcal{F}'_s$ not of Type~(V), such that an element of 
$\hat{\mathbf{G}}_{\tilde{\mu}}^F \setminus \mathbf{G}_{\tilde{\mu}}^F$
conjugates $s_{\tilde{\mu}}$ to $-s_{\tilde{\mu}}$. 
By Lemma~\ref{ElementsInConformalGroups}(c), this~$\mu$ must have Type~(I).
In particular,~$\alpha$ is a square in~$\mathbb{F}_q$. Hence~$-1$
is a square in $\mathbb{F}_q$ and thus $q \equiv 1\,(\text{mod\ } 4)$.

(iii) Suppose that $F = F''$. Then~$\ell$ is odd, which implies that~$-\alpha$
is not a square in $\mathbb{F}_q$. It follows that
$\sum_{j \in I_{\text{\rm (V)}}} d_jk_j$ is even by Lemma~\ref{MuLemma}(b), and
thus~$m$ is odd. As in the proof of~(ii) we conclude that~$\mathcal{F}_s$
contains an element of Type~(I).
Thus $\alpha$ is a square in $\mathbb{F}_q$ and~$-1$ is not a square, hence
$q \equiv 3\,(\text{mod\ } 4)$.

Suppose now that, in the respective cases, the parameters are as in
Tables~\ref{ConjugateToNegativeOrthogonalA},~\ref{ConjugateToNegativeOrthogonalB}
and~\ref{ConjugateToNegativeOrthogonalC}.
Assume first that~$m$ is even. Then $F = F'$ and~$\ell$ and
$\sum_{j \in I_{\text{\rm (V)}}} d_jk_j$ are
even. By applying Lemma~\ref{ElementsInConformalGroups} to each of the pairs
$(\mathbf{V}_{\tilde{\mu}}(s), s_{\tilde{\mu}})$, we find that~$s$ is conjugate
to $-s$ in~$G$ and that ${\tilde{A}_{\mathbf{G}}( s )}^F$ is as asserted.
Next assume that~$m$ is odd. Then~$\alpha$ is a square
in~$\mathbb{F}_q$ and $\ell + \sum_j d_jk_j$ is odd. Moreover, $\ell = 0$ if
$F = F'$, and~$\ell$ is odd if $F = F''$. First consider the special
case that $m = k + \ell + \sum_{j \in I_{\text{\rm (V)}}} d_jk_j$, i.e.\ that there
are no $\mu \in \mathcal{F}_s$ of Type~(IV). Let $\zeta \in \mathbb{F}_q$
denote a root of $X^2 - \alpha$. Let $\zeta_{k + 1}, \ldots , \zeta_{m - \ell}$
be the roots of $\prod_{j \in I_{\text{\rm (V)}}} \mu_j^{k_j}$, counted with
multiplicities. If $\ell \neq 0$, i.e.\ if $F = F''$, let $\zeta_{k + \ell + 1}$
be a root of $X^2 + \alpha$, and put $\zeta_i := \zeta_{k + \ell + 1}$ for
$k + \ell + 2 \leq i \leq m$. Then~$s$ is conjugate in~$\mathbf{G}$ to
$$t_1 := h( \zeta, -\zeta, \ldots, \zeta, -\zeta, \zeta_{k+1}, \ldots , \zeta_{m}; \alpha)$$
with~$k/2$ occurrences of $\zeta$, or to $t_2 := t_1^{\sigma_m}$ if
$\ell \neq 0$. Let $t \in \{ t_1, t_2 \}$. Then there is
$w \in \langle \tau_{k+1}, \ldots , \tau_{m - \ell - 1} \rangle$ such that 
$$u := (\sigma_1 \cdot \tau_1 \tau_3 \cdots \tau_{k - 1}) \cdot (\sigma_{k + 1}
\cdots \sigma_{m - \ell} \cdot w) \cdot (\sigma_{m - \ell + 1} \cdots \sigma_{m})$$
conjugates $t$ to~$-t$.
Now $u \in W$ and lifts to an $F$-stable element~$\dot{u}$ of 
$N_{\mathbf{G}}( \mathbf{T} )$. Moreover, $\dot{u}^2 \not\in 
C_{\mathbf{G}}^\circ( t )$ by Lemma~\ref{NonConnectedCentralizersInCSO},
as 
$u^2 = \sigma_1\sigma_2\cdot(\sigma_{k + 1} \cdots \sigma_{m - \ell } \cdot w)^2$. 
This implies our claim in the special case. The general case follows from this 
with Lemma~\ref{ElementsInConformalGroups} applied to each $\mu \in 
\mathcal{F}_s'$ of Type~(IV).
\end{prf}

\begin{landscape}
\begin{table}[f]
\caption{\label{ConjugateToNegativeSymplectic} The critical elements in
$\CSp_{2m}(q)$, $q$ odd (explanations in the paragraph preceding 
Proposition~\ref{CorollaryElementsInConformalGroups}
and in Remark~\ref{RemarkForOddDimensionalSpinGroups})}
$$
\begin{array}{cccccc||cc}\\ \hline\hline
\text{Type} & 
k_\mu( {s} ) & 
d_\mu &
\mbox{\rm dim}( \mathbf{V}_{\tilde{\mu}}( {s} ) ) & 
\text{Condt's} & 
C_{[\mathbf{G}_{\tilde{\mu}},\mathbf{G}_{\tilde{\mu}}]}( {s}_{\tilde{\mu}} )^F &
\lambda & \lambda^a \rule[- 7pt]{0pt}{ 20pt} \\ \hline\hline
\text{I} & k & 1 & 2k & 
k \text{\ even}^1 &
\Sp_{k}( q ) \times \Sp_{k}( q ) & (\Lambda_1,\Lambda_2) & (\Lambda_2,\Lambda_1)
\rule[ 0pt]{0pt}{ 13pt} \\ \hline
\text{II} & k & 2 & 2k & 
k \text{\ even}^1 &
\Sp_{k}( q^2 ) & \Lambda & \Lambda
\rule[ 0pt]{0pt}{ 13pt} \\ \hline
\text{III} & \ell & 2 & 2\ell & 
 &
\GU_{\ell}( q ) & \kappa & \kappa
\rule[ 0pt]{0pt}{ 13pt} \\ \hline
\text{IV}.j, j \in I_{\text{\rm (IV)}} & \ell_{j} & 2e_j & 4 e_j \ell_{j} & 
 &
\GU_{\ell_{j}}( q^{e_j} ) \times \GU_{\ell_{j}}( q^{e_j} ) & 
(\kappa_{1,j},\kappa_{2,j}) & (\kappa_{2,j},\kappa_{1,j})
\rule[- 6 pt]{0pt}{ 18pt} \\ \hline
\text{V}.j, j \in I_{\text{\rm (V)}} & k_{j} & d_j & 2 d_j k_{j} & 
 &
\GL_{k_{j}}( q^{d_j} ) & \omega_j & \omega_j
\rule[- 7pt]{0pt}{22pt} \\ \hline\hline
\end{array}
$$
\bigskip
\noindent \begin{tabular}{ll}
1 & At most one of the multiplicities in Types I and II is non-zero. \\
\end{tabular}
\end{table}
\end{landscape}

\begin{landscape}
\begin{table}[f]
\caption{\label{ConjugateToNegativeOrthogonalA} The critical elements in
$\CSO^+_{2m}(q)$, $m$ even, $q$ odd (explanations in the paragraph preceding 
Proposition~\ref{CorollaryElementsInConformalGroups}
and in Remark~\ref{RemarkForEvenDimensionalSpinGroups})}
$$
\begin{array}{cccccc||ccc}\\ \hline\hline
\text{Type} & 
k_\mu(s) & 
d_\mu &
\mbox{\rm dim}( \mathbf{V}_{\tilde{\mu}}( s ) ) & 
\text{Condt's} & 
C^\circ_{[\mathbf{G}_{\tilde{\mu}},\mathbf{G}_{\tilde{\mu}}]}( s_{\tilde{\mu}} )^F
& \lambda & \lambda^a & \lambda^b \rule[- 7pt]{0pt}{ 20pt} \\ \hline\hline
\text{I}.1 & k & 1 & 2k &  k \text{\ even}^1 &
\SO_{k}^+( q ) \times \SO_{k}^+( q )  
& (\Lambda_1,\Lambda_2) & (\Lambda_2,\Lambda_1) & (\Lambda_1',\Lambda_2') 
\rule[ 0pt]{0pt}{ 19pt} \\ \hline
\text{I}.2 & k & 1 & 2k &  k  \text{\ even}^1 & 
\SO_{k}^-( q ) \times \SO_{k}^-( q ) 
& (\Lambda_1,\Lambda_2) & (\Lambda_2,\Lambda_1) & (\Lambda_1,\Lambda_2)
\rule[ 0pt]{0pt}{ 19pt} \\ \hline
\text{II}.1,2 & k & 2 & 2k & k \text{\ even}^1 &
\SO_{k}^+( q^2 ) 
& \Lambda & \Lambda & \Lambda'
\rule[ 0pt]{0pt}{ 19pt} \\ \hline
\text{III} & \ell & 2 & 2\ell & \ell \text{\ even} &
\GU_{\ell}( q ) & \kappa & \kappa & \kappa \rule[ 0pt]{0pt}{ 13pt} \\ \hline
\text{IV}.j, j \in I_{\text{\rm (IV)}} & \ell_{j} & 2e_j & 4 e_j \ell_{j} & &
\GU_{\ell_{j}}( q^{e_j} ) \times \GU_{\ell_{j}}( q^{e_j} ) 
& (\kappa_{1,j},\kappa_{2,j}) & (\kappa_{2,j},\kappa_{1,j}) & (\kappa_{1,j},\kappa_{2,j})
\rule[-6pt]{0pt}{ 18pt} \\ \hline
\text{V}.j, j \in I_{\text{\rm (V)}} & k_{j} & d_j & 2 d_j k_{j} & 
\sum_{j} d_jk_j \text{\rm\ even} & \GL_{k_{j}}( q^{d_j} ) 
& \omega_j & \omega_j & \omega_j \rule[- 7pt]{0pt}{22pt} \\ \hline\hline
\end{array}
$$
\bigskip
\noindent \begin{tabular}{ll}
1 & At most one of the multiplicities in Types I.1, I.2, II.1 and II.2 is non-zero. \\
\end{tabular}
\end{table}
\end{landscape}

\begin{landscape}
\begin{table}[f]
\caption{\label{ConjugateToNegativeOrthogonalB} The critical elements in
$\CSO^+_{2m}(q)$, $m$ odd, $q \equiv 1\,(\text{mod\ } 4)$
(explanations in the paragraph preceding
Proposition~\ref{CorollaryElementsInConformalGroups}
and in Remark~\ref{RemarkForEvenDimensionalSpinGroups})}
$$
\begin{array}{cccccc||cc}\\ \hline\hline
\text{Type} & 
k_\mu(s) & 
d_\mu &
\mbox{\rm dim}( \mathbf{V}_{\tilde{\mu}}( s ) ) & 
\text{Condt's} & 
C^\circ_{[\mathbf{G}_{\tilde{\mu}},\mathbf{G}_{\tilde{\mu}}]}( s_{\tilde{\mu}} )^F
& \lambda & \lambda^a \rule[- 7pt]{0pt}{ 20pt} \\ \hline\hline
\text{I}.1 & k & 1 & 2k & k \geq 2 \text{\ even}^1 &
\SO_{k}^+( q ) \times \SO_{k}^+( q ) & (\Lambda_1,\Lambda_2) & (\Lambda_2,\Lambda_1') 
\rule[ 0pt]{0pt}{ 19pt} \\ \hline
\text{I}.2 & k & 1 & 2k & k \geq 2 \text{\ even}^1 & 
\SO_{k}^-( q ) \times \SO_{k}^-( q ) & (\Lambda_1,\Lambda_2) & (\Lambda_2,\Lambda_1)
\rule[ 0pt]{0pt}{ 19pt} \\ \hline
\text{IV}.j, j \in I_{\text{\rm (IV)}} & \ell_{j} & 2e_j & 4 e_j \ell_{j} & 
& \GU_{\ell_{j}}( q^{e_j} ) \times \GU_{\ell_{j}}( q^{e_j} ) 
& (\kappa_{1,j},\kappa_{2,j}) & (\kappa_{2,j},\kappa_{1,j})
\rule[- 6 pt]{0pt}{ 18pt} \\ \hline
\text{V}.j, j \in I_{\text{\rm (V)}} & k_{j} & d_j & 2 d_j k_{j} & 
\sum_{j} d_jk_j \text{\rm\ odd} & \GL_{k_{j}}( q^{d_j} ) 
& \omega_j & \omega_j \rule[- 7pt]{0pt}{22pt} \\ \hline\hline
\end{array}
$$
\bigskip
\noindent \begin{tabular}{ll}
1 & Exactly one of the multiplicities in Types I.1 and I.2 is non-zero.\\
\end{tabular}
\end{table}

\begin{table}[f]
\caption{\label{ConjugateToNegativeOrthogonalC} The critical elements in
$\CSO^-_{2m}(q)$, $m$ odd, $q \equiv 3\,(\text{mod\ } 4)$
(explanations in the paragraph preceding
Proposition~\ref{CorollaryElementsInConformalGroups}
and in Remark~\ref{RemarkForEvenDimensionalSpinGroups})}
$$
\begin{array}{cccccc||cc}\\ \hline\hline
\text{Type} & 
k_\mu(s) & 
d_\mu &
\mbox{\rm dim}( \mathbf{V}_{\tilde{\mu}}( s ) ) & 
\text{Condt's} & 
C^\circ_{[\mathbf{G}_{\tilde{\mu}},\mathbf{G}_{\tilde{\mu}}]}( s_{\tilde{\mu}} )^F
& \lambda & \lambda^a \rule[- 7pt]{0pt}{ 20pt} \\ \hline\hline
\text{I}.1 & k & 1 & 2k & k \geq 2 \text{\ even}^1 &
\SO_{k}^+( q ) \times \SO_{k}^+( q ) & (\Lambda_1,\Lambda_2) & (\Lambda_2,\Lambda_1')
\rule[ 0pt]{0pt}{ 19pt} \\ \hline
\text{I}.2 & k & 1 & 2k & k \geq 2 \text{\ even}^1 & 
\SO_{k}^-( q ) \times \SO_{k}^-( q ) & (\Lambda_1,\Lambda_2) & (\Lambda_2,\Lambda_1)
\rule[ 0pt]{0pt}{ 19pt} \\ \hline
\text{III} & \ell & 2 & 2\ell & 
\ell \text{\ odd} &
\GU_{\ell}( q ) & \kappa & \kappa \rule[ 0pt]{0pt}{ 19pt} \\ \hline
\text{IV}.j, j \in I_{\text{\rm (IV)}} & \ell_{j} & 2e_j & 4 e_j \ell_{j} & 
& \GU_{\ell_{j}}( q^{e_j} ) \times \GU_{\ell_{j}}( q^{e_j} ) 
& (\kappa_{1,j},\kappa_{2,j}) & (\kappa_{2,j},\kappa_{1,j})
\rule[- 6 pt]{0pt}{ 18pt} \\ \hline
\text{V}.j, j \in I_{\text{\rm (V)}} & k_{j} & d_j & 2 d_j k_{j} & 
\sum_{j} d_jk_j \text{\rm\ even} & \GL_{k_{j}}( q^{d_j} ) 
& \omega_j & \omega_j \rule[- 7pt]{0pt}{22pt} \\ \hline\hline
\end{array}
$$
\bigskip
\noindent \begin{tabular}{ll}
1 & Exactly one of the multiplicities in Types I.1, I.2 is non-zero. \\
\end{tabular}
\end{table}
\end{landscape}

\subsection{Centralizers and Levi subgroups}
\label{ClassicalGroupsCentralizers}

In this section we let $(\mathbf{G},F)$, be one of the pairs of a reductive 
algebraic group~$\mathbf{G}$ and a Frobenius morphism~$F$ introduced 
in~\ref{GeneralUnitaryGroups} and~\ref{ConformalGroupsSection}. We also include
the case that $\mathbf{G} = \GL_n( \mathbb{F} )$ and~$F$ is the standard
Frobenius map raising every entry of an element of~$\mathbf{G}$ to its $q$th power.
We let~$\mathbf{V}$ and~$V$ denote the natural vector spaces for~$\mathbf{G}$ 
respectively $G = \mathbf{G}^F$. We study the containment of centralizers
of semisimple elements of~$G$ in Levi subgroups of~$\mathbf{G}$ via the action 
of~$G$ on~$V$. Notice that we do not, in general, assume~$q$ to be odd here.

\begin{lem}
\label{CentralizersActIrreducibly}
Let $s \in G$ be semisimple, and let~$\mu$ be an irreducible factor of the
minimal polynomial of~$s$. Then~$C^\circ_{\mathbf{G}}(s)^F$ acts irreducibly 
on~$V_\mu(s)$, unless~$\mathbf{G}$ is the conformal special orthogonal group 
and~$V_\mu(s)$ is $2$-dimensional and~$\mu$ is of Type~{\rm (I)}, or~$V_\mu(s)$
is $4$-dimensional and~$\mu$ is of Type~{\rm (II)} (in which case~$q$ is odd). 
\end{lem}
\begin{prf}
If $G = \GL_n(q)$, the result is trivial, as then the restriction of 
$C_{\mathbf{G}}(s)^F$ to $\GL(V_\mu(s))$ contains a Singer cycle. Thus assume 
that~$G$ is one of the other groups. The same argument as above works if~$\mu$ 
has Type~(V). We may thus assume that $\mu = \mu^\dagger$ in the unitary case 
or that $\mu = \mu^{*\alpha}$ in the other cases (where~$\alpha$ denotes the 
multiplier of~$s$). We may further assume that $V = V_\mu(s)$. The claim then 
follows from Lemmas~\ref{CentralizersInUnitaryGroups}
and~\ref{SpecialMinimalPolynomial}, as the given centralizers, respectively 
their subgroups $C^\circ_{[\mathbf{G},\mathbf{G}]}(s)^F$, do indeed act 
irreducibly on the spaces~$V_\mu(s)$.
\end{prf}
%%%%%%%%%%%%%%%%%%%%%%%%%%%%%%%%%%%%%%%%%%%%%%%%%%%%%%%%%%%%%%%%%%%%%%%%%%%%%%%%
%%%%%%%%%%%%%%%%%%%%%%%%%%%%%%%%%%%%%%%%%%%%%%%%%%%%%%%%%%%%%%%%%%%%%%%%%%%%%%%%
%%
%% Für mehr Details zum Beweis der Irreduzibilitaet siehe Notizen vom 17.10.2016, 
%% Blaetter 1 - 3
%%
%%%%%%%%%%%%%%%%%%%%%%%%%%%%%%%%%%%%%%%%%%%%%%%%%%%%%%%%%%%%%%%%%%%%%%%%%%%%%%%%
%%%%%%%%%%%%%%%%%%%%%%%%%%%%%%%%%%%%%%%%%%%%%%%%%%%%%%%%%%%%%%%%%%%%%%%%%%%%%%%%

\smallskip

In case~$G$ is a general linear or a unitary group, the statement of 
Lemma~\ref{CentralizersActIrreducibly} follows from the results presented, 
without proof, in \cite[Proposition~(1A)]{fosri1}. In the other cases, if~$q$ is 
odd, one may use \cite[(1.13)]{fosri2}. As we have shown, analogous results also 
hold for~$q$ even, but there does not seem to be a convenient reference in the 
literature.
%%%%%%%%%%%%%%%%%%%%%%%%%%%%%%%%%%%%%%%%%%%%%%%%%%%%%%%%%%%%%%%%%%%%%%%%%%%%%%%%
%%%%%%%%%%%%%%%%%%%%%%%%%%%%%%%%%%%%%%%%%%%%%%%%%%%%%%%%%%%%%%%%%%%%%%%%%%%%%%%%
%%
%% Für mehr Details siehe Notizen vom 11.08.2016, Blatt 1
%%
%%%%%%%%%%%%%%%%%%%%%%%%%%%%%%%%%%%%%%%%%%%%%%%%%%%%%%%%%%%%%%%%%%%%%%%%%%%%%%%%
%%%%%%%%%%%%%%%%%%%%%%%%%%%%%%%%%%%%%%%%%%%%%%%%%%%%%%%%%%%%%%%%%%%%%%%%%%%%%%%%

\begin{lem}
\label{StabilizerOfSplitLevi}
Let $s \in G$ be semisimple and let $H \leq G$ fixing a totally isotropic 
subspace $U \leq V$. Suppose that there is $\mu \in \mathcal{F}_s$ such 
that~$V_\mu(s)$ is $H$-invariant and that~$H$ acts irreducibly on $V_\mu(s)$.
If $U_\mu(s) \neq \{0\}$, then $V_\mu(s) \leq U$. In particular,~$V_\mu(s)$ is 
totally isotropic.
\end{lem}
\begin{prf} As $U_\mu(s) \leq V_\mu(s)$, we obtain
$V_\mu(s) = \langle U_\mu(s)H \rangle \leq U$.
\end{prf}

\begin{lem}
\label{CentralizerContainedInSplitLeviUnitary}
Let $(\mathbf{G},F)$ be unitary as in~{\rm \ref{GeneralUnitaryGroups}}
and let~$s \in G$ be semisimple. Then the following statements hold.

{\rm (a)} Suppose that $C_{\mathbf{G}}(s)^F$ is contained in a proper
split $F$-stable Levi subgroup~$\mathbf{L}$ of~$\mathbf{G}$. Then there
exists $\mu \in \mathcal{F}_{s}$ such that~$\mu^\dagger$ does
not lie in the $\tilde{A}_{\mathbf{L}}( s )^F$-orbit of~$\mu$. In particular,
$\mu \neq \mu^\dagger$.

{\rm (b)} Suppose that there is $\mu \in \mathcal{F}_{s}$ such
that~$\mu^\dagger$ does not lie in the $\tilde{A}_{\mathbf{G}}(s)^F$-orbit
of~$\mu$. Then there is a proper split $F$-stable Levi
subgroup~$\mathbf{L}$ of $\mathbf{G}$ with
$\tilde{C}_{\mathbf{G}}(s)^F\,C_{\mathbf{G}}(s) \leq \mathbf{L}$.

{\rm (c)} The centralizer $C_{\mathbf{G}}( s )$ is contained in a proper split
$F$-stable Levi subgroup of~$\mathbf{G}$ if and only if there is $\mu \in
\mathcal{F}_s$ with $\mu \neq \mu^\dagger$.
\end{lem}
\begin{prf}
(a) Let~$U$ be a non-zero totally isotropic subspace of~$V$ stabilized by~$L$.
Let $\mu \in \mathcal{F}_s$ with $U_\mu(s) \neq 0$.
It follows from Lemmas~\ref{CentralizersActIrreducibly}
and~\ref{StabilizerOfSplitLevi}, the latter applied with $H = 
C_{\mathbf{G}}(s)^F \leq L$, that $V_\mu( s ) \leq U$. In particular, $\mu \neq 
\mu^\dagger$.  As $V_\mu( s )x \leq U$ for all $x \in 
\tilde{C}_{\mathbf{L}}( s )^F \leq L$, the claim follows.

(b) Let $g \in \tilde{C}_{\mathbf{G}}( s )^F$ such that its image in 
$\tilde{A}_{\mathbf{G}}( s )^F$ generates this group (see the first paragraph 
of the proof of Lemma~\ref{CriticalCaseInGUn}).
Suppose that the $\tilde{A}_{\mathbf{G}}(s)^F$-orbit $\mathcal{O}$ of~$\mu$ 
does not contain~$\mu^\dagger$. Then $\mathbf{V}_1 := \sum_{\nu \in \mathcal{O}}
\mathbf{V}_{\nu}( s )$ and $\mathbf{V}_2 := \sum_{\nu \in \mathcal{O}} 
\mathbf{V}_{\nu^\dagger}( s )$ are a pair of non-zero, complementary, totally
isotropic subspaces of~$\mathbf{V}$ fixed by~$g$. Let $\mathbf{L}$ denote the
stabilizer in~$\mathbf{G}$ of~$\mathbf{V}_1$ and~$\mathbf{V}_2$. Then
$\mathbf{L}$ is a proper split $F$-stable Levi subgroup of~$\mathbf{G}$
containing~$g$ and $C_{\mathbf{G}}( s )$. It follows that
$\tilde{C}_{\mathbf{G}}(s)^F = \langle C_{\mathbf{G}}( s )^F, g \rangle 
\leq \mathbf{L}$.

(c) This is a direct consequence of~(a) and~(b).
\end{prf}

We single out $F$-stable semisimple elements of $\CSO_{2m}( \mathbb{F} )$ of a 
special kind.

\begin{dfn}
\label{ExceptionalElements}
{\rm
Let~$\mathbf{G} = \CSO_{2m}( \mathbb{F} )$ and let $s \in G$ be semisimple
with multiplier~$\alpha$. We call~$s$ \textit{exceptional}, if the following 
conditions are satisfied:~$m \geq 2$, every $\nu \in \mathcal{F}_s$ satisfies 
$\nu = \nu^{*\alpha}$, there is $\nu \in \mathcal{F}_s$ of Type~{\rm (I)} with 
$V_\nu(s)$ of dimension~$2$ and $(\mathbf{G}_\nu(s),F)$ of plus-type, and 
if~$q$ is odd, $X^2 - \alpha$ does not divide the minimal polynomial of~$s$.
(Here,~$\mathbf{G}_\nu(s)$ is as defined 
in~{\rm \ref{CriticalSemisimpleElements}}.)
}
\end{dfn}

The latter four conditions can also be phrased as follows: the 
multiplier~$\alpha$ has a square root~$\zeta \in \mathbb{F}_q^*$, the 
$\zeta$-eigenspace $V_{X-\zeta}(s)$ of~$s$ on~$V$ is $2$-dimensional and 
contains a non-zero isotropic vector, and if~$q$ is odd, $-\zeta$ is not an 
eigenvalue of~$s$.

\begin{lem}
\label{CentralizerContainedInSplitLeviConformal}
Let $(\mathbf{G},F)$ be one of the groups considered 
in~{\rm \ref{ConformalGroupsSection}}. Assume that $m \geq 2$ if $\mathbf{G} =
\CSO_{2m}( \mathbb{F} )$. Let $s \in G$ be semisimple with 
multiplier~$\alpha$. Then the following statements hold.

{\rm (a) (i)} Suppose that $C^\circ_{\mathbf{G}}( s )^F$ is contained in a proper 
split $F$-stable Levi subgroup of~$\mathbf{G}$. Then one of the following cases
occurs. 

{\rm (i.1)} There is $\mu \in \mathcal{F}_s$ with $\mu \neq \mu^{*\alpha}$.

{\rm (i.2)} We have $\mathbf{G} = \CSO_{2m}( \mathbb{F} )$, there is 
$\nu \in \mathcal{F}_s$ of Type~{\rm (I)} or Type~{\rm (II)} with 
multiplicity~$2$ in the characteristic polynomial of~$s$ and 
$(\mathbf{G}_\nu(s),F)$ is of plus-type. Moreover, 
$C^\circ_{\mathbf{G}}(s)$ is contained in  proper split
$F$-stable Levi subgroup of~$\mathbf{G}$ stabilizing a pair of complementary,
totally isotropic subspaces $\mathbf{U}', \mathbf{U}'' \leq \mathbf{V}$ such 
that $\mathbf{U}' \oplus \mathbf{U}'' = \mathbf{V}_\nu(s)$. Finally,
$C_{\mathbf{G}}(s)^F$ is contained in a proper split $F$-stable Levi 
subgroup of~$\mathbf{G}$ if and only if~$s$ is exceptional.

{\rm (a) (ii)} Suppose that~$q$ is odd, and that there is $h \in G$ with 
$hsh^{-1} = -s$. If~$s$ is as in~{\rm (i.1)} and
$\langle C^\circ_{\mathbf{G}}( s )^F, h \rangle$ is contained in a proper 
split $F$-stable Levi subgroup of~$\mathbf{G}$, then there is 
$\mu \in \mathcal{F}_s$ with $\mu \neq \mu^{*\alpha} \neq \mu'$.
If~$s$ is as in~{\rm (i.2)} and~$\nu$ is of Type~{\rm(I)}, 
then there is no proper split $F$-stable Levi subgroup 
of~$\mathbf{G}$ containing $\langle C^\circ_{\mathbf{G}}( s )^F, h \rangle$.

{\rm (b)} Suppose that there is $\mu \in \mathcal{F}_s$ with 
$\mu \neq \mu^{*\alpha}$. Then $C_{\mathbf{G}}( s )$ is contained in a proper 
split $F$-stable Levi subgroup of~$\mathbf{G}$. The analogous conclusion holds 
for $\tilde{C}_{\mathbf{G}}( s )$ if $\mu \neq \mu^{*\alpha} \neq \mu'$.

{\rm (c)} The centralizer $C_{\mathbf{G}}( s )$ is contained in a proper split
$F$-stable Levi subgroup of~$\mathbf{G}$ if and only if there is $\mu \in
\mathcal{F}_s$ with $\mu \neq \mu^{*\alpha}$, or $\mathbf{G} = 
\CSO_{2m}( \mathbb{F} )$ and $s$ is exceptional.
\end{lem}
\begin{prf}
(a) For the proof of~(i) let $h = 1$, and for the proof of~(ii) let~$h$ be as
in the assertion. Suppose that~$U$ is a non-zero totally isotropic subspace 
of~$V$ fixed by $\langle C^\circ_{\mathbf{G}}( s )^F, h \rangle$. Suppose first 
that there is $\mu \in \mathcal{F}_s$ such that $U_\mu(s) \neq \{0\}$ and such 
that $C^\circ_{\mathbf{G}}( s )^F$ acts irreducibly on~$V_\mu(s)$. 
Then~$V_\mu(s)$ is totally isotropic by Lemma~\ref{StabilizerOfSplitLevi}.
In particular, $\mu \neq \mu^{*\alpha}$, yielding~(i.1). To prove~(ii) in this 
case, assume 
that $\mu' = \mu^{*\alpha}$. Then $V_{\mu^{*\alpha}}(s) = V_{\mu'}(s) = 
V_{\mu}(s)h \leq Uh \leq U$. It follows that~$U$ contains the
non-degenerate subspace $V_{\mu}(s) \oplus V_{\mu^{*\alpha}}(s)$, a
contradiction.

Now assume that $C^\circ_{\mathbf{G}}( s )^F$ acts reducibly on~$V_\mu(s)$ for 
all $\mu \in \mathcal{F}_s$ with $U_\mu(s) \neq \{0\}$. Let $\nu \in 
\mathcal{F}_s$ with $U_\nu(s) \neq \{0\}$. 
Lemma~\ref{CentralizersActIrreducibly} implies that 
$\mathbf{G} = \CSO_{2m}( \mathbb{F} )$ and that~$\nu$ is of Type~(I) or~(II) 
and occurs with multiplicity~$2$ in the characteristic polynomial of~$s$. 
It follows that $U = U_\nu(s) \lneq V_\nu(s)$. 

Let us prove the other assertions of~(i.2). Put $\mathbf{G}_\nu := \mathbf{G}_\nu(s)$,
and write~$s_\nu$ for the restriction of~$s$ to~$\mathbf{V}_\nu(s)$. Our 
assumption implies that $C^\circ_{\mathbf{G}_\nu}(s_\nu)^F$ fixes~$U$. We claim
that $(\mathbf{G}_\nu,F)$ is of plus-type. Indeed, if~$\nu$ has Type~(I) and
$(\mathbf{G}_\nu,F)$ is of minus type, then 
$V_\nu(s)$ does not have any non-zero totally isotropic subspace. Now suppose
that~$\nu$ has Type~(II). Then $\alpha$ is not a square in~$\mathbb{F}_q^*$, and
in particular~$q$ is odd. Assume that~$F$ acts as~$F''$ on $\mathbf{G}_\nu \cong 
\CSO_4( \mathbb{F} )$. The $F$-stable Borel subgroups are the only proper split 
$F$-stable Levi subgroups of~$\mathbf{G}_\nu$. The $p'$-part of the 
order of a Borel subgroup of $\mathbf{G}_\nu^{F}$ equals $(q - 1)(q^2 - 1)$. 
On the other hand, the order of $C^\circ_{\mathbf{G}_\nu}( s_\nu )^F$ is 
divisible by $q^2 + 1$ (see Lemma~\ref{SpecialMinimalPolynomial}). Thus 
$C^\circ_{\mathbf{G}_\nu}(s_\nu)^F$ cannot fix~$U$, a contradiction. Our claim
is proved. Any non-trivial minimal $s$-invariant subspace of $V_\nu(s)$ is
totally isotropic by Lemma~\ref{SpecialMinimalPolynomial}. Thus~$s_\nu$ lies in 
a split Levi subgroup of~$\mathbf{G}_\nu$. It is now easy to see 
that $C^\circ_{\mathbf{G}_\nu}( s_\nu )$ fixes 
a pair of complementary, totally isotropic subspace $\mathbf{U}', \mathbf{U}'' 
\leq \mathbf{V}_\nu(s)$ with $\mathbf{U}' \oplus \mathbf{U}'' = \mathbf{V}_\nu(s)$
whose stabilizer in~$\mathbf{G}$ is $F$-stable. Then~$C^\circ_{\mathbf{G}}(s)$ 
also fixes~$\mathbf{U}'$ and $\mathbf{U}''$. It remains to prove the last 
assertion. 
%%%%%%%%%%%%%%%%%%%%%%%%%%%%%%%%%%%%%%%%%%%%%%%%%%%%%%%%%%%%%%%%%%%%%%%%%%%%%%%%
%%%%%%%%%%%%%%%%%%%%%%%%%%%%%%%%%%%%%%%%%%%%%%%%%%%%%%%%%%%%%%%%%%%%%%%%%%%%%%%%
%%
%% Für mehr Details zu diesem Beweisteil siehe Notizen vom 14.11.2016, Seiten 1 - 3.
%%
%%%%%%%%%%%%%%%%%%%%%%%%%%%%%%%%%%%%%%%%%%%%%%%%%%%%%%%%%%%%%%%%%%%%%%%%%%%%%%%%
%%%%%%%%%%%%%%%%%%%%%%%%%%%%%%%%%%%%%%%%%%%%%%%%%%%%%%%%%%%%%%%%%%%%%%%%%%%%%%%%
If~$s$ is exceptional,
then $C_\mathbf{G}(s)$ is connected by Lemma~\ref{NonConnectedCentralizersInCSO} 
and we are done. If~$s$ is not exceptional,~$q$ is odd and $X^2 - \alpha$ divides 
the minimal polynomial of~$s$. 
Let $\tilde{\nu} := \nu\nu'$ if~$\nu$ is of Type~(I), and let $\tilde{\nu} := 
\nu$, otherwise. We claim that there exists $g_{\tilde{\nu}} \in 
C_{\mathbf{G}_{\tilde{\nu}}}(s_{\tilde{\nu}})^F$ mapping some non-zero vector 
$v \in U$ to a vector $v'$ such that $(v,v')$ is a hyperbolic pair 
in~$V_{\nu}(s)$. Indeed, if~$\nu$ has Type~(I), there is an $F$-stable 
element in $\hat{\mathbf{G}}_\nu^F \setminus \mathbf{G}_\nu^F$ achieving 
this; take any element in $\hat{\mathbf{G}}_{\nu'}^F \setminus 
\mathbf{G}_{\nu'}^F$ and multiply the two elements to 
obtain~$g_{\tilde{\nu}}$. If~$\nu$ has Type~(II), identify
$(\mathbf{G}_\nu,F)$ with $(\CSO_{4}( \mathbb{F} ), F' )$ and~$s_\nu$ with
a suitable $\CSO_{4}( \mathbb{F} )$-conjugate of 
$t := h( \zeta, -\zeta; \alpha )$, where $\zeta$ is a square root of~$\alpha$. 
Then some conjugate of the element~$\dot{b}$ introduced in
the proof of Lemma~\ref{ElementsInConformalGroups} satisfies the requirement,
and our claim follows. This implies that~$g_{\tilde{\nu}}$ does not 
fix~$U$. Clearly, $g_{\tilde{\nu}}$ can be extended to an element 
$g \in C_{\mathbf{G}}(s)^F$, proving the last assertion of~(i.2).
%%%%%%%%%%%%%%%%%%%%%%%%%%%%%%%%%%%%%%%%%%%%%%%%%%%%%%%%%%%%%%%%%%%%%%%%%%%%%%%%
%%%%%%%%%%%%%%%%%%%%%%%%%%%%%%%%%%%%%%%%%%%%%%%%%%%%%%%%%%%%%%%%%%%%%%%%%%%%%%%%
%%
%% Für mehr Details zur vorletzten Aussage siehe Notizen vom 08.10.2016, 
%% Seiten 1, 2. Die vorletzte Aussage ist: 
%% It is now easy to see that $C^\circ_{\mathbf{G}_\mu}( s_\mu )$ fixes a 
%% non-zero totally isotropic subspace $U' \leq V_\mu(s)$.
%%
%%%%%%%%%%%%%%%%%%%%%%%%%%%%%%%%%%%%%%%%%%%%%%%%%%%%%%%%%%%%%%%%%%%%%%%%%%%%%%%%
%%%%%%%%%%%%%%%%%%%%%%%%%%%%%%%%%%%%%%%%%%%%%%%%%%%%%%%%%%%%%%%%%%%%%%%%%%%%%%%%

We turn to the proof of~(ii) in this case. 
Assume that the assertion is false. As $\nu \neq \nu' = \nu^{*\alpha}$, we obtain 
a contradiction as in the first paragraph of the proof of~(a). 

(b) Let $\mu \in \mathcal{F}_s$ with $\mu \neq \mu^{*\alpha}$. For the proof of
the first statement, put $\tilde{\mu} := \mu$. For the proof of the second 
statement, set~$\tilde{\mu}$  to be equal to~$\mu$ if $\mu = \mu'$, and to 
$\mu\mu'$ if $\mu \neq \mu'$. Then $\mathbf{V}_{\tilde{\mu}}( s )$ and 
$\mathbf{V}_{\tilde{\mu}^{*\alpha}}( s )$ 
form a pair of complementary, non-zero, totally isotropic subspaces 
of~$\mathbf{V}$ which are fixed by $C_{\mathbf{G}}( s )$, and also by
$\tilde{C}_{\mathbf{G}}( s )$ if $\mu' \neq \mu^{*\alpha}$.

(c) By Lemma~\ref{NonConnectedCentralizersInCSO}, $C_{\mathbf{G}}(s)$ is 
connected if $\mathbf{G} = \CSp_{2m}( \mathbb{F} )$ or $\mathbf{G} = 
\CSO_{2m}( \mathbb{F} )$ and~$s$ is exceptional,
so that the \textit{if} part of the assertion follows from~(b) and~(a)(i.2).
For the \textit{only if} part use (a)(i).
\end{prf}

\section{The quasisimple groups of Lie type}
\label{QuasisimpleGroupsOfLieType}

\subsection{Introduction}
\label{QuasisimpleGroupsOfLieTypeIntroduction}
We are now going to describe the Harish-Chandra imprimitive ordinary absolutely
irreducible characters of all finite quasisimple groups of Lie type. Except for 
finitely many such groups, an absolutely irreducible character is 
Harish-Chandra imprimitive if and only if it is imprimitive 
(see~\cite[Theorem~$6.1$]{HiHuMa}). As we do not exclude the groups providing
exceptions right away, we will stick to the notion of Harish-Chandra 
imprimitivity. The above task has already been 
achieved for groups with an exceptional Schur multiplier and for groups 
with two distinct defining characteristics in \cite[Chapter~5]{HiHuMa}, for the 
Tits simple group in \cite[Chapter~3]{HiHuMa}, and for the exceptional series
$G_2(q)$, ${^3\!D}_4(q)$, ${^2\!G}_2( 3^{2m+1} )$, ${^2\!F}_4( 2^{2m+1} )$ and
${^2\!B}_2( 2^{2m+1} )$ in \cite[Section~10.2]{HiHuMa}. The groups $F_4(q)$ and
$E_8(q)$ arise from algebraic groups with connected centers, and are treated 
in~\cite[Section~10.1]{HiHuMa}.

We are thus left with the groups~$G$ listed in 
Subsection~\ref{PrimitivityUnderRestriction}, as well as with the symplectic and 
orthogonal groups $\Sp_{2m}( q )$ and $\Omega_{2m}^\pm( q )$, where~$q$ is even. 
The latter groups arise from algebraic groups with connected center and will 
not be considered here, as their Harish-Chandra imprimitive characters 
have been classified in terms of centralizers of semisimple elements of the dual 
groups (see \cite[Proposition~9.5]{HiHuMa}).

We take this opportunity to correct some inaccuracies in the statement
of~\cite[Proposition~9.5]{HiHuMa} and its proof. Firstly, to establish the 
equivalence of~(b) and~(c) in \cite[Proposition~9.5]{HiHuMa}, we cite Fong and 
Srinivasan \cite[Proposition~(1A)]{fosri1} and
\cite[(1.13)]{fosri2} for the structure of the centralizers in classical groups. 
In the latter reference, the results are only formulated for odd~$q$. However, 
analogous statements also hold for even~$q$, as we have now sketched in 
Lemmas~\ref{CentralizersInUnitaryGroups} and~\ref{SpecialMinimalPolynomial} 
above. Secondly, and more seriously, 
statements (a) and (b) of~\cite[Proposition~9.5]{HiHuMa} are not equivalent
if~$G$ is an orthogonal group. In~(b) one also has to allow for special 
elements. Also, in order to show that~(b) implies~(a), we implicitly assumed 
a result which we have formulated as Lemma~\ref{StabilizerOfSplitLevi} here. 
The proof of this implication as well as the correct statement~(b) follow at 
once from our Lemmas~\ref{CentralizerContainedInSplitLeviUnitary}(c) 
and~\ref{CentralizerContainedInSplitLeviConformal}(c) above.
Finally, if~$G$ is a conformal group in~\cite[Proposition~9.5]{HiHuMa}, one
should also take into account the multipliers of the semisimple elements 
of~$G$ and replace~$\mu^*$ by~$\mu^{*\alpha}$.

To remedy these deficiencies, we present the following corrected version of
the relevant parts of~\cite[Proposition~9.5]{HiHuMa}.

\begin{prop}
\label{PropositionNineFive}
Let $p = 2$ and let $\mathbf{G}$ be one of the groups $\Sp_{2m}( \mathbb{F} )$
with $m \geq 1$ or $\SO_{2m}( \mathbb{F} )$ with $m \geq 2$, where 
$\SO_{2m}( \mathbb{F} )$ is as defined in \cite[p.~$160$]{Taylor}. Let~$F$
denote a Frobenius endomorphism of~$\mathbf{G}$ arising from an 
$\mathbb{F}_q$-structure on~$\mathbf{G}$, and let $(\mathbf{G}^*, F)$
be dual to $(\mathbf{G},F)$.  Then $G = \mathbf{G}^F$ is one of 
$\Sp_{2m}(q)$, $m \geq 1$ or $\SO^\pm_{2m}(q)$, $m \geq 2$. In particular,~$G$
is quasisimple except for $G \in \{ \Sp_2(2), \Sp_4(2), \SO_4^+(2) \}$. 
If~$\mathbf{G}$ is orthogonal and $G \neq \SO_4^+(2)$, then $G = 
\Omega^\pm_{2m}(q)$.

Let $s \in G^*$ be semisimple. Then the following statements are equivalent.

{\rm (a)} Every element of $\mathcal{E}( G, [s] )$ is Harish-Chandra primitive.

{\rm (b)} Every $\mu \in \mathcal{F}_s$ satisfies $\mu = \mu^*$, and if
$\mathbf{G} = \SO_{2m}( \mathbb{F} )$, then~$s$ is not exceptional.
\end{prop}
\begin{prf}
The stated facts about~$G$ are standard (see \cite[Section~$11$]{Taylor}).
The equivalence of~(a) and~(b) follows from 
\cite[Theorem~$7.3$, Theorem~$8.4$]{HiHuMa}
in conjunction with~\ref{CentralizerContainedInSplitLeviConformal}(c), using the 
fact that the conformal groups corresponding to~$\mathbf{G}^*$ are direct 
products of their derived subgroups (isomorphic to~$\mathbf{G}^*$) with their 
centers.
\end{prf}

\subsection{The special linear groups}
\label{SpecialLinearGroups}
Let $G = \SL_n( q )$, $n \geq 2$. In this case we take $\mathbf{G} = 
\SL_n(\mathbb{F})$ with the standard Frobenius map~$F$ raising every
matrix entry of an element of~$\mathbf{G}$ to its $q$th power. We also take
$\tilde{\mathbf{G}} = \tilde{\mathbf{G}}^* = \GL_n( \mathbb{F} )$, acting on 
the natural vector space $\mathbf{V} = \mathbb{F}^n$. Finally, 
$\mathbf{G}^* = \PGL_n( \mathbb{F} )$ and~$i^*$ is the canonical epimorphism.

Let $s \in \PGL_n( q ) = G^*$ be semisimple and let $\chi \in 
\mathcal{E}( G, [s] )$. Our aim is to decide whether~$\chi$ is Harish-Chandra 
primitive, in terms of the $A_{\mathbf{G}^*}(s)^F$-orbit $[\lambda]$ of unipotent 
characters of~$C^\circ_{G^*}( s )$ associated to the $\tilde{G}$-orbit~$[\chi]$.
For this purpose choose $\tilde{s} \in \tilde{G}^* = \GL_n(q)$ with 
$i^*( \tilde{s} ) = s$ and let $\tilde{\chi} \in 
\mathcal{E}( \tilde{G}, [\tilde{s}] )$ such that~$\chi$ is a constituent of 
$\Res_G^{\tilde{G}}( \tilde{\chi} )$.

Let $\mathcal{F}_{\tilde{s}} = \{ \mu_0, \ldots , \mu_{e-1} \}$, and for
$0 \leq i \leq e - 1$, let $d_i$ denote the degree of~$\mu_i$ and $k_i$ the 
multiplicity of~$\mu_i$ in the characteristic polynomial of~$\tilde{s}$. 
(For the definition of $\mathcal{F}_{\tilde{s}}$ see 
Subsection~\ref{RemarksOnTransformations}.) Then
$$C_{\tilde{G}^*}( \tilde{s} ) \cong \GL_{k_0}( q^{d_0} ) \times \cdots \times 
\GL_{k_{e-1}}( q^{d_{e-1}} ).$$
Let $\lambda$ be the unipotent character of $C_{\tilde{G}^*}( \tilde{s} )$
such that $\tilde{\chi}$ corresponds to $(\tilde{s},\lambda)$ in Lusztig's
Jordan decomposition of characters. Then $[\chi] \leftrightarrow [\lambda]$
(if $\lambda$ is viewed as a character of $C^\circ_{G^*}( s )$).
Let $\lambda$ be labelled by the $e$-tuple $(\pi_0, \ldots , \pi_{e-1})$ of 
partitions $\pi_i \vdash k_i$ for $0 \leq i \leq e - 1$.

Let $\beta \in \mathbb{F}_q^*$ denote a generator of the setwise stabilizer of
$\mathcal{F}_{\tilde{s}}$ in $\mathbb{F}_q^*$ (see \ref{SubsectionPreliminaries}).
%\addtocounter{thm}{2}
\begin{thm}
\label{SLnTheorem}
Assume the notation introduced in~{\rm \ref{SpecialLinearGroups}}.
Then $\chi$ is Harish-Chandra primitive if and only 
if~$\langle \beta \rangle$ permutes $\mathcal{F}_{\tilde{s}}$ transitively,
$k_0 = k_1 = \cdots = k_{e-1}$ and $\pi_0 = \pi_1 = \cdots = \pi_{e-1}$.

In particular, {\rm Theorem~\ref{MainResult}(b)} is satisfied for 
$G = \SL_n(q)$.
\end{thm}
\begin{prf}
Let us use the notation introduced in~\ref{JordanDecomposition}. In particular,
$\tilde{C}_{\tilde{\mathbf{G}}^*}( \tilde{s} )$ is the inverse image 
under~$i^*$ of $C_{\mathbf{G}^*}( s )$, i.e.\
$$\tilde{C}_{\tilde{\mathbf{G}}^*}( \tilde{s} ) =
\{ \tilde{g} \in \tilde{\mathbf{G}}^* \mid [\tilde{g}, \tilde{s} ] = \gamma  I_n
\text{\ for some\ } \gamma \in \mathbb{F}^* \}.$$ 
Recall that~$i^*$ induces an equivariant isomorphism between 
$A_{\mathbf{G}^*}( s )^F$ and
$\tilde{C}_{\tilde{\mathbf{G}}^*}( \tilde{s} )^F/C_{\tilde{\mathbf{G}}^*}( \tilde{s} )^F$,
the former group acting on $\mathcal{E}( C^\circ_{G^*}( s ), [1] )$, the
latter group on $\mathcal{E}( C_{\tilde{G}^*}( \tilde{s} ), [1] )$.
It follows from~\ref{RemarksOnTransformations} that
$\tilde{C}_{\tilde{\mathbf{G}}^*}( \tilde{s} )$ permutes the
set $\{ \mathbf{V}_{\mu_i}( \tilde{s} ) \mid 0 \leq i \leq e - 1 \}$.
 
Suppose first that either the action of~$\langle \beta \rangle$ on 
$\mathcal{F}_{\tilde{s}}$ is intransitive, or that $k_i \neq k_j$ for some
$0 \leq i \neq j \leq e - 1$.
Then $\tilde{C}_{\tilde{\mathbf{G}}^*}( \tilde{s} )^F$ does not permute the
spaces $\mathbf{V}_{\mu_i}( \tilde{s} )$ transitively. Hence
$\tilde{C}_{\tilde{\mathbf{G}}^*}( \tilde{s} )^F\,{C}_{\tilde{\mathbf{G}}^*}( \tilde{s} )$ 
is contained in a proper split $F$-stable Levi subgroup 
$\tilde{\mathbf{L}}^*$ of~$\tilde{\mathbf{G}}^*$,
which implies that $C_{\mathbf{G}^*}( s )^F\,C^\circ_{\mathbf{G}^*}( s )  
\leq \mathbf{L}^*$. This is exactly Condition~(\ref{InclusionCondition}) 
of Theorem~\ref{MainResult}.
It follows from Corollary~\ref{SufficientImprimitivityCondition} that~$\chi$ 
is Harish-Chandra induced from~$L$.

Now suppose that~$\langle \beta \rangle$ permutes $\mathcal{F}_{\tilde{s}}$ 
transitively, and that $k_0 = k_1 = \cdots = k_{e-1}$. We claim that in this 
case $\tilde{C}_{\tilde{\mathbf{G}}^*}( \tilde{s} )^F$ acts transitively on 
the set $\{ \mathbf{V}_{\mu_i}( \tilde{s} ) \mid 0 \leq i \leq e - 1 \}$. 
To see this, let~$k$ be the common value of the $k_i$'s. Then the
characteristic polynomial of~$\tilde{s}$ equals $(\mu_0 \cdots \mu_{e-1})^k$.
As~$\beta$ stabilizes $\{ \mu_0, \ldots , \mu_{e-1} \}$, the element
$\beta\tilde{s}$ has the same characteristic polynomial as~$\tilde{s}$,
and thus $\tilde{s}$ and~$\beta \tilde{s}$ are conjugate by an element of
$\tilde{\mathbf{G}}^*$. As~$\tilde{s}$ and~$\beta \tilde{s}$ are $F$-stable and the
centralizers of semisimple elements in~$\tilde{\mathbf{G}}^*$ are connected, there
is $\tilde{g} \in \tilde{G}^*$ with $\tilde{g} \tilde{s} 
\tilde{g}^{-1} = \beta \tilde{s}$. By definition,
$\tilde{g} \in \tilde{C}_{\tilde{\mathbf{G}}^*}( \tilde{s} )^F$ and our claim
follows from~\ref{RemarksOnTransformations}. As a consequence, 
$\tilde{C}_{\tilde{\mathbf{G}}^*}( \tilde{s} )^F$ is not contained in any
proper split $F$-stable Levi subgroup of~$\tilde{\mathbf{G}}^*$ (use the fact
that $C_{\tilde{\mathbf{G}}^*}( \tilde{s} )^F$ acts irreducibly on each
$V_{\mu_i}( \tilde{s} )$; see Lemma~\ref{CentralizersActIrreducibly}). In turn, 
the same statement holds for $C_{\mathbf{G}^*}( s )^F$ by~\ref{RegularEmbeddings}.

If two entries of the label $( \pi_0 , \ldots , \pi_{e-1} )$ of~$\lambda$ 
are distinct, the 
stabilizer~$\tilde{C}_{\tilde{\mathbf{G}}^*}( \tilde{s} )_\lambda^F$
of~$\lambda$ in~$\tilde{C}_{\tilde{\mathbf{G}}^*}( \tilde{s} )^F$ is not
transitive on $\{ \mathbf{V}_{\mu_i}( \tilde{s} ) \mid 0 \leq i \leq e - 1 \}$ 
any longer (as a proper subgroup of a transitive cyclic group acting on a set
with more than~$1$ element is not transitive any more), and thus lies in a 
proper split $F$-stable Levi subgroup~$\tilde{\mathbf{L}}^*$ 
of~$\tilde{\mathbf{G}}^*$, such that~$\tilde{\mathbf{L}}^*$ also contains 
$C_{\tilde{\mathbf{G}}^*}( \tilde{s} )$. Again, Condition~(\ref{InclusionCondition})
of Theorem~\ref{MainResult} is satisfied.
It follows from Corollary~\ref{SufficientImprimitivityCondition}
that $\chi$ is Harish-Chandra induced from~$L$.

Finally suppose that $\pi_0 = \pi_1 = \cdots = \pi_{e-1}$. 
Then~$\tilde{C}_{\tilde{\mathbf{G}}^*}( \tilde{s} )^F$ stabilizes~$\lambda$,
and Corollary~\ref{SufficientPrimitivityConditionNew}(b) implies that~$\chi$ is 
Harish-Chandra primitive.
\end{prf}

%\addtocounter{subsection}{2}
\subsection{The special unitary groups}

As in the case of the special linear groups, we take $\tilde{\mathbf{G}} = 
\tilde{\mathbf{G}}^* = \GL_n( \mathbb{F} )$, acting on the natural vector space 
$\mathbf{V} = \mathbb{F}^n$. Also, $\mathbf{G}^* = \PGL_n( \mathbb{F} )$ 
and~$i^*$ is the canonical epimorphism. We let~$F$ denote the Frobenius
endomorphism of~$\mathbf{G}$ introduced in~\ref{SpecialUnitaryGroupsI}.
Then $\tilde{G} = \tilde{\mathbf{G}}^F = (\tilde{\mathbf{G}}^*)^F = \tilde{G}^* 
= \GU_n(q) \leq \GL_n(q^2)$, $G = \mathbf{G}^F = \SU_n( q )$, and $G^* = 
{\mathbf{G}^*}^F = \PGU_n( q )$. We assume that $n \geq 3$.

To present the main theorem of this subsection, we use the notation and concepts
introduced in Lemmas~\ref{CentralizersInUnitaryGroups},~\ref{CriticalCaseInGUn}. 
Suppose that $s \in {G}^*$ and $\tilde{s} \in \tilde{G}^*$ 
are semisimple with $i^*(\tilde{s}) = s$. As in \ref{SpecialUnitaryGroupsI} we 
write $\tilde{C}_{\tilde{\mathbf{G}}^*}( \tilde{s} ) := 
(i^*)^{-1}( C_{\mathbf{G}^*}( s ) ) = 
\{ \tilde{g} \in \tilde{\mathbf{G}}^* \mid \tilde{g}\tilde{s}\tilde{g}^{-1} = 
\gamma s \text{\ for some\ } \gamma \in \mathbb{F}^* \}$
and put $\tilde{A}_{\tilde{\mathbf{G}}^*}( \tilde{s} ) := 
\tilde{C}_{\mathbf{G}^*}( \tilde{s} )/C_{\tilde{\mathbf{G}}^*}( \tilde{s} )$.
Recall from~\ref{JordanDecomposition}, that $i^*$ induces an isomorphism
$\tilde{A}_{\tilde{\mathbf{G}}^*}( \tilde{s} )^F
\rightarrow A_{\mathbf{G}^*}( s )^F$ which commutes with the actions of
$\tilde{A}_{\tilde{\mathbf{G}}^*}( \tilde{s} )^F$ on
$\mathcal{E}( C_{\tilde{\mathbf{G}}^*}( \tilde{s} )^F, [1] )$ and of
$A_{\mathbf{G}^*}( s )^F$ on $\mathcal{E}( C^\circ_{\mathbf{G}^*}( s )^F, [1] )$,
respectively, and induces, by transport of structure, an action of
$A_{\mathbf{G}^*}( s )^F$ on~$\mathcal{F}_{\tilde{s}}$.

Notice that the structure of the centralizer 
$C_{\tilde{\mathbf{G}}^*}( \tilde{s} )^F$ and the action 
of~$\tilde{A}_{\tilde{\mathbf{G}}^*}( \tilde{s} )^F$ on this centralizer can 
be determined with Lemma~\ref{CriticalCaseInGUn}(b).
The automorphism group of a finite general unitary group fixes every 
unipotent character of the group (see \cite[Remarks on p.~$159$]{Lu}).
Therefore, the stabilizer 
$\tilde{A}_{\tilde{\mathbf{G}}^*}( \tilde{s} )^F_{\tilde{\lambda}}$ of a 
unipotent character $\tilde{\lambda}$ of 
$C_{\tilde{\mathbf{G}}^*}( \tilde{s} )^F$ can easily be found by inspection. 
By the remarks above, 
$\tilde{A}_{\tilde{\mathbf{G}}^*}( \tilde{s} )^F_{\tilde{\lambda}} \cong 
A_{\mathbf{G}^*}(s)^F_\lambda$, if~$\tilde{\lambda}$ is obtained from
$\lambda \in \mathcal{E}( C_{\mathbf{G}^*}( s )^F, [ 1 ] )$ by inflation.

\begin{thm}\label{RestrictionToSUn}
Let $s \in G^* = \PGU_n(q)$ be semisimple and choose $\tilde{s} \in \GU_n(q)$ 
with $i^*( \tilde{s} ) = s$. Let $\chi \in \mathcal{E}( G, [s] )$ and $\lambda 
\in \mathcal{E}( C^\circ_{\mathbf{G}^*}( s )^F, [ 1 ] )$ such that $[\chi] \leftrightarrow 
[\lambda]$. 

{\rm (a)} If there are $A_{\mathbf{G}^*}(s)^F$-orbits on 
$\mathcal{F}_{\tilde{s}}$ of type {\rm (ls)}, then~$\chi$ is Harish-Chandra 
imprimitive.

{\rm (b)} If every $A_{\mathbf{G}^*}(s)^F$-orbit on $\mathcal{F}_{\tilde{s}}$
is of type {\rm (u)}, then~$\chi$ is Harish-Chandra primitive.

{\rm (c)} Suppose that there is at least one $A_{\mathbf{G}^*}(s)^F$-orbit of 
type {\rm (lt)} on $\mathcal{F}_{\tilde{s}}$. Then $q$ is odd. Suppose further
that there are no $A_{\mathbf{G}^*}(s)^F$-orbits on $\mathcal{F}_{\tilde{s}}$ 
of type {\rm (ls)}. Define~$e$ such that~$2e$ is the length of an orbit of 
type {\rm (lt)} on~$\mathcal{F}_{\tilde{s}}$, and such that $\nu_2(2e) \leq 
\nu_2(2e')$ for all other orbit lengths $2e'$ of this type.
Then~$\chi$ is Harish-Chandra primitive if and only if 
$$\nu_2 \left( 
\frac{2e|A_{\mathbf{G}^*}(s)_\lambda^F|}{|A_{\mathbf{G}^*}(s)^F|}
\right) > 0$$ 
(here, $\nu_2$ denotes the $2$-adic valuation of~$\mathbb{Q}$).

{\rm (d)} {\rm Theorem~\ref{MainResult}(b)} is satisfied for $G = \SU_n(q)$.
\end{thm}
\begin{prf}
Part~(d) will be proved as we go along.

(a) This is an immediate consequence of 
Lemma~\ref{CentralizerContainedInSplitLeviUnitary}(b) and
Corollary~\ref{SufficientImprimitivityCondition}. Notice that in this case
Condition~(\ref{InclusionCondition}) of Theorem~\ref{MainResult} is satisfied.

(b) Suppose that~$\chi$ is Harish-Chandra imprimitive. 
Let~$\mathbf{L}^*$ be a proper split $F$-stable Levi subgroup of~$\mathbf{G}^*$ 
and let~$\vartheta$ be a character of~$L$ such that
$R_L^G( \vartheta ) = \chi$. We may assume that $s \in L^*$ (see
\cite[Proposition 15.7]{CaEn}). Let~$\nu$ denote the unipotent character of
$C_{\mathbf{L}^*}^\circ( s )^F$ such that the $\tilde{L}$-orbit of~$\vartheta$
corresponds to the $A_{\mathbf{L}^*}(s)^F$-orbit of~$\nu$ under Lusztig's
generalized Jordan decomposition of characters. It follows
from~(\ref{OrbitLength}) and Theorem~\ref{SameComponentGroupConverse} that
$|A_{\mathbf{G}^*}( s )_\lambda^F| \leq |A_{\mathbf{L}^*}( s )_\nu^F|$.

We claim that $C_{\tilde{\mathbf{G}}^*}( \tilde{s} )^F = 
C_{\tilde{\mathbf{L}}^*}( \tilde{s} )^F$ and that $\nu \in [\lambda]$. (The
second assertion of this claim will only be used in the proof of~(c).) Suppose
this claim has been proved. Then $C_{\tilde{\mathbf{G}}^*}( \tilde{s} )^F \leq 
\tilde{\mathbf{L}}^*$. By Lemma~\ref{CentralizerContainedInSplitLeviUnitary}(a), 
there is $\mu \in \mathcal{F}_{\tilde{s}}$ with $\mu \neq \mu^\dagger$. This, 
however, contradicts our hypothesis.

It remains to prove the claim. Let $\tilde{\vartheta} \in
\mathcal{E}( \tilde{L}, [\tilde{s}] )$ denote an irreducible constituent of
$\Ind_{L}^{\tilde{L}}( \vartheta )$ (see \cite[Proposition 11.7(b)]{CeBo2}).
Put $C_{\tilde{G}} := C_{\tilde{\mathbf{G}}^*}( \tilde{ s } )^F$
and $C_{\tilde{L}} := C_{\tilde{\mathbf{L}}^*}( \tilde{ s } )^F$. We
let~$\tilde{\lambda}$ and~$\tilde{\nu}$ denote the unipotent characters of 
$C_{\tilde{G}}$ and $C_{\tilde{L}}$, respectively, which are obtained 
from~$\lambda$, respectively~$\nu$, by inflation.
As $\tilde{\mathbf{G}}^* = \GL_n( \mathbb{F} )$, the centralizers
$C_{\tilde{\mathbf{G}}^*}( \tilde{ s } )$ and
$C_{\tilde{\mathbf{L}}^*}( \tilde{ s } )$ are regular subgroups
of~$\tilde{\mathbf{G}}^*$. Hence there is a linear character $\hat{s}$ of
$C_{\tilde{G}}$, associated to~$\tilde{s}$ as in \cite[(8.19)]{CaEn}, such that
$$\tilde{\vartheta} = \pm R_{C_{\tilde{L}}}^{\tilde{L}}( \hat{s} \cdot \tilde{\nu} ),$$
where~$\hat{s}$ is viewed as a character of $C_{\tilde{L}} \leq 
C_{\tilde{G}}$ through restriction (see \cite[Proposition 15.10(ii)]{CaEn}).
By the transitivity of twisted induction and by \cite[(8.20)]{CaEn}, we find that
\begin{eqnarray*}
R_{\tilde{L}}^{\tilde{G}}( \tilde{\vartheta} ) & = & 
\pm R_{\tilde{L}}^{\tilde{G}}( R_{C_{\tilde{L}}}^{\tilde{L}}( \hat{s} \cdot \tilde{\nu} ) ) \\
& = & \pm R_{C_{\tilde{L}}}^{\tilde{G}}( \hat{s} \cdot \tilde{\nu} ) \\
& = & \pm R_{C_{\tilde{G}}}^{\tilde{G}}( R_{C_{\tilde{L}}}^{C_{\tilde{G}}}( \hat{s} \cdot \tilde{\nu} ) ) \\
& = & \pm R_{C_{\tilde{G}}}^{\tilde{G}}( \hat{s} \cdot R_{C_{\tilde{L}}}^{C_{\tilde{G}}}( \tilde{\nu} ) ).
\end{eqnarray*}
Let~$\tilde{\rho}$ be an irreducible constituent
of~$R_{C_{\tilde{L}}}^{C_{\tilde{G}}}( \tilde{\nu} )$. By the above equation,
$$\pm R_{C_{\tilde{G}}}^{\tilde{G}}( \hat{s} \cdot \tilde{\rho} )$$
is an irreducible constituent
of~$R_{\tilde{L}}^{\tilde{G}}( \tilde{\vartheta} )$, and hence, by Lemma~\ref{7}, 
of $\Ind_G^{\tilde{G}}( \chi )$.
As $[\chi] \leftrightarrow [\lambda]$, we have $\rho \in [\lambda]$, 
where~$\rho$ is the unipotent character of $C_{\mathbf{G}^*}(s)^F$ such 
that~$\tilde{\rho}$ is inflated from~$\rho$ (cf.~\ref{JordanDecomposition}).
Now $C_{\tilde{\mathbf{L}}^*}( \tilde{ s } )$ is a split Levi subgroup of
$C_{\tilde{\mathbf{G}}^*}( \tilde{ s } )$, so that the map
$R_{C_{\tilde{L}}}^{C_{\tilde{G}}}$ is Harish-Chandra induction.
%%%%%%%%%%%%%%%%%%%%%%%%%%%%%%%%%%%%%%%%%%%%%%%%%%%%%%%%%%%%%%%%%%%%%%%%%%%%%%%%
%%%%%%%%%%%%%%%%%%%%%%%%%%%%%%%%%%%%%%%%%%%%%%%%%%%%%%%%%%%%%%%%%%%%%%%%%%%%%%%%
%% 
%% The fact that $C_{\tilde{\mathbf{L}}^*}( \tilde{ s } )$ is a split Levi 
%% subgroup of $C_{\tilde{\mathbf{G}}^*}( \tilde{ s } )$ follows from
%% [Carter, Proposition 2.8.9]. The intersection of a parabolic subgroup P with 
%% a Levi subgroup L is a parabolic subgroup, whose Levi complement is an
%% intersection of a Levi subgroup of P with L. If P and L are F-stable, so is
%% their intersection.
%%
%%%%%%%%%%%%%%%%%%%%%%%%%%%%%%%%%%%%%%%%%%%%%%%%%%%%%%%%%%%%%%%%%%%%%%%%%%%%%%%%
%%%%%%%%%%%%%%%%%%%%%%%%%%%%%%%%%%%%%%%%%%%%%%%%%%%%%%%%%%%%%%%%%%%%%%%%%%%%%%%%
As $R_{\tilde{L}}^{\tilde{G}}( \tilde{\vartheta} )$ is multiplicity free
by~(\ref{EquationFromLemma7}) and Lusztig's result (see \cite[Section~$10$]{Lu}
and \cite[Proposition 15.11]{CaEn}), the same is true for
$R_{C_{\tilde{L}}}^{C_{\tilde{G}}}( \tilde{\nu} )$. Now~$C_{\tilde{G}}$ is a direct
product $C_1 \times \cdots \times C_l$, where the ~$C_i$ correspond to the
$A_{\mathbf{G}^*}(s)^F$-orbits on~$\mathcal{F}_{\tilde{s}}$. Each of the~$C_i$
is of the form~(\ref{Casel}) or~(\ref{Caseu}), and
$A_{\mathbf{G}^*}(s)^F$ acts on~$C_{\tilde{G}}$ by transitively permuting the
direct factors of the~$C_i$. Apply Lemma~\ref{DirectProductHCInduction} to each
of the~$C_i$. The claim follows from this in conjunction
with~\cite[Lemma 8.2]{HiHuMa}.

(c) The fact that~$q$ is odd follows from Lemma~\ref{CriticalCaseInGUn}(b).
Let $\mu \in \mathcal{F}_{\tilde{s}}$ with $\mu \neq \mu^\dagger$ and let~$2e'$
denote the length of the $A_{\mathbf{G}^*}(s)^F$-orbit of~$\mu$. Thus $2e' = 
|A_{\mathbf{G}^*}(s)^F\colon\!S|$, where~$S$ denotes the stabilizer of~$\mu$ 
in $A_{\mathbf{G}^*}(s)^F$. Therefore,
$$\frac{2e'|A_{\mathbf{G}^*}(s)_\lambda^F|}{|A_{\mathbf{G}^*}(s)^F|} =
\frac{|A_{\mathbf{G}^*}(s)_\lambda^F|}{|S|}.$$
As~$A_{\mathbf{G}^*}(s)^F$ is cyclic, 
$\nu_2( |A_{\mathbf{G}^*}(s)_\lambda^F|/|S|) \leq 0$, if and only if~$S$ 
contains a Sylow $2$-subgroup of~$A_{\mathbf{G}^*}(s)_\lambda^F$. This is the
case if and only if the length 
$|A_{\mathbf{G}^*}(s)_\lambda^F|/|S \cap A_{\mathbf{G}^*}(s)_\lambda^F|$
of the $A_{\mathbf{G}^*}(s)_\lambda^F$-orbit of~$\mu$ is odd. In turn, this
is equivalent to the statement that this orbit does not contain $\mu^\dagger$.
%%%%%%%%%%%%%%%%%%%%%%%%%%%%%%%%%%%%%%%%%%%%%%%%%%%%%%%%%%%%%%%%%%%%%%%%%%%%%%%%
%%%%%%%%%%%%%%%%%%%%%%%%%%%%%%%%%%%%%%%%%%%%%%%%%%%%%%%%%%%%%%%%%%%%%%%%%%%%%%%%
%%
%% Für die letzte Äquivalenz siehe die Notizen vom 14.08.2013, Seite 5,
%% sowie vom 17.08.2016, Seite 1.
%%
%%%%%%%%%%%%%%%%%%%%%%%%%%%%%%%%%%%%%%%%%%%%%%%%%%%%%%%%%%%%%%%%%%%%%%%%%%%%%%%%
%%%%%%%%%%%%%%%%%%%%%%%%%%%%%%%%%%%%%%%%%%%%%%%%%%%%%%%%%%%%%%%%%%%%%%%%%%%%%%%%

Suppose first that 
$\nu_2( 2e|A_{\mathbf{G}^*}(s)_\lambda^F|/|A_{\mathbf{G}^*}(s)^F|) \leq 0$
and choose $\mu$ such that its orbit has length $2e$. By the above,
the $A_{\mathbf{G}^*}(s)_\lambda^F$-orbit of~$\mu$ does not 
contain~$\mu^\dagger$.
Let~$\mathbf{\tilde{L}}^*$ denote the stabilizer of the corresponding pair
of totally isotropic subspaces of~$\mathbf{V}$. Then $C^\circ_{\mathbf{G}^*}( s ) 
\leq \mathbf{L}^*$ and
$A_{\mathbf{G}^*}(s)_\lambda^F \leq A_{\mathbf{L}^*}(s)^F$, i.e.\ 
$C_{\mathbf{G}^*}( s )^F_\lambda \leq \mathbf{L}^*$. If follows from 
Corollary~\ref{SufficientImprimitivityCondition},
that~$\chi$ is Harish-Chandra induced from~$L$. Again, 
Condition~(\ref{InclusionCondition}) of Theorem~\ref{MainResult} is satisfied.

Now suppose that
$\nu_2( 2e|A_{\mathbf{G}^*}(s)_\lambda^F|/|A_{\mathbf{G}^*}(s)^F|) > 0$.
If follows from our choice of~$e$ and the above considerations,
that for every $\mu \in \mathcal{F}_{\tilde{s}}$ with 
$\mu \neq \mu^\dagger$ there is an element in $A_{\mathbf{G}^*}(s)_\lambda^F$
mapping $\mu$ to $\mu^\dagger$. We aim to show that under these 
conditions~$\chi$ is Harish-Chandra primitive.

Suppose that this is not the case and adopt the notation of the proof of~(b).
By the claim in that proof, $C_{\tilde{\mathbf{G}}^*}( \tilde{s} )^F = 
C_{\tilde{\mathbf{L}}^*}( \tilde{s} )^F$ and $\nu \in [\lambda]$. This implies 
that $|A_{\mathbf{L}^*}( s )_\nu^F| \leq
|A_{\mathbf{G}^*}( s )_\nu^F| = |A_{\mathbf{G}^*}( s )_\lambda^F|$. As we
also have $|A_{\mathbf{G}^*}( s )_\lambda^F| \leq 
|A_{\mathbf{L}^*}( s )_\nu^F|$, it follows that
$|A_{\mathbf{G}^*}( s )_\lambda^F| = |A_{\mathbf{L}^*}( s )_\nu^F|$ and this
number divides $|A_{\mathbf{L}^*}( s )^F|$. Hence 
$A_{\mathbf{G}^*}( s )_\lambda^F \leq A_{\mathbf{L}^*}( s )^F$, as each of 
these groups is a subgroup of the cyclic group $A_{\mathbf{G}^*}( s )^F$. Now 
$C_{\tilde{\mathbf{G}}^*}( \tilde{s} )^F \leq \tilde{\mathbf{L}}^*$, 
and~$\mu^\dagger$ lies in the $A_{\mathbf{G}^*}( s )_\lambda^F$-orbit of~$\mu$ 
for every $\mu \in \mathcal{F}_{\tilde{s}}$ with $\mu \neq \mu^\dagger$. This 
contradicts Lemma~\ref{CentralizerContainedInSplitLeviUnitary}(a).
\end{prf}

%\addtocounter{subsection}{1}
\subsection{The symplectic groups}
\label{TheSymplecticGroups}

Let $G = \Sp_{2m}( q )$ with $m \geq 2$ and~$q$ odd (see 
Proposition~\ref{PropositionNineFive}). In this case we take 
$\mathbf{G} = \Sp_{2m}(\mathbb{F})$ and $\mathbf{G}^* = 
\SO_{2m+1}( \mathbb{F} )$ with the 
standard Frobenius map~$F$ raising every matrix entry of an element 
of~$\mathbf{G}$, respectively $\mathbf{G}^*$, to its $q$th power. 
By~$\mathbf{V}$ we denote the natural $\mathbb{F}$-vector space for 
$\mathbf{G}^*$, i.e.\ $\mathbf{V} = \mathbb{F}^{2m+1}$.

%\addtocounter{thm}{1}
\begin{prop}
\label{CriterionSOoddDimension}
Let $s \in G^* = \SO_{2m+1}( q )$ be semisimple. Then the following statements 
are equivalent.

{\rm (1)} The centralizer $C_{\mathbf{G^*}}( s )$ is contained in a
proper split $F$-stable Levi subgroup of~$\mathbf{G}^*$.

{\rm (2)} The connected centralizer $C^\circ_{\mathbf{G^*}}( s )$ is
contained in a proper split $F$-stable Levi subgroup of~$\mathbf{G}^*$.

{\rm (3)} The minimal polynomial of~$s$ (viewed as a linear transformation 
on~$\mathbb{F}_{q}^{2m+1}$) has an irreducible factor~$\mu$ with 
$\mu \neq \mu^*$.
\end{prop}
\begin{prf}
Trivially,~(1) implies~(2). 

Suppose that (2) holds. It follows from 
\cite[(1.13) and subsequent remarks]{fosri2} that 
$H := C^\circ_{\mathbf{G}^*}(s)^F$ satisfies the
hypotheses of Lemma~\ref{StabilizerOfSplitLevi}.
Hence there is an irreducible factor $\mu$ of the
minimal polynomial of~$s$ such that $V_{\mu}(s)$ is totally isotropic.
In particular, $\mu \neq \mu^*$. Thus~(2) implies~(3).

Now assume~(3). Then the subspaces $\mathbf{V}_\mu(s)$ and 
$\mathbf{V}_{\mu^*}(s)$ are totally isotropic and fixed by
$C_{\mathbf{G^*}}( s )$. As the stabilizer of this pair of subspaces is a
proper split $F$-stable Levi subgroup of~$\mathbf{G}^*$, this implies~(1).
\end{prf}

It is well known and easy to see that the centralizer $C_{\mathbf{G^*}}( s )$
of a semisimple element $s \in G^*$ is connected if and only if~$-1$ is not an
eigenvalue of~$s$ (in its action on~$\mathbf{V}$).

\begin{cor}
Let $s \in \SO_{2m+1}( q )$ be semisimple. Then every element of 
$\mathcal{E}( \Sp_{2m}( q ), [ s ] )$ is Harish-Chandra imprimitive, if and
only if the minimal polynomial of~$s$ has an irreducible factor~$\mu$ with 
$\mu \neq \mu^*$.

Otherwise, every element of 
$\mathcal{E}( \Sp_{2m}( q ), [ s ] )$ is Harish-Chandra primitive.
\end{cor}
\begin{prf}
This follows from Proposition~\ref{CriterionSOoddDimension} with the help of
Corollaries~\ref{SufficientImprimitivityCondition} 
and~\ref{ThirdLusztigSeriesResult}.
\end{prf}

%\addtocounter{subsection}{2}
\subsection{The odd dimensional spin groups}

Let $G = \Spin_{2m+1}( q )$, $m \geq 2$, and assume that~$q$ is odd. 
(As $\Spin_{2m+1}( q ) \cong \Omega_{2m+1}( q ) \cong \Sp_{2m}( q )$ if~$q$ 
is even, we may assume that~$q$ is odd here; see 
Proposition~\ref{PropositionNineFive}.) In this case we take
$\mathbf{G} = \Spin_{2m+1}(\mathbb{F})$ and $\tilde{\mathbf{G}}$ the 
connected component of the corresponding Clifford group (see 
\cite[\S 9, $\text{n}^\circ 5$]{Bourbaki2} for the definition of the 
Clifford group). Then $\tilde{\mathbf{G}}^* =
\CSp_{2m}( \mathbb{F} )$ (see \cite[8.1]{lusz2}) and
%%%%%%%%%%%%%%%%%%%%%%%%%%%%%%%%%%%%%%%%%%%%%%%%%%%%%%%%%%%%%%%%%%%%%%%%%%%%%%%%
%%%%%%%%%%%%%%%%%%%%%%%%%%%%%%%%%%%%%%%%%%%%%%%%%%%%%%%%%%%%%%%%%%%%%%%%%%%%%%%%
%%
%% Zu Clifford groups siehe Notizen vom 17.08.2016 und 18.08.2016
%%
%%%%%%%%%%%%%%%%%%%%%%%%%%%%%%%%%%%%%%%%%%%%%%%%%%%%%%%%%%%%%%%%%%%%%%%%%%%%%%%%
%%%%%%%%%%%%%%%%%%%%%%%%%%%%%%%%%%%%%%%%%%%%%%%%%%%%%%%%%%%%%%%%%%%%%%%%%%%%%%%%
$\mathbf{G}^* = \PCSp_{2m}( \mathbb{F} )$
with~$i^*$ the natural epimorphism. We write $\mathbf{V} := \mathbb{F}^{2m}$
for the natural vector space of~$\tilde{\mathbf{G}}^*$, and assume 
that~$\mathbf{V}$ is equipped with a non-degenerate symplectic form defined 
over~$\mathbb{F}_q$. Thus $\tilde{\mathbf{G}}^* \leq \GL_{2m}( \mathbb{F} )$
is invariant under the standard Frobenius map~$F$ raising every matrix entry 
of an element of~$\tilde{\mathbf{G}}^*$ to its $q$th power; we let the Frobenius 
map on~$\tilde{\mathbf{G}}$ be dual to the latter.

\begin{thm}
\label{RestrictionToSpinOdd}
Let~$s \in G^*$ be semisimple and let $\tilde{s} \in \tilde{G}^*$ with $s = 
i^*( \tilde{s} )$. Suppose that~$\tilde{s}$ has multiplier~$\alpha$. Then one 
of the following occurs.

{\rm(a)} For all $\mu \in \mathcal{F}_{\tilde{s}}$ we have
$\mu = \mu^{*\alpha}$; in this case every element of $\mathcal{E}( G, [s] )$
is Harish-Chandra primitive.

{\rm(b)} There exists $\mu \in \mathcal{F}_{\tilde{s}}$ with $\mu \neq 
\mu^{*\alpha} \neq \mu'$; in this case every element of $\mathcal{E}( G, [s] )$
is Harish-Chandra imprimitive.

{\rm(c)} Every $\mu \in \mathcal{F}_{\tilde{s}}$ satisfies $\mu = \mu^{*\alpha}$ 
or $\mu' = \mu^{*\alpha}$, and there exists $\mu \in \mathcal{F}_{\tilde{s}}$ 
with $\mu \neq \mu^{*\alpha}$. Here, we distinguish two cases.

{\rm (i)} If~$\tilde{s}$ is not conjugate to~$-\tilde{s}$ in~$\tilde{G}^*$, then 
every element of $\mathcal{E}( G, [s] )$ is Harish-Chandra imprimitive.

{\rm (ii)} Suppose that~$\tilde{s}$ is conjugate to~$-\tilde{s}$ 
in~$\tilde{G}^*$. Then $|A_{\mathbf{G}^*}( s )^F| = 2$. 
Let~$\chi \in \mathcal{E}( G, [s] )$, and 
let $\lambda \in \mathcal{E}( C^\circ_{{G}^*}( s ), [1] )$ be such that 
$[\chi] \leftrightarrow [\lambda]$ under Lusztig's generalized Jordan 
decomposition of characters (see {\rm \ref{JordanDecomposition}}).
Then~$\chi$ is Harish-Chandra imprimitive if and only if
$A_{\mathbf{G}^*}( s )^F_\lambda$ is trivial.
\end{thm}
\begin{prf}
(a) By Lemmas~\ref{CentralizerContainedInSplitLeviConformal}(c)
and~\ref{ComparingCentralizers}, the hypothesis implies that 
$C^\circ_{\mathbf{G}^*}( s ) = i^*( C_{\tilde{\mathbf{G}}^*}( \tilde{s} ) )$ is 
not contained in any proper split $F$-stable Levi subgroup of~$\mathbf{G}^*$. 
It follows from Corollary~\ref{ThirdLusztigSeriesResult} that every element of 
$\mathcal{E}( G, [s] )$ is Harish-Chandra primitive.

(b) In this case, $C_{\mathbf{G}^*}( s ) = 
i^*( \tilde{C}_{\tilde{\mathbf{G}}^*}( \tilde{s} ) )$ is contained in a proper 
split $F$-stable Levi subgroup of~$\mathbf{G}^*$ by 
Lemmas~\ref{CentralizerContainedInSplitLeviConformal}(b) and~\ref{ComparingCentralizers}.
The claim follows from \cite[Theorem~$7.3$]{HiHuMa}.

(c) By assumption, there is $\mu \in \mathcal{F}_{\tilde{s}}$ with $\mu \neq 
\mu^{*\alpha}$. Hence $C^\circ_{\mathbf{G}^*}( s ) = 
i^*( C_{\tilde{\mathbf{G}}^*}( \tilde{s} ) )$
is contained in a proper split $F$-stable Levi subgroup of~$\mathbf{G}^*$
by Lemmas~\ref{CentralizerContainedInSplitLeviConformal}(c) 
and~\ref{ComparingCentralizers}. If~$\tilde{s}$
is not conjugate to~$-\tilde{s}$ in $\tilde{G}^*$, then
$\tilde{C}_{\tilde{\mathbf{G}}^*}( \tilde{s} ) = 
C_{\tilde{\mathbf{G}}^*}( \tilde{s} )$, and thus $C_{\mathbf{G}^*}( s ) = 
C^\circ_{\mathbf{G}^*}( s )$, as $(i^*)^{-1}( C_{\mathbf{G}^*}( s ) ) =
\tilde{C}_{\tilde{\mathbf{G}}^*}( \tilde{s} )$. The result follows as in~(b).

Now suppose that~$\tilde{s}$ is conjugate to~$-\tilde{s}$ by some $\tilde{h} 
\in \tilde{G}^*$. Then~$\tilde{s}$ satisfies the conditions
of Proposition~\ref{CorollaryElementsInConformalGroups}(a). In particular, 
$A_{\mathbf{G}^*}( s )^F$ has order~$2$ (see \ref{JordanDecomposition}).
Moreover, Lemma~\ref{CentralizerContainedInSplitLeviConformal}(a)(ii) implies 
that $C_{\mathbf{G}^*}( s )^F = 
i^*( \tilde{C}_{\tilde{\mathbf{G}}^*}( \tilde{s} )^F )$ is not contained in 
any proper split $F$-stable Levi subgroup of~$\mathbf{G}^*$. The assertion 
follows from Corollaries~\ref{SufficientPrimitivityConditionNew}(b)
and~\ref{SufficientImprimitivityCondition}.
\end{prf}

We finally explain how to determine~$A_{\mathbf{G}^*}( s )^F_\lambda$ in the
situation of Theorem~\ref{RestrictionToSpinOdd}(c)(ii).

\begin{rem}
\label{RemarkForOddDimensionalSpinGroups}
{\rm 
Assume the hypotheses of Theorem~\ref{RestrictionToSpinOdd}(c)(ii).
Table~{\rm \ref{ConjugateToNegativeSymplectic}} lists the various possibilities
for~$\tilde{s}$ by the types of the elements of $\mathcal{F}_{\tilde{s}}$ (the
elements are denoted by~$s$ in that table). The table also gives the labels for 
the unipotent characters $\lambda \in 
\mathcal{E}( C_{[\tilde{\mathbf{G}}^*,\tilde{\mathbf{G}}^*]}( \tilde{s} )^F, [1] )$ 
and the labels of their conjugates~$\lambda^a$ where~$a$ is a generator
of~$\tilde{A}_{\tilde{\mathbf{G}}^*}( \tilde{s} )^F$ (for the definition of the 
latter group see~\ref{SemisimpleInConformal}). The conjugates are determined 
from the action given in Lemma~\ref{ElementsInConformalGroups}(b), and the fact 
that the unipotent character of the nearly simple components involved in 
$C_{[\tilde{\mathbf{G}}^*,\tilde{\mathbf{G}}^*]}( \tilde{s} )^F$ are invariant 
under automorphisms (see \cite[Remarks on p.~$159$]{Lu} and 
\cite[Theorem~$2.5$]{MalleExt}). As labels for the unipotent 
characters we use symbols as defined in \cite[Appendix]{luszbuch} (where the 
condition $k \geq 4$ is imposed; but the symbols can also be defined and used 
for $k = 2$), respectively partitions. 

From this information it is easy to read
off~$A_{\mathbf{G}^*}( s )^F_\lambda$ as follows.
The unipotent characters of $C^\circ_{\mathbf{G}^*}( s )^F$ may be
identified with those of $C_{\tilde{\mathbf{G}}^*}( \tilde{s} )^F$,
and the latter with those of
$C_{[\tilde{\mathbf{G}}^*,\tilde{\mathbf{G}}^*]}( \tilde{s} )^F$. This
identification is compatible with the action 
of~$A_{\mathbf{G}^*}( \tilde{s} )^F$ on the first of these sets, and with 
the action of $\tilde{A}_{\tilde{\mathbf{G}}^*}( \tilde{s} )^F$ on the latter.
}
\end{rem}

\subsection{The even dimensional spin groups}

Assume that $\text{char}(\mathbb{F})$ is odd (see 
Proposition~\ref{PropositionNineFive}). Let $\mathbf{G} = 
\Spin_{2m}(\mathbb{F})$, $m \geq 2$, defined with respect to a non-degenerate 
quadratic form on $\mathbf{V} = \mathbb{F}^{2m}$. Let $\check{\mathbf{G}}$ 
denote  the connected component of the Clifford group with respect to this 
form (see \cite[\S 9, $\text{n}^\circ 5$]{Bourbaki2}), and choose a regular 
embedding $k : \check{\mathbf{G}} \rightarrow \tilde{\mathbf{G}}$. Denote
by~$j$ the embedding $j : \mathbf{G} \rightarrow \check{\mathbf{G}}$. Then
$i := k \circ j : \mathbf{G} \rightarrow \tilde{\mathbf{G}}$ is a regular 
embedding.  We have $\check{\mathbf{G}}^* = \CSO_{2m}( \mathbb{F} )$ and 
$\mathbf{G}^* = \PCSO_{2m}( \mathbb{F} )$, the quotient group of 
$\CSO_{2m}( \mathbb{F} )$ modulo its center. Moreover, 
$j^* : \check{\mathbf{G}}^* \rightarrow \mathbf{G}^*$ is the canonical 
epimorphism, and~$i^*$ 
factors as 
$$
\begin{xy}
\xymatrix{ \tilde{\mathbf{G}}^* \ar[r]^{k^*} \ar[dr]_{i^*} &
\check{\mathbf{G}}^* \ar[d]^{j^*} \\
 & \mathbf{G}^*
}
\end{xy}
$$
%%%%%%%%%%%%%%%%%%%%%%%%%%%%%%%%%%%%%%%%%%%%%%%%%%%%%%%%%%%%%%%%%%%%%%%%%%%%%%%%
%%%%%%%%%%%%%%%%%%%%%%%%%%%%%%%%%%%%%%%%%%%%%%%%%%%%%%%%%%%%%%%%%%%%%%%%%%%%%%%%
%%
%% Zu Clifford groups siehe Notizen vom 17.08.2016 und 18.08.2016
%%
%%%%%%%%%%%%%%%%%%%%%%%%%%%%%%%%%%%%%%%%%%%%%%%%%%%%%%%%%%%%%%%%%%%%%%%%%%%%%%%%
%%%%%%%%%%%%%%%%%%%%%%%%%%%%%%%%%%%%%%%%%%%%%%%%%%%%%%%%%%%%%%%%%%%%%%%%%%%%%%%%
Let~$F'$ and $F''$ denote Frobenius morphisms of~$\tilde{\mathbf{G}}$ such that 
$\mathbf{G}^{F'} = \Spin^+_{2m}(q)$ and $\mathbf{G}^{F''} = \Spin^-_{2m}(q)$ and
such that the induced morphisms on $\check{\mathbf{G}}^*$ are as 
in~\ref{ConformalGroups}. Let~$F$ be one of~$F'$ or~$F''$.
The groups~$\tilde{\mathbf{G}}$ and $\tilde{\mathbf{G}}^*$ are only used in the
proof, but not in the statement of Theorem~\ref{RestrictionToSpinEven} below.

\begin{lem}
\label{CentralizerGeneration}
Let~$s \in G^*$ be semisimple and let $\check{s} \in \check{G}^*$ with $s = 
j^*( \check{s} )$. Then
$$C^\circ_{\mathbf{G}^*}( s ) = 
j^*( C^\circ_{\check{\mathbf{G}}^*}( \check{s} ) )$$ and
$$C^\circ_{\check{\mathbf{G}}^*}( \check{s} ) = 
(j^*)^{-1}( C^\circ_{\mathbf{G}^*}( s ) ).$$
Moreover, putting 
\begin{equation}
\label{TildeC}
\tilde{C}_{\check{\mathbf{G}}^*}( \check{s} ) := \{ \check{g} \in \check{\mathbf{G}}^* 
\mid \check{g} \check{s} \check{g}^{-1} = \pm \check{s} \}
\end{equation}
(cf.\ {\rm \ref{SemisimpleInConformal}}) and
\begin{equation}
\label{TildeA}
\tilde{A}_{\check{\mathbf{G}}^*}( \check{s} ) := 
\tilde{C}_{\check{\mathbf{G}}^*}( \check{s} )/C^\circ_{\check{\mathbf{G}}^*}( \check{s} ),
\end{equation}
then $j^*$ induces an $F$-equivariant isomorphism 
$j^* : \tilde{A}_{\check{\mathbf{G}}^*}( \check{s} ) \rightarrow
A_{\mathbf{G}^*}( s )$, compatible with the actions of these groups on the
unipotent characters of $C^\circ_{\check{\mathbf{G}}^*}( \check{s} )^F$,
respectively $C^\circ_{\mathbf{G}^*}(s)^F$.
\end{lem}
\begin{prf}
Let $\tilde{s} \in \tilde{G}^*$ with $k^*( \tilde{s} ) = \check{s}$. We have
$j^*( k^*( C_{\tilde{\mathbf{G}}^*}( \tilde{ s } ) ) ) =
i^*( C_{\tilde{\mathbf{G}}^*}( \tilde{ s } ) ) = C^\circ_{\mathbf{G}^*}( s )$,
and $k^*( C_{\tilde{\mathbf{G}}^*}( \tilde{ s } ) ) =  
C^\circ_{\check{\mathbf{G}}^*}( \check{s} )$ (see \cite[p.~36]{CeBo2}). Hence
$j^*( C^\circ_{\check{\mathbf{G}}^*}( \check{s} ) ) = C^\circ_{\mathbf{G}^*}( s )$.
As the kernel of~$j^*$ is contained in
$C^\circ_{\check{\mathbf{G}}^*}( \check{s} )$, we obtain the second claim.
The third follows from the second, as $(j^*)^{-1}( C_{\mathbf{G}^*}( s ) ) =
\tilde{C}_{\check{\mathbf{G}}^*}( \check{s} )$.
\end{prf}

\smallskip

\noindent
In the theorem below and its proof we will use the following notation. Let 
$s \in G^*$, choose $\check{s} \in \check{G}^*$ with $j^*(\check{s}) = 2$, and 
choose $\tilde{s} \in \tilde{G}^*$ with $k^*(\tilde{s}) = \check{s}$. By
\cite[p.~$36$]{CeBo2} and Lemma ~\ref{CentralizerGeneration}, we have surjective 
homomorphisms 
$$C_{\tilde{\mathbf{G}}^*}( \tilde{s} )^F 
\stackrel{k^*}{\longrightarrow} C^\circ_{\check{\mathbf{G}}^*}( \check{s} )^F
\stackrel{j^*}{\longrightarrow} C^\circ_{\mathbf{G}^*} ( s )^F.$$
Next, for $\chi \in \mathcal{E}( G, [s] )$, let~$\tilde{\chi}$ denote an
irreducible constituent of $\Ind_{G}^{\tilde{G}}( \chi )$ lying in 
$\mathcal{E}( \tilde{G}, [\tilde{s}] )$, and let~$\check{\chi}$ denote an
irreducible constituent of $\Ind_G^{\check{G}}( \chi )$ such that~$\tilde{\chi}$
occurs in $\Ind_{\check{G}}^{\tilde{G}}( \check{\chi} )$. Then $\check{\chi} \in 
\mathcal{E}( \check{G}, [\check{s}] )$ (see the results summarized 
in~\ref{JordanDecomposition}). Fix a Jordan decomposition of characters
between $\mathcal{E}( \tilde{G}, [\tilde{s}] )$ and 
$\mathcal{E}( C_{\tilde{\mathbf{G}}^*}( \tilde{s} )^F, [ 1 ] )$ and 
let~$\tilde{\lambda}$ denote the unipotent character of 
$C_{\tilde{\mathbf{G}}^*}( \tilde{s} )^F$ corresponding to~$\tilde{\chi}$.
Then the kernel of~$\tilde{\lambda}$ contains the kernels of $k^*|_{\tilde{G}^*}$
and of $i^*|_{\tilde{G}^*}$, and we let~$\check{\lambda}$ und~$\lambda$ denote
the unipotent characters of $C^\circ_{\check{\mathbf{G}}^*}( \check{s} )^F$,
respectively $C^\circ_{\mathbf{G}^*}( s )^F$, obtained from~$\tilde{\lambda}$
by deflation. Then $[\check{\chi}] \leftrightarrow [\check{\lambda}]$ and
$[\chi] \leftrightarrow [\lambda]$ in Lusztig's generalized Jordan decomposition
of characters. Notice that $\check{\lambda}$ is uniquely determined by its 
restriction to 
$C^\circ_{[\check{\mathbf{G}}^*,\check{\mathbf{G}}^*]}( \check{s} )^F$, and we
will usually identify~$\check{\lambda}$ with this restriction.
The configuration considered is depicted in 
Table~\ref{NotationSpinEvenTheorem}.

\begin{table}[h]
\caption{\label{NotationSpinEvenTheorem} The notation in 
Theorem~\ref{RestrictionToSpinEven}
(explanations in the paragraph following Lemma~\ref{CentralizerGeneration})}
$$
\begin{xy}
\xymatrix{\tilde{G} & \mathcal{E}( \tilde{G}, [\tilde{s}] ) \ni \tilde{\chi} \ar@{<->}[rr] & &
\tilde{\lambda} \in \mathcal{E}( C_{\tilde{\mathbf{G}}^*}( \tilde{s} )^F, [1] ) \ar@<-3.9em>[d] & 
\tilde{G}^* \ar[d]^{k^*} \\
\check{G} \ar[u]^k & 
\mathcal{E}( \check{G}, [\check{s}] ) \supseteq [\check{\chi}] \ar@<-2.5em>[u] \ar@{<->}[rr] & &
[\check{\lambda}] \subseteq \mathcal{E}( C_{\check{\mathbf{G}}^*}( \check{s} )^F, [1] ) \ar@<-3.9em>[d] &
\check{G}^* \ar[d]^{j^*} \\
G \ar[u]^j & \mathcal{E}( G, [ s ] ) \supseteq [\chi] \ar@<-2.5em>[u] \ar@{<->}[rr] & &
[\lambda] \subseteq \mathcal{E}( C_{\mathbf{G}^*}( s )^F, [1] ) & G^*
}
\end{xy}
$$
\end{table}

\noindent
The main result of this subsection is formulated in terms of the $F$-stable 
semisimple elements of the conformal group~$\check{\mathbf{G}}^*$, and we will 
make use of the notation introduced in~\ref{CriticalSemisimpleElements} 
and~\ref{GeneralSemisimpleElements} with respect to this group and the action
on its natural vector space. We will, of course, also use the notation of 
Lemma~\ref{CentralizerGeneration}. The following remark describes how to
determime  the action of $A_{\mathbf{G}^*}(s)^F$ on 
$\mathcal{E}( C^\circ_{\mathbf{G}^*}( s )^F, [ 1 ] )$. This information will
be used at various places in the proof without further comment.

\begin{rem}
\label{RemarkForEvenDimensionalSpinGroups}
{\rm
Tables~{\rm \ref{ConjugateToNegativeOrthogonalA}},~{\rm \ref{ConjugateToNegativeOrthogonalB}}
and~{\rm \ref{ConjugateToNegativeOrthogonalC}} list the various possibilities
for~$\check{s}$ by the types of the elements of $\mathcal{F}_{\check{s}}$ (the
elements are denoted by~$s$ there). These tables also give the labels for the
elements of
$\mathcal{E}( C^\circ_{[\check{\mathbf{G}}^*,\check{\mathbf{G}}^*]}( \check{s} )^F, [1] )$
and the labels of their conjugates
under~$\tilde{A}_{\check{\mathbf{G}}^*}( \check{s} )^F$ (for the definition of
the latter group see~(\ref{TildeC}) and~(\ref{TildeA})). Again, the labels are
symbols as in~\cite[Appendix]{luszbuch}. (There, the condition $k \geq 8$ is
imposed; but the symbols can also be defined and used for $k = 2, 4, 6$. For
example, if $k = 2$, there is exactly one relevant symbol, corresponding to the
fact that the unique unipotent character of $\SO_2^\pm(q)$ is the trivial
character.) A pair of symbols labels the two factors of an outer tensor product 
of unipotent characters in a direct product of groups such as, e.g.\ 
$\SO_k^+( q ) \times \SO_k^+ ( q )$.

In the situation of
Table~{\rm \ref{ConjugateToNegativeOrthogonalA}} we have
$\tilde{A}_{\check{\mathbf{G}}^*}( \check{s} )^F = \langle a, b \rangle$ if
$X^2 - \alpha$ divides the minimal polynomial of~$\check{s}$ (i.e.\ if
$k \neq 0$), and
$\tilde{A}_{\check{\mathbf{G}}^*}( \check{s} )^F = \langle a \rangle$, otherwise.
In Tables~{\rm \ref{ConjugateToNegativeOrthogonalB}}
and~{\rm \ref{ConjugateToNegativeOrthogonalC}} we have
$\tilde{A}_{\check{\mathbf{G}}^*}( \check{s} )^F = \langle a \rangle$.
The conjugates are determined from
the action given in Lemma~{\rm \ref{ElementsInConformalGroups}(b)}, and the
fact that the unipotent characters of the nearly simple components involved in
$C^\circ_{[\tilde{\mathbf{G}}^*,\check{\mathbf{G}}^*]}( \check{s} )^F$ are invariant
under automorphisms, except for graph automorphisms of components equal to
$\SO^+_{k}( q )$ for~$k$ divisible by~$4$ (see \cite[Remarks on p.~$159$]{Lu}
and \cite[Theorem~$2.5$]{MalleExt}, which also gives the action of the graph
automorphisms in the latter case). The unipotent characters of $\SO^+_{k}( q )$
for even~$k$ (including $k = 2$, where only the trivial character is
unipotent) are labelled by symbols~$\Lambda$ and copies~$\Lambda'$
thereof, where we use the convention that~$\Lambda$ and~$\Lambda'$ label the
same unipotent character unless~$\Lambda$ is degenerate (in which case~$k$
is divisible by~$4$). Then the graph automorphism of~$\SO^+_{k}( q )$
swaps the characters labelled by $\Lambda$ and~$\Lambda'$ and fixes the
other unipotent characters (see \cite[Theorem~$2.5$]{MalleExt}).
From this information it is easy to read off~$A_{\mathbf{G}^*}( s )^F_\lambda$
as explained in Remark~\ref{RemarkForOddDimensionalSpinGroups}.
}
\end{rem}

\noindent
We now come to our main result for the even dimensional spin groups.

\begin{thm}
\label{RestrictionToSpinEven}
Let~$s \in G^*$ be semisimple and let $\check{s} \in \check{G}^*$ with 
$s = j^*( \check{s} )$. Suppose that~$\check{s}$ has multiplier~$\alpha$. 
Let $\chi \in \mathcal{E}( G, [s] )$, and let 
$\lambda \in \mathcal{E}( C^\circ_{{G}^*}( s ), [1] )$ be such that
$[\chi] \leftrightarrow [\lambda]$ under Lusztig's generalized Jordan
decomposition of characters (see {\rm \ref{JordanDecomposition}}).
Let $\check{\lambda} \in 
\mathcal{E}( C^\circ_{\check{\mathbf{G}}^*}( \check{s} )^F, [ 1 ] )$ denote 
the inflation of~$\lambda$ over the kernel of~$j^*|_{\check{G}^*}$. 
Then one of the following occurs.

{\rm(a)} For all $\mu \in \mathcal{F}_{\check{s}}$ we have 
$\mu = \mu^{*\alpha}$ and~$\check{s}$ is not as in 
Lemma {\rm \ref{CentralizerContainedInSplitLeviConformal}(a)(i.2)}
Then~$\chi$ is Harish-Chandra primitive.

{\rm(b)} There exists $\mu \in \mathcal{F}_{\check{s}}$ with $\mu \neq 
\mu^{*\alpha} \neq \mu'$; in this case~$\chi$ 
is Harish-Chandra imprimitive.

{\rm(c)} Either~$\check{s}$ is as in
{\rm Lemma~\ref{CentralizerContainedInSplitLeviConformal}(a)(i.2)} or every 
$\mu \in \mathcal{F}_{\check{s}}$ satisfies $\mu = \mu^{*\alpha}$
or $\mu' = \mu^{*\alpha}$, and there exists $\mu \in \mathcal{F}_{\check{s}}$
with $\mu \neq \mu^{*\alpha}$. Then there is a proper split $F$-stable Levi 
subgroup~$\mathbf{L}^*$ of~$\mathbf{G}^*$ such that $C^\circ_{\mathbf{G}^*}( s ) 
\leq \mathbf{L}^*$. 
Let~$\mathbf{L}$ denote a split $F$-stable Levi subgroup of~$\mathbf{G}$ dual 
to~$\mathbf{L}^*$.

Then~$\chi$ is Harish-Chandra induced from~$L$ if and only
if $A_{\mathbf{G}^*}( s )^F_\lambda \leq A_{\mathbf{L}^*}( s )^F$.
To investigate this latter condition more closely, we distinguish three cases.

{\rm (i)} Suppose that~$\check{s}$ is not conjugate to~$-\check{s}$ in
$\check{G}^*$. If~$\check{s}$ is as in 
Lemma~{\rm \ref{CentralizerContainedInSplitLeviConformal}(a)(i.2)},
then~$\chi$ is Harish-Chandra induced from~$L$ if either~$\check{s}$ is 
exceptional, or if 
the following conditions are satisfied: There is $\zeta \in \mathbb{F}_q$ with 
$\alpha = \zeta^2$ such that $\nu := X - \zeta$ and $\nu' := X + \zeta$ occur 
with multiplicity~$2$ and $4k' > 0$, respectively, in the characteristic 
polynomial of~$\check{s}$. Moreover, 
$(\check{\mathbf{G}}^*_{\nu}( \check{s} ), F)$
and $(\check{\mathbf{G}}^*_{\nu'}( \check{s} ), F)$ are of plus-type, and
the factor of $\check{\lambda}$ corresponding to $V_{\nu'}( \check{s} )$
is labelled by a degenerate symbol. Otherwise,~$\chi$ is Harish-Chandra primitive.

If~$\check{s}$ is not as in
Lemma~{\rm \ref{CentralizerContainedInSplitLeviConformal}(a)(i.2)},
then  $C_{\mathbf{G}^*}( s ) 
\leq \mathbf{L}^*$ and thus $A_{\mathbf{G}^*}( s )^F_\lambda \leq 
A_{\mathbf{L}^*}( s )^F$. Hence~$\chi$ is Harish-Chandra induced 
from~$L$. 

{\rm (ii)} Suppose that~$\check{s}$ is conjugate to~$-\check{s}$ 
in~$\check{G}^*$ and that $X^2 - \alpha$ does not divide the minimal polynomial 
of~$\check{s}$.  Then $|A_{\mathbf{G}^*}( s )^F| = 2$ and 
$|A_{\mathbf{L}^*}( s )^F| = 1$. In particular,~$\chi$ is Harish-Chandra induced 
from~$L$ if $|A_{\mathbf{G}^*}( s )^F_\lambda| = 1$, and is Harish-Chandra 
primitive, otherwise.

{\rm (iii)} Suppose that~$\check{s}$ is conjugate to~$-\check{s}$ in
$\check{G}^*$ and that $X^2 - \alpha$ divides the minimal 
polynomial of~$\check{s}$. Then $|A_{\mathbf{G}^*}( s )^F| = 4$.

If $\mu = \mu^{*\alpha}$ for all $\mu \in \mathcal{F}_{\check{s}}$, then~$\chi$
is not Harish-Chandra induced from~$L$. If there is 
$\mu \in \mathcal{F}_{\check{s}}$ with $\mu \neq \mu^{*\alpha}$, then 
$|A_{\mathbf{L}^*}( s )^F|= 2$  and the following statements hold.

{\rm (iii.1)}
If~$m$ is odd, then~$\chi$ is Harish-Chandra induced from~$L$
if and only if~$\lambda$ is not fixed by $|A_{\mathbf{G}^*}( s )^F|$.

{\rm (iii.2)}
Suppose that~$m$ is even. Then $G = \Spin_{2m}^+( q )$ and
$A_{\mathbf{G}^*}( s )^F = \langle a, b \rangle$ is a Klein four group. We may
choose notation such that $b \in A_{\mathbf{L}^*}( s )^F$, and~$a$ denotes
the image of an element of $\check{G}^*$ which conjugates 
$\check{s}$ to~$-\check{s}$. Then~$\chi$ is Harish-Chandra induced
from~$L$, if and only if neither~$a$ nor~$ab$ fix~$\lambda$.

{\rm (d)} Suppose that~$\check{s}$ and~$\mathbf{L}^*$ are as in~{\rm (c)} and
that~$\chi$ is not Harish-Chandra induced from~$L$. Then~$\chi$ is 
Harish-Chandra primitive. In particular, {\rm Theorem~\ref{MainResult}(b)} holds 
for~$G$.
\end{thm}
\begin{prf}
(a) Suppose that~$\chi$ is Harish-Chandra imprimitive. Let~$\mathbf{L}^*$ be a 
proper split $F$-stable Levi subgroup of~$\mathbf{G}^*$, and let~$\mathbf{L}$
be an $F$-stable Levi subgroup of~$\mathbf{G}$ dual to~$\mathbf{L}^*$.
Let~$\vartheta$ be a character of~$L$ such that $R_L^G( \vartheta ) = \chi$. 
We may assume that $s \in L^*$ and that $\vartheta \in \mathcal{E}( L, [s] )$ 
(see \cite[Proposition 15.7]{CaEn}). Using the notation and the statements of 
Theorem~\ref{SameComponentGroupConverse}, we find that $c(\vartheta) \geq 
c(\chi)$. Suppose first that $c(\vartheta) = c(\chi)$. Then 
$C^\circ_{\mathbf{G}^*}(s) \leq \mathbf{L}^*$ by 
Theorem~\ref{SameComponentGroupConverse}, and we conclude with 
Lemmas~\ref{CentralizerGeneration} and~\ref{ComparingCentralizers} that 
$C^\circ_{\check{\mathbf{G}}^*}( \check{s} ) = 
(j^*)^{-1}( C^\circ_{\mathbf{G}^*}(s) )$ is contained in a proper split 
$F$-stable Levi subgroup of~$\check{\mathbf{G}}^*$. As~$\check{s}$ is not as
in Lemma~\ref{CentralizerContainedInSplitLeviConformal}(a)(i.2), we have
arrived at a contradiction.

Hence $c(\vartheta) > c(\chi)$. Now $c(\chi) \leq |A_{\mathbf{G}^*}(s)^F| \leq 
|Z(\mathbf{G})/Z^\circ(\mathbf{G})| \leq 4$, where the first inequality follows 
from~(\ref{OrbitLength}), and by Corollary~\ref{CorollaryToWeylGroupInduction} 
and Theorem~\ref{SameComponentGroupConverse} we cannot have $c(\vartheta)/c(\chi) 
= 2$. Thus $c(\chi) = 1$ and $c(\vartheta) = 4$. As $c(\vartheta) \leq 
|A_{\mathbf{L}^*}(s)^F| \leq |A_{\mathbf{G}^*}(s)^F|$, we have in particular 
that $|A_{\mathbf{G}^*}(s)^F| = 4$. Applying Lemma~\ref{CentralizerGeneration},
we find that $\tilde{A}_{\check{\mathbf{G}}^*}(\check{s})^F$ has order~$4$,
so that~$\check{s}$ is conjugate to $-\check{s}$ in~$\check{G}$.
By~(\ref{OrbitLength}) we have $1 = c(\chi) = 
|A_{\mathbf{G}^*}(s)_\lambda^F|$, i.e.\ the orbit of~$\lambda$ under 
$A_{\mathbf{G}^*}(s)^F$ has length~$4$. 
Transferring the situation to~$\check{\mathbf{G}}^*$ with the help of
Lemma~\ref{CentralizerGeneration}, we see that the orbit of~$\check{\lambda}$ 
under $\tilde{A}_{\check{\mathbf{G}}^*}(\check{s})^F$ has length~$4$.
With the notation introduced in~\ref{CriticalSemisimpleElements}, put
$\check{\mathbf{G}}^*_{\tilde{\mu}} := \check{\mathbf{G}}^*_{\tilde{\mu}}( \check{s} )$ 
and $\check{C}_{\tilde{\mu}} := 
C^\circ_{\check{\mathbf{G}}^*_{\tilde{\mu}}}( \check{s}_{\tilde{\mu}} )^F$, where
$\mu \in \mathcal{F}_{\check{s}}$, and $\tilde{\mu} = \mu$ if $\mu = \mu'$,
and $\tilde{\mu} = \mu\mu'$, otherwise. 
Then $C^\circ_{\check{\mathbf{G}}^*}( \check{s} )^F$ is the direct product of
the groups $\check{C}_{\tilde{\mu}}$, and $\check{\lambda}$ is the outer
tensor product of unipotent characters $\check{\lambda}_{\tilde{\mu}}$ of
$\check{C}_{\tilde{\mu}}$. 
Moreover, $\tilde{A}_{\check{\mathbf{G}}^*}(\check{s})^F$ fixes each of the
groups $\check{C}_{\tilde{\mu}}$.
By Poposition~\ref{CorollaryElementsInConformalGroups}, the characters
$\check{\lambda}_{\tilde{\mu}}$ lie in 
$\tilde{A}_{\check{\mathbf{G}}^*}(\check{s})^F$ orbits of lenghts~$1$,~$2$ 
or~$4$. As~$\check{\lambda}$ lies in an orbit of length~$4$, one of the
$\check{\lambda}_{\tilde{\mu}}$ must lie in an orbit of length~$4$. 
We are thus in the following situation: $F = F'$ and $\mathcal{F}_{\check{s}}$ 
contains an element $\nu = X - \zeta$ of Type~(I), which occurs with even 
multiplicity $k \geq 2$ in the characteristic polynomial of~$\check{s}$
(in particular, $\alpha = \zeta^2$ is a square in $\mathbb{F}_q$),
and the restriction of~$\check{\lambda}_{\tilde{\nu}}$ to 
$C^\circ_{[\check{\mathbf{G}}^*_{\tilde{\nu}},\check{\mathbf{G}}^*_{\tilde{\nu}}]}
( \check{s}_{\tilde{\nu}} )^F \cong \SO^+_k(q) \times \SO^+_k( q )$ is labelled by a 
pair of symbols $(\Lambda_1,\Lambda_2)$, where $\Lambda_1$ and $\Lambda_2$ are 
both degenerate and label two different unipotent characters of $\SO_k^+( q )$ 
(in particular, $k \geq 4$). Here, of course, the pair $(\Lambda_1,\Lambda_2)$
labels the outer tensor product of the characters of $\SO_k^+( q )$
corresponding to~$\Lambda_1$ and $\Lambda_2$, respectively.

As $c(\vartheta)/c(\chi) = 4$, we conclude from 
Theorem~\ref{SameComponentGroupConverse} that there is an irreducible 
constituent $\tilde{\vartheta} \in \mathcal{E}( \tilde{L}, [\tilde{s}] )$ of
$\Ind_L^{\tilde{L}}( \vartheta )$ such that 
$R_{\tilde{L}}^{\tilde{G}}( \tilde{\vartheta} )$ is a sum of four irreducible 
characters of equal degrees. Now put 
$C_{\tilde{G}} := C_{\tilde{\mathbf{G}}^*}( \tilde{s} )^F$ and
$C_{\tilde{L}} := C_{\tilde{\mathbf{L}}^*}( \tilde{s} )^F$
Let $\tilde{\kappa} \in \mathcal{E}( C_{\tilde{L}}, [1] )$ correspond 
to~$\tilde{\vartheta}$ under Lusztig's Jordan decomposition of characters. 
The latter commutes with Harish-Chandra induction (see 
\cite[p.~$1049$--$1050$]{En08}) 
and thus $R_{C_{\tilde{L}}}^{C_{\tilde{G}}}( \tilde{\kappa} )$
is a sum of four irreducible characters, and all of these have the same degree.
Now $k^*(C_{\tilde{G}}) = C^\circ_{\check{\mathbf{G}}^*}( \check{s} )^F =: 
C_{\check{G}}$ by Lemma~\ref{CentralizerGeneration}.
Similarly, $k^*( C_{\tilde{L}} ) = C^\circ_{\check{\mathbf{L}}^*}( \check{s} )^F 
=: C_{\check{L}}$. Let~$\check{\kappa}$ denote the unipotent character of
$C_{\check{L}}$ such that $\tilde{\kappa}$ is the inflation of $\check{\kappa}$
over the kernel of $k^*|_{\tilde{G}}$. As the unipotent characters 
of~$\tilde{G}$ and~$C_{\tilde{G}}$ have this kernel in their kernels, we find 
that $R_{C_{\check{L}}}^{C_{\check{G}}}( \check{\kappa} )$ is a sum of four 
unipotent characters, one of which is~$\check{\lambda}$. As $c(\vartheta) = 4$, 
we conclude from~(\ref{OrbitLength}) that~$\check{\kappa}$ is invariant under 
the action of $\tilde{A}_{\check{\mathbf{L}}^*}( \check{s} )^F = 
\tilde{A}_{\check{\mathbf{G}}^*}( \check{s} )^F$. It follows that the four
constituents of $R_{C_{\check{L}}}^{C_{\check{G}}}( \check{\kappa} )$ are 
exactly the four characters in the 
$\tilde{A}_{\check{\mathbf{G}}^*}( \check{s} )^F$-orbit of~$\check{\lambda}$. 
Now $C_{\check{L}}$ decomposes
into direct summands $C_{\check{L},\tilde{\mu}}$, and each of these is a
Levi subgroup of $\check{C}_{\tilde{\mu}}$. As Harish-Chandra induction is
compatible with this direct product decomposition, we conclude that for 
each~$\tilde{\mu}$, the orbit sums of the $\check{\lambda}_{\tilde{\mu}}$ are 
Harish-Chandra induced. This applies in particular for $\tilde{\nu} = 
X^2 - \alpha$. To show that the latter is impossible, we may assume that 
$\mathcal{F}_{\check{s}} = \{ X - \zeta, X + \zeta \}$.

Assuming this, we put $\mathbf{C} := 
C^\circ_{\check{\mathbf{G}}^*}( \check{s} )$. We then have two embeddings
$\mathbf{C}_1' \times \mathbf{C}_2' \stackrel{\gamma}{\rightarrow} \mathbf{C} 
\stackrel{\delta}{\rightarrow} \mathbf{C}_1 \times \mathbf{C}_2$, where 
$\mathbf{C}_1$ and
$\mathbf{C}_2$ denote the special conformal orthogonal groups acting on
the $\zeta$ eigenspace and the $-\zeta$ eigenspace, respectively, 
of~$\check{s}$. Also, $\mathbf{C}_i'$ denotes the set of elements of
$\mathbf{C}_i$ with multipliers~$1$, $i = 1, 2$, and 
$\mathbf{C} = \{ (g_1,g_2) \in \mathbf{C}_1 \times \mathbf{C}_2 \mid \alpha_{g_1} = 
\alpha_{g_2} \}$. As $\mathbf{C}_1 \times \mathbf{C}_2$ has connected center,
the maps $\delta$ and $\delta \circ \gamma$ are regular embeddings. Thus
the Levi subgroups of $\mathbf{C}$ are the intersections of the Levi subgroups
of $\mathbf{C}_1 \times \mathbf{C}_2$ with~$\mathbf{C}$, and similarly for
the Levi subgroups of $\mathbf{C}_1' \times \mathbf{C}_2'$. Now 
let~$\mathbf{M}$ be an $F$-stable split Levi subgroup of~$\mathbf{C}$.
By the above considerations, $\mathbf{M}' := \mathbf{M} \cap (\mathbf{C}_1'
\times \mathbf{C}_2')$ is a split $F$-stable Levi subgroup of $\mathbf{C}_1'
\times \mathbf{C}_2'$. We also find that $C = (C_1' \times C_2')M$, as the Levi
subgroups of~$\mathbf{C_i}$ contain elements with arbitrary multipliers.
It follows that $\Res^C_{C_1'\times C_2'}( R_M^C( \beta ) )
= R_{M'}^{C_1' \times C_2'}( \Res_{M'}^M( \beta ) )$ for $\beta \in \Irr(M)$. 
Applying this to our situation above, we find that the 
$\tilde{A}_{\check{\mathbf{G}}^*}( \check{s} )^F$-orbit sum of the restriction 
of~$\check{\lambda}$ to $\SO^+_k(q) \times \SO^+_k(q)$ is Harish-Chandra 
induced. As this orbit contains the four different characters labelled by
$(\Lambda_1,\Lambda_2), (\Lambda_2,\Lambda_1), (\Lambda_1',\Lambda_2'), 
(\Lambda_2',\Lambda_1')$, this is impossible, a contradiction.

(b) The proof is the same as the one of Theorem~\ref{RestrictionToSpinOdd}(b).

(c) Let us begin with some preliminary remarks.
Lemma~\ref{CentralizerContainedInSplitLeviConformal} implies that
there is a proper split $F$-stable Levi subgroup $\check{\mathbf{L}}^*$ of
$\check{\mathbf{G}}^*$ containing $C^\circ_{\check{\mathbf{G}}^*}( \check{s} )$,
and even $C_{\check{\mathbf{G}}^*}( \check{s} )$ if there is $\mu \in
\mathcal{F}_{\check{s}}$ with $\mu \neq \mu^{*\alpha}$ or if~$\check{s}$ is
exceptional. Putting $\mathbf{L}^* := j^*( \check{\mathbf{L}}^* )$, we 
obtain $C^\circ_{\mathbf{G}^*}( s ) = 
j^*( C^\circ_{\check{\mathbf{G}}^*}( \check{s} ) ) \leq \mathbf{L}^*$ from
Lemma~\ref{CentralizerGeneration}. This implies in particular 
$C^\circ_{\mathbf{G}^*}(s) = C^\circ_{\mathbf{G}^*}(s) \cap \mathbf{L}^*
= C^\circ_{\mathbf{L}^*}(s)$ (see \ref{TheComponentGroup}), a fact that will be
used throughout the proof.
As~$\mathbf{L}^*$ is a proper split $F$-stable Levi subgroup of~$\mathbf{G}^*$
by Lemma~\ref{ComparingCentralizers}, our claim about~$\chi$ follows from
Corollary~\ref{SufficientImprimitivityCondition}.

(i) Suppose that~$\check{s}$ is not conjugate to $-\check{s}$ in~$\check{G}^*$. 
Then $\tilde{C}_{\check{\mathbf{G}}^*}( \check{s} ) = 
C_{\check{\mathbf{G}}^*}( \check{s} )$. Assume first that there is $\mu \in 
\mathcal{F}_{\check{s}}$ with $\mu \neq \mu^{*\alpha}$ or that~$\check{s}$ is 
exceptional. By Lemma~\ref{CentralizerGeneration} and the preliminary remarks, 
we conclude that $C_{\mathbf{G}^*}(s) = 
{j}^*( \tilde{C}_{\check{\mathbf{G}}^*}( \check{s} ) ) = 
{j}^*( C_{\check{\mathbf{G}}^*}( \check{s} ) ) \leq \mathbf{L}^*$. 
In particular, $A_{\mathbf{G}^*}( s )^F_\lambda \leq 
A_{\mathbf{G}^*}( s )^F = A_{\mathbf{L}^*}( s )^F$, proving our claim.

Now assume that~$\check{s}$ is as in
Lemmas~\ref{CentralizerContainedInSplitLeviConformal}(a)(i.2) and 
that~$\check{s}$ is not exceptional. In particular, $X^2 - \alpha$ divides the
minimal polynomial of~$\check{s}$. It follows from
Lemma~\ref{NonConnectedCentralizersInCSO} that
$\tilde{A}_{\check{\mathbf{G}}^*}(\check{s})^F$ has order~$2$ and acts as a graph 
automorphism on the components of $C^\circ_{\check{\mathbf{G}}^*}( \check{s} )^F$ 
corresponding to the Type~(I) or Type~(II) elements of $\mathcal{F}_{\check{s}}$, 
and acts trivially on the other components of this centralizer. From this it is 
easy to see (cf.\ Remark~\ref{RemarkForEvenDimensionalSpinGroups}) that 
$\tilde{A}_{\check{\mathbf{G}}^*}(\check{s})^F =
\tilde{A}_{\check{\mathbf{G}}^*}(\check{s})_{\lambda}^F$, unless~$\check{s}$ 
and~$\check{\lambda}$ are as in the statement, in which case 
$A_{\mathbf{G}^*}(s)_{\lambda}^F$ is trivial.
Corollary~\ref{SufficientImprimitivityCondition} implies that~$\chi$ is Harish 
Chandra induced from~$L$ if $A_{\mathbf{G}^*}(s)_{\lambda}^F$ is trivial.
Otherwise,~$\chi$ is Harish-Chandra primitive by 
Corollary~\ref{SufficientPrimitivityConditionNew} and
Lemma~\ref{CentralizerContainedInSplitLeviConformal}(a)(i.2). This completes 
the proof of~(i).

(ii) Now suppose that~$\check{s}$ is not exceptional, that~$\check{s}$ is 
conjugate to~$-\check{s}$ in~$\check{G}^*$ and that $X^2 - \alpha$ does not 
divide the minimal polynomial of~$\check{s}$. This implies, first of all,
that~$\check{s}$ is not as in
Lemma~\ref{CentralizerContainedInSplitLeviConformal}(a)(i.2).
Also, $C^\circ_{\check{\mathbf{G}}^*}( \check{s} ) = 
C_{\check{\mathbf{G}}^*}( \check{s} ) \lneq 
\tilde{C}_{\check{\mathbf{G}}^*}( \check{s} )$ by 
Lemma~\ref{NonConnectedCentralizersInCSO} and the definition 
of~$\tilde{C}_{\check{\mathbf{G}}^*}( \check{s} )$ in~(\ref{TildeC}). Thus 
$|A_{\mathbf{G}^*}( s )^F| = 2$ by Lemma~\ref{CentralizerGeneration}.
Lemmas~\ref{CentralizerContainedInSplitLeviConformal}(a)(ii) 
and~\ref{ComparingCentralizers} imply that 
$C_{\mathbf{G}^*}( s )^F =
j^*( \tilde{C}_{\check{\mathbf{G}}^*}( \check{s} )^F )$ 
is not contained in any proper split $F$-stable 
Levi subgroup of~$\mathbf{G}^*$. In particular, $|A_{\mathbf{L}^*}( s )^F| = 1$. 
The last two assertions follow from 
Corollary~\ref{SufficientImprimitivityCondition}, respectively
Corollary~\ref{SufficientPrimitivityConditionNew}.

(iii) Now suppose that~$\check{s}$ is conjugate to~$-\check{s}$ 
in $\check{G}^*$ and that $X^2 - \alpha$ divides the minimal polynomial 
of~$\check{s}$. Then $|A_{\mathbf{G}^*}( s )^F| =
|\tilde{A}_{\check{\mathbf{G}}^*}(\check{s})^F| = 4$ by
Lemma~\ref{CentralizerGeneration} and
Proposition~\ref{CorollaryElementsInConformalGroups}.

Suppose first that $\mu = \mu^{*\alpha}$ for
all $\mu \in \mathcal{F}_{\check{s}}$, i.e.\ that~$\check{s}$ is as in
Lemma~\ref{CentralizerContainedInSplitLeviConformal}(a)(i.2). 
As the multiplicity of~$\nu$ in the characteristic polynomial of~$\check{s}$
is~$2$, the orbits of $A_{\mathbf{G}^*}( s )^F$ on
$\mathcal{E}( C^\circ_{\mathbf{G}^*}( s )^F, [1] )$ have lengths~$1$ or~$2$ 
by Proposition~\ref{CorollaryElementsInConformalGroups}.
Hence $|A_{\mathbf{G}^*}( s )_\lambda^F| \in \{ 4, 2 \}$.  
Also,~$\check{\mathbf{L}}^*$ is the stabilizer of two complementary, totally 
isotropic subspaces of~$\mathbf{V}_\nu( s )$ whose sum 
equals~$\mathbf{V}_\nu( s )$.
If~$\nu$ is of Type~(I), this easily implies that 
$C_{\check{\mathbf{L}}^*}( \check{s} ) = 
C^\circ_{\check{\mathbf{G}}^*}( \check{s} )$, hence $A_{\mathbf{L}^*}(s)^F$ is 
trivial by Lemma~\ref{CentralizerGeneration}. In particular,
$A_{\mathbf{G}^*}( s )_\lambda^F \lneq A_{\mathbf{L}^*}( s )^F$, proving
our assertion. If~$\nu$ has Type~(II), then
$\tilde{A}_{\check{\mathbf{L}}^*}( \check{s} )^F$ has order~$2$ and acts by 
conjugating~$\check{s}$ to~$-\check{s}$. 
%%%%%%%%%%%%%%%%%%%%%%%%%%%%%%%%%%%%%%%%%%%%%%%%%%%%%%%%%%%%%%%%%%%%%%%%%%%%%%%%
%%%%%%%%%%%%%%%%%%%%%%%%%%%%%%%%%%%%%%%%%%%%%%%%%%%%%%%%%%%%%%%%%%%%%%%%%%%%%%%%
%%
%% Für mehr Details zu der Aktion von $\tilde{A}_{\check{\mathbf{L}}^*}( \check{s} )^F$
%% siehe Notizen vom 14.11.2016, Seite 3.
%%
%%%%%%%%%%%%%%%%%%%%%%%%%%%%%%%%%%%%%%%%%%%%%%%%%%%%%%%%%%%%%%%%%%%%%%%%%%%%%%%%
%%%%%%%%%%%%%%%%%%%%%%%%%%%%%%%%%%%%%%%%%%%%%%%%%%%%%%%%%%%%%%%%%%%%%%%%%%%%%%%%
Now if 
$A_{\mathbf{G}^*}( s )_\lambda^F \leq A_{\mathbf{L}^*}( s )^F$, 
an element of $\check{G}^*$ conjugating~$\check{s}$ to~$-\check{s}$ would
stabilize~$\check{\lambda}$. But then 
$|A_{\mathbf{G}^*}( s )_\lambda^F| = 
|\tilde{A}_{\check{\mathbf{G}}^*}( \check{s} )^F_{\check{\lambda}}| = 4$ by 
Proposition~\ref{CorollaryElementsInConformalGroups}, a contradiction.

Now suppose that there is $\mu \in \mathcal{F}_{\check{s}}$ such that
$\mu \neq \mu^{*\alpha}$. By the preliminary remark, 
$C_{\check{\mathbf{G}}^*}( \check{s} ) \leq \check{\mathbf{L}}^*$. We have
$|(C_{\check{\mathbf{G}}^*}( \check{s} )/C^\circ_{\check{\mathbf{G}}^*}( \check{s} ))^F| = 2$,
as $|\tilde{A}_{\check{\mathbf{G}}^*}(\check{s})^F| = 4$.
In turn, 
$|(j^*( C_{\check{\mathbf{G}}^*}( \check{s} ) )/C^\circ_{\mathbf{G}^*}(s))^F| = 2$.
Now $j^*( C_{\check{\mathbf{G}}^*}( \check{s} ) ) \leq 
C_{\mathbf{G}^*}(s) \cap \mathbf{L}^* = C_{\mathbf{L}^*}( s )$, and
thus $A_{\mathbf{L}^*}( s )^F$ is nontrivial. We 
cannot have $|A_{\mathbf{L}^*}( s )^F| = 4$ by 
Lemma~\ref{CentralizerContainedInSplitLeviConformal}(a)(ii),
hence $|A_{\mathbf{L}^*}( s )^F| = 2$ as claimed, and we may assume that 
$b \in A_{\mathbf{L}^*}( s )^F$. 

(iii.1) If~$m$ is odd, $|A_{\mathbf{G}^*}( s )^F|$ is cyclic by
Proposition~\ref{CorollaryElementsInConformalGroups}. In particular,
$A_{\mathbf{G}^*}( s )^F_\lambda \leq A_{\mathbf{L}^*}( s )^F$, if and only if
$A_{\mathbf{G}^*}( s )^F_\lambda \lneq A_{\mathbf{G}^*}( s )^F$.

(iii.2) Now suppose that $m$ is even.
The structure of $A_{\mathbf{G}^*}( s )$ and the description of its elements
follow from Proposition~\ref{CorollaryElementsInConformalGroups}. We have
$A_{\mathbf{G}^*}( s )^F_\lambda \leq A_{\mathbf{L}^*}( s )^F$, if and only
if~$a$ does not fix~$\lambda$.

(d) The assertion has already been proved for elements~$\check{s}$ as in~(c)(i) 
or~(c)(ii). Thus let~$\check{s}$ is as in~(c)(iii) and assume that~$\chi$ is not
Harish-Chandra induced from~$L$. Then $|A_{\mathbf{G}^*}( s )^F_\lambda| > 1$.
If $|A_{\mathbf{G}^*}( s )^F_\lambda| = 4$, our assertion follows from
Corollary~\ref{SufficientPrimitivityConditionNew}(b), as $C_{\mathbf{G}^*}(s)^F$
is not contained in any proper split $F$-stable Levi subgroup of~$\mathbf{G}^*$ 
by Lemmas~\ref{CentralizerContainedInSplitLeviConformal}(a),~\ref{CentralizerGeneration}
and~\ref{ComparingCentralizers}. Hence $|A_{\mathbf{G}^*}( s )^F_\lambda| = 2$,
and by~(iii.2) we have $A_{\mathbf{G}^*}( s )^F_\lambda = \langle a \rangle$ or
$A_{\mathbf{G}^*}( s )^F_\lambda = \langle ab \rangle$. 
Proposition~\ref{CorollaryElementsInConformalGroups} excludes the case 
that~$\check{s}$ is as in 
Lemma~\ref{CentralizerContainedInSplitLeviConformal}(a)(i.2) and $\nu$ is of 
Type~(II) (recall that the multiplicity of~$\nu$ in the characteristic polynomial 
of~$\check{s}$ in this case is~$2$, and that $\SO_2^+( q^2 )$ has only one
unipotent character).

Suppose that~$\chi$ is Harish-Chandra induced from $M = \mathbf{M}^F$, 
where~$\mathbf{M}$ is a proper split Levi subgroup of~$\mathbf{G}$. 
Let~$\mathbf{M}^*$ be an $F$-stable Levi subgroup of~$\mathbf{G}^*$ dual 
to~$\mathbf{M}$. Let~$\vartheta$ be an irreducible character of~$M$ such that 
$\chi = R_M^G( \vartheta )$.
By Theorem~\ref{SameComponentGroupConverse} (with $\tilde{\mathbf{L}}$ replaced 
by $\tilde{\mathbf{M}}$), we have $c(\chi) \leq c(\vartheta)$.
Now~$c(\chi) = |A_{\mathbf{G}^*}( s )^F_\lambda| = 2$ by~(\ref{OrbitLength}), 
and $c(\vartheta) \leq |A_{\mathbf{M}^*}( s )^F|$, where 
$|A_{\mathbf{M}^*}( s )^F|$ divides 
$|A_{\mathbf{G}^*}( s )^F| = 4$. 

If $c(\chi) = c(\vartheta)$, 
Theorem~\ref{SameComponentGroupConverse} implies that 
$C^\circ_{\mathbf{G}^*}( s ) \leq \mathbf{M}^*$. But then
$A_{\mathbf{G}^*}( s )^F_\lambda \leq A_{\mathbf{M}^*}( s )^F$ by
Corollary~\ref{SufficientImprimitivityCondition}. As
$A_{\mathbf{G}^*}( s )^F_\lambda$ contains~$a$ or~$ab$, it follows that
$(j^*)^{-1}( C_{\mathbf{G}^*}( s )_\lambda^F ) \leq (j^*)^{-1}( \mathbf{M}^* )$ 
contains an element that conjugates~$\check{s}$ to~$-\check{s}$. As 
$(j^*)^{-1}( \mathbf{M}^* )$ is a proper split $F$-stable Levi subgroup 
of~$\check{\mathbf{G}}^*$ by Lemma~\ref{ComparingCentralizers}, this
contradicts Lemma~\ref{CentralizerContainedInSplitLeviConformal}(a)(ii).

Thus suppose that $c(\vartheta) = 4$. By 
Theorem~\ref{SameComponentGroupConverse}, there is an irreducible 
character~$\tilde{\vartheta}$ of~$\tilde{M}$ such that 
$R_{\tilde{M}}^{\tilde{G}}( \tilde{\vartheta} )$ has exactly two 
irreducible constituents of the same degree. This contradicts 
Corollary~\ref{CorollaryToWeylGroupInduction}(a).

Condition~(\ref{InclusionCondition}) of Theorem~\ref{MainResult} is trivially 
satisfied under the hypothesis on~$\check{s}$ in~(b). We have just shown that 
if~$\check{s}$ is as in~(c), then either Condition~(\ref{InclusionCondition}) 
is satisfied for~$\mathbf{L}^*$, or~$\chi$ is Harish-Chandra primitive. Thus 
Theorem~\ref{MainResult}(b) holds for~$G$.
\end{prf}

\smallskip
\noindent
We summarize the conditions for Harish-Chandra imprimitivity of an
irreducible character of $\Spin_{2m}^\pm(q)$.

\begin{cor}
Let the notation be as in {\rm Theorem~\ref{RestrictionToSpinEven}} and 
{\rm Tables~\ref{ConjugateToNegativeOrthogonalA},~\ref{ConjugateToNegativeOrthogonalB}}
and~{\rm \ref{ConjugateToNegativeOrthogonalC}}.
Then~$\chi$ is Harish-Chandra imprimitive exactly in the following cases.

{\rm (a)} For all $\mu \in \mathcal{F}_{\check{s}}$ we have $\mu = 
\mu^{*\alpha}$, there is $\zeta \in \mathbb{F}_q$ with $\zeta^2 = \alpha$ and
$X - \zeta$, respectively $X + \zeta$, occur with multiplicity~$2$, 
respectively~$4 k'$, in the characteristic polynomial of~$\check{s}$. 
Moreover, the eigenspaces $V_{X \pm \zeta}( \check{s} )$ have maximal Witt index
and if~$k' > 0$, the component of $\check{\lambda}$ corresponding 
to~$V_{X + \zeta}( \check{s} )$ is labelled by a degenerate symbol.

{\rm(b)} There exists $\mu \in \mathcal{F}_{\check{s}}$ with $\mu \neq 
\mu^{*\alpha} \neq \mu'$.

{\rm (c)} Every $\mu \in \mathcal{F}_{\check{s}}$ satisfies $\mu = \mu^{*\alpha}$
or $\mu' = \mu^{*\alpha}$, there exists $\mu \in \mathcal{F}_{\check{s}}$
with $\mu \neq \mu^{*\alpha}$, and if $\check{s}$ is conjugate to $-\check{s}$ 
in $\check{G}^*$, the following conditions are satisfied: There is~$j$
such that $\kappa_{1,j} \neq \kappa_{2,j}$ or~$m$ is even and $\Lambda_1 \neq 
\Lambda_2, \Lambda_2'$ or~$m$ is odd and $\Lambda_1 \neq \Lambda_2$.
\end{cor}

\subsection{The exceptional groups} \label{TheExceptionalGroups}
Here, we present the results for the 
quasisimple exceptional groups of Lie type arising from algebraic groups with
non-connected centers. These are the groups $E_6( \mathbb{F} )_{\text{sc}}$
for $\text{char}( \mathbb{F} ) \neq 3$ and $E_7( \mathbb{F} )_{\text{sc}}$
for $\text{char}( \mathbb{F} ) \neq 2$. Let~$\mathbf{G}$ denote one of these
groups. Then $\mathbf{G}^* = E_6( \mathbb{F} )_{\text{ad}}$, respectively
$\mathbf{G}^* = E_7( \mathbb{F} )_{\text{ad}}$. In the first
case, $|Z(\mathbf{G})| = 3$ and in the second case, $|Z(\mathbf{G})| = 2$.
If $\mathbf{G} = E_6( \mathbb{F} )_{\text{sc}}$, let~$F'$ and~$F''$ denote 
Frobenius morphisms of~$\mathbf{G}$ such that $\mathbf{G}^{F'} = 
E_6(q)_{\text{sc}}$ and $\mathbf{G}^{F''} = {^2\!E}_6(q)_{\text{sc}}$. Let~$F$ 
be one of~$F'$ or~$F''$. Then $G^* = E_6(q)_{\text{ad}}$ if $F = F'$, and 
$G^* = {^2\!E}_6(q)_{\text{ad}}$, otherwise. If $\mathbf{G} = 
E_7( \mathbb{F} )_{\text{sc}}$, let~$F$ denote a Frobenius morphism
of~$\mathbf{G}$ such that $\mathbf{G}^{F} = E_7(q)_{\text{sc}}$. Then 
$G^* = E_7(q)_{\text{ad}}$.

In our situation Theorem~\ref{MainResult}(b) holds by 
Corollary~\ref{ProofOfTheorem11bInSpecialCaseNew}. We thus have to decide,
for semisimple elements $s \in G^*$, the containment of
$C_{\mathbf{G}^*}(s)_\lambda^F C^\circ_{\mathbf{G}^*}(s)$ in proper split 
$F$-stable Levi subgroups of~$\mathbf{G}^*$ (in the notation of 
Theorem~\ref{MainResult}). For the purpose of this investigation we introduce 
one further piece of 
notation and recall some facts from a paper by Brou{\'e} and Malle~\cite{BrouMa}.
For a positive integer~$i$, we let $\Phi_i \in \mathbb{Z}[X]$ denote the~$i$th 
cyclotomic polynomial. In \cite[D{\'e}finition~$1.9$, Lemma~$3.1$]{BrouMa},
Brou{\'e} and Malle associate to an $F$-stable torus~$\mathbf{T}^*$ 
of~$\mathbf{G}^*$ an order polynomial $f \in \mathbb{Z}[X]$, such that 
$|{(\mathbf{T}^*)}^F| = f(q)$. Moreover,~$f$ is a product of $\Phi_i$'s for 
certain values of~$i$, and $\mathbf{T}^*$ contains a nontrivial split 
$F$-stable subtorus, if and only if~$\Phi_1$ divides~$f$. Suppose that 
$\mathbf{H}^*$ is a closed $F$-stable subgroup of~$\mathbf{G}^*$ satisfying
$C_{\mathbf{G}^*}( \mathbf{H}^* ) \leq \mathbf{H}^*$. (The latter condition
is satisfied if~$\mathbf{H}^*$ contains a maximal torus of~$\mathbf{G}^*$.)
Then $\mathbf{H}^*$ is contained
in a proper split $F$-stable Levi subgroup of~$\mathbf{G}^*$, if and only if
the order polynomial of $Z^\circ( \mathbf{H}^* )$ is divisible by~$\Phi_1$.
Indeed, $Z^\circ( \mathbf{H}^* )$ is an $F$-stable torus of~$\mathbf{G}^*$,
and $\mathbf{H}^*$ is contained in a proper split $F$-stable Levi subgroup 
of~$\mathbf{G}^*$, if and only if $\mathbf{H}^*$ centralizes a nontrivial 
split $F$-stable torus.

In the following remark we introduce the cases that may arise, the notation and
labels used in the tables below.

\begin{rem}\label{ExceptionalGroupsRemark}
{\rm
Let $s \in G^*$ be semisimple and let $\chi \in \mathcal{E}( G, [s] )$. Let
$\lambda \in \mathcal{E}( C^\circ_{\mathbf{G}^*}(s)^F, [1] )$, such 
that~$[\chi]$ corresponds to $[\lambda]$ under Lusztig's generalized Jordan 
decomposition of characters. Put $\mathbf{C} := C_{\mathbf{G}^*}( s )$ and 
$\mathbf{Z} := Z( \mathbf{C} )$. The following cases arise.

Case~$\ast$: Here,~$\mathbf{C}^\circ$ is not contained in any proper split 
$F$-stable Levi subgroup of~$\mathbf{G}^*$. By the remarks above, this happens 
if and only if the order polynomial of $Z^\circ( \mathbf{C}^\circ )$ is not 
divisible by $\Phi_1$. In this case, every element of $\mathcal{E}( G, [s] )$ 
is Harish-Chandra primitive by Corollary~\ref{ThirdLusztigSeriesResult}.  

Case~$\dagger$: Here,~$\mathbf{C}^\circ$ is contained in some proper split 
$F$-stable Levi subgroup of~$\mathbf{G}^*$ but~$\mathbf{C}$ is not contained in 
any such subgroup. This is the case if and only if the order polynomial of 
$Z^\circ( \mathbf{C}^\circ )$ is divisible by $\Phi_1$, while the order 
polynomial of $\mathbf{Z}^\circ$ is not divisible by~$\Phi_1$. In this case, all 
elements of $\mathcal{E}( G, [s] )$ are Harish-Chandra primitive if 
$|A_{\mathbf{G}^*}(s)_\lambda^F| \neq 1$. Otherwise, all elements of 
$\mathcal{E}( G, [s] )$ are Harish-Chandra imprimitive. Indeed, if 
$|A_{\mathbf{G}^*}(s)_\lambda^F| = 1$, then $\mathbf{C}_\lambda^F
\leq \mathbf{C}^\circ$, and the claim follows from Theorem~\ref{MainResult}(a). 
Conversely, suppose that $\chi$ is Harish-Chandra imprimitive. Then, 
by Theorem~\ref{MainResult}(b), there is some proper split $F$-stable Levi 
subgroup $\mathbf{L}^* \leq \mathbf{G}^*$ such that 
$\mathbf{C}_\lambda^F\mathbf{C}^\circ \leq \mathbf{L}^*$. Aiming at a 
contradiction, we assume that $|A_{\mathbf{G}^*}(s)_\lambda^F| \neq 1$. Then 
$A_{\mathbf{G}^*}(s)^F = A_{\mathbf{G}^*}(s)$, and thus $\mathbf{C} = 
\mathbf{C}^F\mathbf{C}^\circ = \mathbf{C}^F_\lambda\mathbf{C}^\circ$, as every 
coset of $\mathbf{C}/\mathbf{C}^\circ$ contains an $F$-stable element. This is 
the desired contradiction.

Case~\Checkmark: Here,~$\mathbf{C}$ is contained in some proper split $F$-stable 
Levi subgroup of~$\mathbf{G}^*$. This is the case, if and only if the order
polynomial of $\mathbf{Z}^\circ$ is divisible by~$\Phi_1$. In this case, every
element of $\mathcal{E}( G, [s] )$ is Harish-Chandra imprimitive by
\cite[Theorem~$7.3$]{HiHuMa}.
}
\end{rem}

Our results for the groups~$G$ considered here rely on a classification of the 
semisimple class 
types of~$G^*$ and their centralizers. (Recall that two semisimple elements 
of~$G^*$ belong to the same class type, if their centralizers in~$\mathbf{G}^*$ 
are conjugate in the finite group~$G^*$.) In the adjoint case we are considering, 
this classification is due to Frank L{\"u}beck and is given in the Tables 
of~\cite{LL}. There, each class type is labelled by a triple of natural numbers. 
Suppose that $s, s' \in G^*$ represent class types. Put $\mathbf{C} := 
C_{\mathbf{G}^*}( s )$ and $\mathbf{C}' := C_{\mathbf{G}^*}( s' )$. Then the 
first index in the triples labelling the class types of~$s$ and~$s'$ are equal, 
if and only if~$\mathbf{C}^\circ$ and~${(\mathbf{C}')}^\circ$ are conjugate 
in~$\mathbf{G}^*$. The first two indices of these triples are equal, if and 
only if~$\mathbf{C}$ and~$\mathbf{C}'$ are conjugate in~$\mathbf{G}^*$. 
Put~$\mathbf{Z} := Z(\mathbf{C})$. Then the entry corresponding to the class 
type of~$s$ gives~$|\mathbf{Z}^F|$, the latter in a form reflecting the order 
polynomial of the torus $\mathbf{Z}^\circ$: the $\varphi_j$ in the table 
(there written as {\tt phi1, phi2, phi3,} \ldots) stands 
for $\Phi_j(q)$, and $\varphi_j^b$ is given as a factor in the entry for 
$|\mathbf{Z}^F|$, if and only if $\Phi_j^b$ divides the order polynomial 
of~$\mathbf{Z}^\circ$. (This follows from the fact that the order formulae 
in~\cite{LL} for ${(\mathbf{Z}^\circ)}^F$ are valid for all~$q$; they are in 
particular valid for all powers of~$F$.) The extra factor for $|\mathbf{Z}^F|$, 
if present, gives $|(\mathbf{Z}/\mathbf{Z}^\circ)^F|$. In 
particular,~$\mathbf{Z}^\circ$ contains a nontrivial split $F$-stable torus, if 
and only if $\varphi_1$ occurs as a factor in $|{(\mathbf{Z}^\circ)}^F|$. In 
case~$\mathbf{C}$ is not connected, the entry corresponding to~$s$ also gives 
the order polynomial of $Z^\circ( \mathbf{C}^\circ )$, following the same 
conventions as for the order polynomial of~$\mathbf{Z}^\circ$. The tables
in~\cite{LL} thus easily allow to assign~$s$ to one of the three cases of
Remark~\ref{ExceptionalGroupsRemark}.

\subsubsection{} \label{ExplanationsE6}
Now let $\mathbf{G} := E_6( \mathbb{F} )_{\text{sc}}$, and let~$F$ be one 
of~$F'$ or~$F''$. Table~\ref{E6} contains a list of those semisimple elements 
$s \in G^*$ such that the $\mathbf{G}^*$-conjugacy class of~$s$ contains at 
least one $F$-stable element~$s'$ with $A_{\mathbf{G}^*}( s' )^F \neq 1$ and 
such that $\mathbf{C}^\circ_{\mathbf{G}^*}( s )$ has semisimple rank at 
least~$2$. 

Let us explain the notation used in Table~\ref{E6}. First of all, we
define $\varepsilon \in \{ \pm 1 \}$ by $\varepsilon = 1$ if $F = F'$, and
$\varepsilon = -1$ if $F = F''$. Now let~$s$ be such an element for which there
is an entry in the table. We then put $\mathbf{C} := C_{\mathbf{G}^*}(s)$ and
$\mathbf{Z} : = Z(\mathbf{C})$. We also write $C := \mathbf{C}^F$ and
$C^\circ := (\mathbf{C}^\circ)^F$; similarly, we put $Z := \mathbf{Z}^F$
and $Z^\circ := (\mathbf{Z}^\circ)^F$. The first column of Table~\ref{E6} just
numbers the cases, and the second column gives the label of the class type
of~$s$ according to~\cite{LL}. The third column gives the Dynkin type of
$[\mathbf{C}^\circ, \mathbf{C}^\circ]$, where~$A_2^3$ denotes three copies of
Type~$A_2$, etc. The fourth column describes
$[\mathbf{C}^\circ,\mathbf{C}^\circ]^F$. As we are only interested in the
unipotent characters of~$C^\circ$, which are insensitive to the center
of~$C^\circ$ and to isogeny, the information here is given in a Chevalley group
type of notation, just presenting the simple components of
$[\mathbf{C}^\circ,\mathbf{C}^\circ]^F$. Again, exponents denote the number of
copies of a specific group, and juxtaposition indicates direct products of
groups. Moreover, a notation such as $A_2( -q )$ stands for a twisted group of
type~$A_2$, defined over the field with~$q^2$ elements, i.e.\ a group with the
same unipotent characters as $\SU_3(q)$. The fifth column gives the order of
$C/C^\circ \cong {(\mathbf{C}/\mathbf{C}^\circ)}^F = A_{\mathbf{G}^*}(s)^F$.
The sixth column gives the order of $Z/Z^\circ \cong 
{(\mathbf{Z}/\mathbf{Z}^\circ)}^F$. The seventh column describes the torus
$\mathbf{Z}^\circ$ by its order polynomial (see the remarks in the introduction
to~\ref{TheExceptionalGroups}). The next column gives the conditions for the
existence of the elements in each row. Finally, the last column yields
information about the containment of $\mathbf{C}^\circ$ and $\mathbf{C}$ in
split $F$-stable Levi subgroups of~$\mathbf{G}^*$, where we use the symbols
introduced in Remark~\ref{ExceptionalGroupsRemark} to label the cases.

\begin{rem}
\label{E6Conclusion}
{\rm
As discussed in Remark~\ref{ExceptionalGroupsRemark}, the question about 
primitivity of the characters in $\mathcal{E}( G, [s] )$ for the classes in 
Table~\ref{E6} can be read off the ``Notes'' column of that table, provided the
entry is one of~$\ast$ or~\Checkmark. In the three cases where the entry is 
a~$\dagger$, we have to determine $A_{\mathbf{G}^*}(s)_\lambda^F$. For this we 
need to know the action of $A_{\mathbf{G}^*}(s)^F = C/C^\circ$ on the components 
of $[\mathbf{C}^\circ,\mathbf{C}^\circ]^F$. In each of the three cases, the 
Frobenius map~$F$ fixes the simple roots 
of~$[\mathbf{C}^\circ,\mathbf{C}^\circ]$, as indicated in the tables 
in~\cite{LL} (or as follows from the structure 
of~$[\mathbf{C}^\circ,\mathbf{C}^\circ]^F$).

\begin{landscape}
\begin{table}[h]
\caption{\label{E6} Some $F$-stable semisimple elements in $\mathbf{G}^* = 
E_6( \mathbb{F} )_{\text{ad}}$, $F \in \{ F', F'' \}$;
explanations in~\ref{ExplanationsE6}}
$$
\begin{array}{rcccccccccc}\\ \hline\hline
\text{No.} & 
\begin{array}{c} \text{Label} \\ \text{in \cite{LL}} \end{array} &
\begin{array}{c} \text{Dynkin} \\ \text{Type} \end{array} & 
[\mathbf{C}^\circ,\mathbf{C}^\circ]^F &
|C/C^\circ| &
|Z/Z^\circ| &
\mathbf{Z}^\circ &
\text{Condition} & 
\multicolumn{2}{c}{\begin{array}{cc} \multicolumn{2}{c}{\text{Notes}} \\ F' & F'' \end{array}}
\rule[-15pt]{0pt}{ 35pt} \\ \hline\hline
1 & [3,2,1] & A_2^3 & A_2(\varepsilon q)^3 & 3 & 3 & 1 & 
3 \mid q - \varepsilon & \ast & \ast \rule[- 3pt]{0pt}{15pt} \\
2 & [3,2,2] & & A_2( q^2 ) A_2(-\varepsilon q) & 1 & 3 & 1 & 
3 \mid q + \varepsilon & \ast & \ast \\
3 & [3,2,3] & & A_2(\varepsilon q^3) & 3 & 3 & 1 & 
3 \mid q - \varepsilon & \ast & \ast \\ \hline
4 & [13,2,1] & A_1^4 & A_1(q)^4 & 3 &  6 & 1 & 
6 \mid q - \varepsilon & \dagger & \ast \rule[- 3pt]{0pt}{15pt} \\
5 & [13,2,2] & & A_1(q)A_1(q^3) & 3 & 6 & 1 & 
6 \mid q - \varepsilon & \ast & \ast \\
6 & [13,2,3] & & A_1(q)^2A_1(q^2) & 1  & 6 & 1 & 
6 \mid q + \varepsilon & \text{\Checkmark} & \text{\Checkmark} \\ \hline
7 & [14,2,1] & D_4 & D_4(q) & 3 & 3 & 1 & 
3 \mid q - \varepsilon & \dagger & \ast \rule[- 3pt]{0pt}{15pt} \\
8 & [14,2,2] &              & {^3\!D}_4(q) & 3 & 3 & 1 & 
3 \mid q - \varepsilon & \ast & \ast \\ 
9 & [14,2,3] &              & {^2\!D}_4(q) & 1 & 3 & 1 & 
3 \mid q + \varepsilon & \text{\Checkmark} & \text{\Checkmark} \\ \hline
10 & [16,2,1] & A_1^3 & A_1(q)^3 & 3 & 3 & \Phi_1 & 
3 \mid q - \varepsilon & \text{\Checkmark} & \text{\Checkmark} \rule[- 3pt]{0pt}{15pt} \\
11 & [16,2,2] &       & A_1(q)^3 & 3 & 3 & \Phi_2 & 
3 \mid q - \varepsilon & \dagger & \ast \\
12 & [16,2,3] &       & A_1(q)A_1(q^2) & 1 & 3 & \Phi_1 & 
3 \mid q + \varepsilon & \text{\Checkmark} & \text{\Checkmark} \\
13 & [16,2,4] &       & A_1(q)A_1(q^2) & 1 & 3 & \Phi_2 & 
3 \mid q + \varepsilon & \text{\Checkmark} & \text{\Checkmark} \\
14 & [16,2,5] &       & A_1(q^3) & 3 & 3 & \Phi_2 & 
3 \mid q - \varepsilon & \ast & \ast \\
15 & [16,2,6] &       & A_1(q^3) & 3 & 3 & \Phi_1 & 
3 \mid q - \varepsilon & \text{\Checkmark} & \text{\Checkmark}  
\rule[- 3pt]{0pt}{13pt} \\ \hline\hline
\end{array}
$$
\end{table}
\end{landscape}

The action of $C/C^\circ$ on the
components of $[\mathbf{C}^\circ,\mathbf{C}^\circ]$ is described in the tables 
of~\cite{LL}.  We find that $C/C^\circ$ acts as follows: In Case~$4$ it fixes 
one component~$A_1(q)$ and acts as a three-cycle on the other components; in 
Case~$7$ it acts as the graph automorphism; in Case~$11$ it acts again by a 
three-cycle. For the action of the graph automorphisms on the unipotent 
characters on a group of type $D_4(q)$ see \cite[Theorem~$2.5$]{MalleExt}.

The classes with $|A_{\mathbf{G}^*}( s )^F| = 3$ not contained in Table~\ref{E6} 
can be treated as follows. For these classes, we always have $C_\lambda = C$.
Thus, as discussed in Remark~\ref{ExceptionalGroupsRemark}, either all elements
of $\mathcal{E}( G, [s] )$ are Harish-Chandra primitive, or all of them are
Harish-Chandra imprimitive. The latter occurs if and only if~$s$ is in 
Case~\Checkmark.
}
\end{rem}

\subsubsection{} \label{ExplanationsE7}
We now consider the case that $G = E_7(q)_{\text{sc}}$. In Table~\ref{E7I} we 
only display those semisimple class types, for which there is potentially a 
nontrivial action of~$A_{\mathbf{G}^*}(s)^F$ on the unipotent characters 
of~$[\mathbf{C}^\circ,\mathbf{C}^\circ]^F$, and which fall into Case~$\dagger$ 
of Remark~\ref{ExceptionalGroupsRemark}. The columns have the same meaning as 
in Table~\ref{E6}, except that we have omitted the column for the order 
of~$A_{\mathbf{G}^*}(s)^F$, since this order is always equal to~$2$. We have 
added a column headed ``Action'' which describes the action of 
$A_{\mathbf{G}^*}(s)^F = C/C^\circ$ on the components of 
$[\mathbf{C}^\circ,\mathbf{C}^\circ]^F$. These components are permuted by the 
action, and we give the cycle lengts of this permutation. If one of these 
components is fixed, so are its unipotent characters, except in the two 
instances where this is a component of type~$D_4(q)$. In this case we replace 
the~$1$ for the corresponding cycle length by ``$\text{g}$'' to indicate that 
$C/C^\circ$ induces the graph automorphism on this component.

\begin{table}[h]
\caption{\label{E7I} Some $F$-stable semisimple elements in $\mathbf{G}^* = 
E_7( \mathbb{F} )_{\text{ad}}$; explanations in~\ref{ExplanationsE7}}
$$
\begin{array}{rccccccc}\\ \hline\hline
\multicolumn{1}{c}{\text{No.}} & 
\text{Case} & 
\begin{array}{c} \text{Dynkin} \\ \text{Type} \end{array} & 
[\mathbf{C}^\circ,\mathbf{C}^\circ]^F &
\text{Action} &
|Z/Z^\circ| & 
|\mathbf{Z}^\circ|  &
\text{Condition} 
\rule[-15pt]{0pt}{ 35pt} \\ \hline\hline
 1 & [12,2,1] & A_2^3 & A_2( q )^3 & (1,2) & 6 & 1 & 6 \mid q - 1 \rule[- 3pt]{0pt}{15pt} \\ 
 2 & [17,2,1] & D_4A_1^2 & D_4( q )A_1(q)^2 & (\text{g},2) & 4 & 1 & 4 \mid q - 1 \\ 
 3 & [25,4,4] & A_3A_1^2 & A_3(-q)A_1(q)^2 & (1,2) & 2 & \Phi_2 & 2 \mid q - 1 \\ 
 4 & [27,3,3] & A_1^5 & A_1(q) A_1(q^2)^2 & (1,2) & 4 & \Phi_2 & 4 \mid q - 1 \\ 
 5 & [27,3,4] & & A_1(q) A_1(q^2)^2 & (1,1,1) & 4 & \Phi_2 & 4 \mid q + 1 \\ 
 6 & [27,3,9] & & A_1(q)^5 & (1,2,2) & 4 & \Phi_2 & 4 \mid q - 1 \\ 
 7 & [27,3,10] & & A_1(q) A_1(q^2)^2 & (1,2) & 4 & \Phi_2 & 4 \mid q + 1 \\ 
 8 & [30,2,6]  & A_2A_1^2 & A_2(-q)A_1(q)^2 & (1,2) & 2 & \Phi_2 & 2 \mid q - 1 \\ 
 9 & [33,4,2)  & A_2^2 & A_2(q)^2 & 2 & 2 & \Phi_3 & 2 \mid q - 1 \\ 
10 & [34,3,6]  & D_4 & D_4(q) & \text{g} & 2 & \Phi_2 & 2 \mid q - 1 \\ 
11 & [36,3,13] & A_1^4 & A_1(q^2)A_1(q)^2  & (1,2) & 2 & \Phi_4 & 2 \mid q - 1 \\
12 & [36,3,20] & & A_1(q^2)^2 & 2 & 2 & \Phi_2^2 & 2 \mid q - 1 \\ 
13 & [36,3,23] & & A_1(q)^4 & (2,2) & 2 & \Phi_2^2 & 2 \mid q - 1 \\ 
14 & [36,3,24] & & A_1(q^2)^2 & (1,1) & 2 & \Phi_2^2 & 2 \mid q - 1 \\ 
15 & [38,3,6]  & A_1^3 & A_1(q)^3 & (1,2) & 2 & \Phi_2^2 & 2 \mid q - 1 \\ 
16 & [41,5,6]  & A_1^2 & A_1(q)^2 & 2 & 2 & \Phi_2\Phi_4 & 2 \mid q - 1 \\
17 & [41,5,17]  & & A_1(q)^2 & 2 & 2 & \Phi_2\Phi_6 & 2 \mid q - 1 \\
18 & [41,5,18]  & & A_1(q)^2 & 2 & 2 & \Phi_2\Phi_4 & 2 \mid q - 1 \\
19 & [41,5,20]   & & A_1(q)^2 & 2 & 2 & \Phi_2^3 & 2 \mid q - 1  
\rule[- 3pt]{0pt}{13pt} \\ \hline\hline
\end{array}
$$
\end{table}
Again, the entries of Table~\ref{E7I} can be extracted from L{\"u}beck's 
tables~\cite{LL}. For the entries in the ``Action'' column one uses the 
information on the $F$-action on the set of simple roots 
of~$[\mathbf{C}^\circ,\mathbf{C}^\circ]$, together with the given information 
on the action of $C/C^\circ$ on this set.

\begin{rem}
\label{E7Conclusion}
{\rm
From the information contained in Table~\ref{E7I} we can decide which
irreducible characters of $\mathcal{E}( G, [s] )$ are Harish-Chandra imprimitive 
(see the explanations in Remark~\ref{ExceptionalGroupsRemark}). The classes 
of semisimple elements $s \in G^*$ with $|A_{\mathbf{G}^*}( s )^F| = 2$ not 
contained in Table~\ref{E7I} can be treated as in the analogous cases for 
$\mathbf{G} = E_6(\mathbb{F})_{\text{sc}}$.
}
\end{rem}

\section*{Acknowledgements}
We very much thank Frank L\"ubeck for innumerous conversations about this
project as well as for preparing tables of semisimple elements, not only
for the exceptional groups but also for some orthogonal groups of small ranks.
We also thank Cedric Bonnaf{\'e} for helpful conversations about his 
monograph~\cite{CeBo2}, and Meinolf Geck and Radha Kessar for pointing out 
reference~\cite{MalleExt}, respectively~\cite{En08}.
We would like to thank the DAAD (Deutscher Akademischer Austauschdienst)
for supporting a two week visit by the first author to Birmingham in May 2013.
This work was completed during a stay of the authors at the CIB at the EPFL 
within the semester ``Local representation theory and simple groups''. We
thank the CIB and the organizers of this semester for their hospitality and
for creating a unique working atmosphere, bringing all the experts in the field 
together. We finally thank Gunter Malle for a careful reading of our manuscript.


\begin{thebibliography}{101}
\bibitem{CeBo2}
{C.~Bonnaf{\'e.}}
\newblock{\em Sur les caract\`eres des groupes r\'eductifs finis a centre
non connexe: applications aux groupes sp\'eciaux lin\'eaires et unitaires.}
\newblock{Ast\'erisque 306, 2006.}

\bibitem{Bourbaki2}
N.~Bourbaki.
\newblock{\em Alg{\`e}bre, Chapitre IX.}
\newblock{Hermann, Paris, 1958.}

\bibitem{BrouMa}
\newblock{M.~Brou\'e and G.~Malle.}
\newblock{Th\'eor\`emes de Sylow g\'en\'eriques pour les groupes
r\'eductifs sur les corps finis.}
\newblock{Math.\ Ann.\ 292:241--262, 1992.}

\bibitem{CaEn}
{M.~Cabanes and M.~Enguehard.}
\newblock{\em Representation Theory of Finite Reductive Groups.}
\newblock{Cambridge University Press, Cambridge, 2004.}

\bibitem{cart}
{R.~W.~Carter.}
\newblock{\em Finite Groups of Lie Type:
              Conjugacy Classes and Complex Characters.}
\newblock{Wiley, New York, 1985.}

\bibitem{DeliLu} 
{P.~Deligne and G.~Lusztig.} 
\newblock{Representations of reductive groups over finite fields.}
\newblock{Ann.\ of Math.\ 103:103--161, 1976.}

\bibitem{DiMi2}
{F.~Digne and J.~Michel.}
\newblock{\em Representations of Finite Groups of Lie Type.}
\newblock{London Math.\ Soc.\
                 Students Texts 21. Cambridge University Press, 1991.}

\bibitem{DM2}
D.~\v{Z}.~Djokovi\'c and J. Malzan.
\newblock{Imprimitive, irreducible complex characters of the alternating group.}
\newblock{Can. J.~Math. 28:1199--1204, 1976.}

\bibitem{En08}   
{M.~Enguehard.} 
\newblock{Vers une d{\'e}composition de Jordan des blocs des groupes 
r{\'e}ductifs finis.}
\newblock{J.~Algebra 319:1035--1115, 2008.}

\bibitem{fosri1}
{P.~Fong and B.~Srinivasan.}
\newblock{The blocks of finite general linear and unitary groups.}
\newblock{Invent.\ Math.\ 69:109--153, 1982.}

\bibitem{fosri2}
{P.~Fong and B.~Srinivasan.}
\newblock{The blocks of finite classical groups.}
\newblock{J.~reine angew.\ Math.\ 396:122--191, 1989.}

\bibitem{GAP4}
The GAP~Group.
\newblock{\emph{{GAP -- Groups, Algorithms, and Programming, 
Version 4.8.5}};}
\newblock{2016,
\url{http://www-gap.dcs.st-and.ac.uk/~gap}.}

\bibitem{chevie}
{M.~Geck, G.~Hiss, F.~L\"ubeck, G.~Malle, and G.~Pfeiffer}.
\newblock{{\sf CHEVIE} --- A system for computing and processing
generic character tables.}
\newblock{AAECC 7:175--210, 1996.}

\bibitem{GP} {M.~Geck and G.~Pfeiffer}. {\em Characters of Finite Coxeter
             Groups and Iwahori-Hecke Algebras}. Oxford University Press,
             Oxford, 2000.

\bibitem{HiHuMa}
{G.~Hiss, W.~J.~Husen, and K.~Magaard.}
\newblock{Imprimitive irreducible modules for finite quasisimple groups.}
{Mem.\ Amer.\ Math.\ Soc.} {\bf 234}, No.~1104, 2015.

\bibitem{HowLeh2}
{R.~B.~Howlett and G.~I.~Lehrer}.
\newblock{Representations of generic algebras and finite groups of Lie type.}
\newblock{Trans. Amer. Math. Soc. 280:753--779, 1983.}

\bibitem{jake}
{G.~D.~James and A.~Kerber.}
\newblock{\em The Representation Theory of the Symmetric Group.}
\newblock{Encyclopedia Math.~16, Addison-Wesley Publishing Company 1981.}

\bibitem{LL}
{F.~L{\"u}beck.}
Centralizers and numbers of semisimple classes in exceptional groups of Lie type.
\url{http://www.math.rwth-aachen.de/~Frank.Luebeck/chev/CentSSClasses/}.

\bibitem{lusz2}  
{G.~Lusztig.}
\newblock{Irreducible representations of finite classical groups.}
\newblock{Invent.\ Math.\ 43:125--175, 1977.}

\bibitem{luszbuch}
{G.~Lusztig.} 
\newblock{\em Characters of Reductive Groups over a Finite Field.}
\newblock{Ann.\ Math.\ Studies  107, Princeton University Press, 1984.}

\bibitem{Lu}
{G.~Lusztig.}
\newblock{On the representations of reductive groups with disconnected centre.}
\newblock{Ast\'erisque 168:157--166, 1988.}

\bibitem{MalleExt}
{G.~Malle.}
Extensions of unipotent characters and the inductive McKay condition. 
J.~Algebra 320:2963--2980, 2008.

\bibitem{NeNoe}
{D.~Nett and F.~Noeske.}
\newblock{The imprimitive faithful complex characters of the Schur covers
of the symmetric and alternating groups.}
\newblock{J.~Group Theory 14:413--435, 2011.}

\bibitem{Taylor}
{D.~E.~Taylor.}
\newblock{\em The Geometry of the Classical Groups.}
\newblock{Heldermann Verlag, Berlin, 1992.}

\end{thebibliography}
\end{document}